\documentclass[11pt,letterpaper]{amsart}

\usepackage{amsmath}

\usepackage{amssymb}

\usepackage[pdftex]{graphicx}

\renewcommand{\le}{\leqslant}
\renewcommand{\ge}{\geqslant}

\usepackage{versions}
\usepackage{hyperref}
\hypersetup{
	colorlinks   = true, 
	urlcolor     = blue, 
	linkcolor    = blue, 
	citecolor   = red 
}
\usepackage{cancel}
\usepackage{yhmath}
\usepackage{enumitem}
\usepackage{stackengine}

\usepackage{lineno}

\newtheorem{theorem}{Theorem}[section]

\newtheorem{lemma}[theorem]{Lemma}

\newtheorem{conjecture}[theorem]{Conjecture}
\newtheorem{corollary}[theorem]{Corollary}

\theoremstyle{definition}
\newtheorem{definition}[theorem]{Definition}

\theoremstyle{remark}

\numberwithin{equation}{section}

\newcommand{\inte}{\mathop{\rm int}}

\newcommand{\diam}{\mathop{\rm diam}}

\newcommand{\dist}{\mathop{\rm dist}}

\newcommand*{\where}{\ \ifnum\currentgrouptype=16 \middle\fi|\ }
\newcommand{\hc}{{\rm HC}} 
\newcommand{\HC}{{\rm HC}} 

\renewcommand{\epsilon}{\varepsilon}
\renewcommand{\phi}{\varphi}
\renewcommand{\kappa}{\varkappa}
\renewcommand{\theta}{\vartheta}

\def\R{{\mathbb R}}

\long\def\forget#1\forgotten{}
\begin{document}

\title{Boxing inequalities in Banach spaces and Riemannian manifolds}

\author{Sergey Avvakumov, Alexander Nabutovsky}
\date{}

\begin{abstract}
We prove the following result: For each closed $n$-dimensional manifold $M^n$ in a finite or infinite-dimensional Banach space $B$,
and each positive real $m\leq n$
there exists a pseudomanifold $W^{n+1}\subset B$ such that $\partial W^{n+1}=M^n$ and $\HC_m(W^{n+1})\leq c(m)\HC_m(M^n)$. Here $\HC_m(X)$ denotes the $m$-dimensional Hausdorff content, i.e the infimum of $\Sigma_i r_i^m$, where the infimum is taken over all coverings of $X$
by a finite collection of metric balls, and $r_i$ denote the radii of these balls.

In the classical
case, when $B=\mathbb{R}^{n+1}$ and $M^n=\partial\Omega$ for a bounded domain $\Omega\subset\mathbb{R}^{n+1}$, this result implies that 
for all $m\in (0,n]$ \
$\HC_m(\Omega)\leq c(m)\HC_m(\partial \Omega)$. This inequality seems to be new despite being well-known and widely used in the case, when $m=n$ (Gustin's boxing inequality, [G]).

The result is a corollary of the following more general theorem: For each compact subset $X$ in a 
Banach space $B$ and a real number $m>0$  such that $\HC_m(X)\not= 0$, 
there exists a homotopy of $X$ into a $(\lceil m\rceil-1)$-dimensional simplicial complex in $B$ such that the lengths of all trajectories of the homotopy do not exceed $c_1(m)\HC_m^{\frac{1}{m}}(X)$, and
$\HC_m$ of the image of the homotopy does not exceed $c_2(m)\HC_m(X)$.

A similar theorem can also be proven in the case when $B$ is a metric space (for example, a Riemannian manifold) such that for some $\Lambda\geq 1$ each metric ball $\beta$ 
is contractible in a concentric ball of radius
$\Lambda\cdot$radius$(\beta)$. 
\forget
where the constants $c_1, c_2$ depend on $m$ and $\Lambda$,
and the conclusion is somewhat weaker: $c_1(m, \Lambda)\HC_m(X)^{1\over m}$ majorizes only the distances between each point  $x\in X$ and the endpoint of the trajectory of $x$ under the homotopy. 
If each metric ball $\beta$ can be contracted to a point within $\Lambda \beta$
by a homotopy where
the lengths of all trajectories do not exceed $\mu\cdot$radius$(\beta)$, then the length of the trajectories can also be bounded by $2\mu c(m,\Lambda)$. 
If the contractibility property holds
only for metric balls of radius $\leq r_0$ for some positive $r_0$, then the assertion of the theorem holds for all compact $X$ such that $\HC_m(X)\leq c_3(m,\Lambda ) r_0^m$. 
The theorem can be generalized to all $(\Lambda,\mu, r_0)$-linearly contractible metric spaces $M$ provided that $X$ is the image of a continuous map of a finite-dimensional polyhedron to $M$.
\forgotten
In particular, our theorem applies to all compact sets $X$ with a controllably small $\HC_m$ in Riemannian manifolds $M^n$
with either the injectivity radius bounded below by a positive number, or the sectional curvature bounded below, the volume bounded below by a positive number, and the diameter bounded above. 
\end{abstract}

\maketitle 

\section{Introduction}

\subsection{Hausdorff content}

Recall that for each nonnegative real $m$ \par\noindent 
$m$-dimensional Hausdorff content $\HC_m(A)$ of a subset $A$ in a metric space $X$ is defined
as the infimum of $\Sigma_i r_i^m$ over all coverings of $A$ by collections of metric balls $B_i$ with radii $r_i$. It is obvious that the result does not depend on whether we assume that the balls are open or closed. Therefore, if $A$ is compact, it is sufficient to consider only finite coverings. It is clear that if $A$ is bounded, then $\HC_m(A)$ is finite for any $m$ (unlike the $m$-dimensional Hausdorff measure, which is infinite, when $m<dim_{Haus}(A)$). In fact, if $A$ is bounded, it can be covered by a single ball.
The Hausdorff contents enjoy several convenient and easy-to-prove
properties: 
\par\noindent
-Monotonicity:  $\HC_m(A)\leq \HC_m(B)$, if $A\subset B$;
\par\noindent
-Semi-additivity: $\HC_m(\cup_i A_i)\leq \sum_i \HC_m(A_i)$;
\par\noindent
-Good behavior under Lipschitz maps: if $f:A\longrightarrow B$ is $L$-Lipschitz, then 
\par\noindent
$\HC_m(f(A))\leq L^m \HC_m(A)$;
\par\noindent
- If $m\in (0,1]$, and $A$ is a bounded and connected subset of a Banach space $X$, then $\HC_m^{\frac{1}{m}}(A)= rad(A),$ the smallest radius of a closed ball in $X$ containing $A$. In fact, if $m\leq 1$ , it is advantageous to merge intersecting metric balls into a larger metric ball. (Since $A$ is connected, 
\par\noindent
- If $A$ is fixed, then $\HC_m^{\frac{1}{m}}(A)$ is a decreasing function of $m$. To see this, observe that if $m\leq k$, then 
$\sum_i r_i^k=\sum_i (r_i^m)^{\frac{k}{m}}\leq (\sum_i r_i^m)^{\frac{k}{m}}.$
\par\noindent
- Comparing the definition of Hausdorff content with Hausdorff measure, we see that if $A$ is a $m$-dimensional rectifiable subset of a Riemannian manifold $M^n$ then 1) $\HC_m(A)\leq\omega_m^{-1} H_m(A)$, where $H_m(A)$ denotes the $m$-dimensional Hausdorff measure of $A$, and $\omega_m$ denotes the volume of the Euclidean unit ball; and 2) If $m=n$ and $A\subset \R^n$, then $\HC_n(A)=\omega_n^{-1} vol_n(A)$, where $vol_n(A)$ is the $n$-dimensional Hausdorff measure, or equivalently, the Lebesgue measure of $A$.
Here, the factor $\omega_m$ appears because $H_m(A)$ is scaled by the factor $\omega_m$, and $HC_m(A)$ is not. (The scaling is chosen so that the $m$-dimensional Hausdorff measure of a Euclidean $m$-ball of radius $1$ is equal to $\omega_m$, while its $m$-dimensional Hausdorff content is equal to $1$.)
\par\noindent
- However, for every $m$ Hausdorff content $\HC_m$ is not additive. For example, if $A$ is a connected bounded subset of a Banach space $X$, $\HC_1(A)$ coincides with the radius $rad(A)$ of $A$.
It is easy to use this observation to construct
examples of disjoint connected sets $A$ and $B$ with $\HC_1(A\bigcup B)< \HC_1(A)+\HC_1(B)$. 

A classical fact about Hausdorff content is the following theorem proven by W. Gustin ([G]) called the boxing inequality:

\begin{theorem} Let $\Omega$ be a bounded open domain in $\mathbb{R}^{n+1}$ with a smooth boundary, $\partial\Omega$. Then $\HC_n(\Omega)\leq c(n)\HC_n(\partial\Omega)$ for some constant $c(n)$.
\end{theorem}

This implies $\HC_n(\Omega)\leq C(n)vol_n(\partial \Omega)$, which is the form of the inequality usually quoted in the literature. As $\HC_n(\Omega)\geq \HC_{n+1}^{\frac{n}{n+1}}(\Omega)=\omega_{n+1}^{-\frac{n}{n+1}}vol_{n+1}(\Omega)$, this inequality is strictly stronger than the non-sharp isoperimetric inequality (that is, the isoperimetric inequality without the exact value of the constant).

\subsection{Main results}
One of our results is an extension of the previous theorem for all $m\leq n$:

\begin{theorem} Let $\Omega$ be a bounded open domain in $\mathbb{R}^{n+1}$ with boundary $\partial\Omega$. Then for each $m\in (0,n]$ $\HC_m(\Omega)\leq c(m)\HC_m(\partial\Omega)$ for some constant $c(m)$ that depends only on the parameter $m$ but not on the dimension of $\Omega$.
\end{theorem}

More generally, the following theorem is true:

\begin{theorem}
Let $M^n$ be a closed submanifold in a Banach space $B$. Then there
is a pseudomanifold $W^{n+1}\subset B$ such that :
\par\noindent
(1) $M^n$ is the boundary of $W^{n+1}$;
\par\noindent
(2) $\HC_m(W^{n+1})\leq c(m)\HC_m(M^n)$;
\par\noindent
(3) For each $w\in W^{n+1}$ there exists a point $x\in M^n$ such that $\Vert w-x\Vert\leq c(m)\HC_m^{\frac{1}{m}}(M^n)$.
\end{theorem}

The last assertion is the upper bound for Gromov's filling radius in terms of $\HC_m(M^n)$. To see that Theorem 1.3 implies Theorem 1.2 apply it to $M^n=\partial\Omega\subset B=\mathbb{R}^{n+1}$. Any $W^{n+1}$
bounded by $\partial\Omega$ must contain $\Omega$.
\par
Recall, that $n$-dimensional {\it pseudomanifold} is a triangulated topological space $X=\vert K\vert$,
such that each $(n-1)$-dimensional simplex in $K$ is in the boundary of exactly one or
two $n$-dimensional simplices; $X$ is the union of all $n$-dimensional simplices of $K$;
and each pair of $n$-dimensional simplices $\sigma$, $\sigma'$ can be connected via
a finite sequence of $n$-dimensional simplices $\sigma_1=\sigma,\ldots, \sigma_{N-1}, \sigma_N=\sigma'$, so that for each $i$ $\sigma_i$ and $\sigma_{i+1}$ share a common $(n-1)$-dimensional face. Informally speaking, the second condition means that $X$ is $n$-dimensional near each of its points; the third condition precludes nearly disconnected
examples like cones over disconnected manifolds. Each smooth manifold or PL manifold is a pseudomanifold, but the cone on a connected manifold $M$ is also a pseudomanifold
and $M$ will be its boundary. So, each manifold is trivially a boundary of a pseudomanifold. On the other hand, not every manifold is a boundary of a manifold,
and this is why we need to use pseudo-manifolds in this theorem.

\par
The following definition of the filling of metric spaces first appeared in the paper by Y. Liokumovich, B. Lishak,
the second named author, and R. Rotman in [LLNR]:

\par\noindent
\begin{definition} Let $X\subset B$ a compact set in a Banach space $B$. For each real $m>0$ a $m$-filling of $X$ consists of a
finite $(\lceil m\rceil -1)$-dimensional simplicial complex $K\subset B$, a Lipschitz map $\phi:X\longrightarrow K\subset B$, and a Lipschitz homotopy $H$ between the identity map of $X$ and $\phi$, where both maps are considered maps from $X$ to $B$.
\end{definition}

The main result of this paper is the following theorem.

\begin{theorem}{\bf (Boxing inequality for compact sets).} For each positive $m$ there exist constants $c_1(m), c_2(m)$ with the following property. Let $X$ be a compact subset of a Banach space $B$ such that $\HC_m(X)\not =0$. Then there exists
a $m$-filling $(K,\phi, H)$ of $X$ such that:
\medskip
\par\noindent
(1) For each $x\in X$ $$\Vert \phi(x)-x\Vert_B\leq c_1(m)\HC_m^{\frac{1}{m}}(X)$$ and, moreover, the length of the trajectory $H(\{x\}\times [0,1])$ of $x$ is bounded by $c_1(m)\HC_m^{\frac{1}{m}}(X)$;
\medskip
\par\noindent
(2) $\HC_m(H(X\times[0,1]))\leq c_2(m)\HC_m(X)$.
\medskip
\par\noindent
If $\HC_m(X)=0$, then for each $\delta>0$  there exists an $m$-filling of $X$ such that for each $x\in X$ $\Vert \phi(x)-x\Vert_B\leq \delta$, and $\HC_m(H(X\times [0,1]))\leq\delta.$

\end{theorem}

Each trajectory of the homotopy $H$ that will be constructed in the proof of this theorem will consist of finitely many straight line segments; its length will be the sum
of the distances between the end points of these segments.

This result strengthens the main result of [LLNR] (Theorem 3.1 in [LLNR]). The difference between our Theorem 1.5 and Theorem 3.1 in [LLNR] 
is that the result in [LLNR] asserted only the inequality
$\HC_{m+1}(H(X\times[0,1]))\leq c_2(m)\HC_m^{\frac{m+1}{m}}(X)$ instead of Theorem 1.5 (2). We can rewrite the inequality in [LLNR] as
$\HC_{m+1}^{\frac{1}{m+1}}(H(X\times[0,1]))\leq c_2(m)\HC_m^{\frac{1}{m}}(X)$. This inequality is strictly weaker than Theorem 1.5(2), because $\HC_t^{\frac{1}{t}}(X)$ is a decreasing function of $t$. 

\medskip\noindent

We believe that the following stronger version of Theorem 1.5 is true:
\begin{conjecture}
\par\noindent
(1) One can choose the constant $c_1(m)$ in Theorem 1.5 (1) so that it does not depend on $m$.
\par\noindent
(2) There exists an absolute constant $c$ that does not depend on $m$ so that Theorem 1.5(2) will be true with $c_2(m)=c^m$.
\end{conjecture}
Similarly to an analogous argument in [LLNR], the validity of Conjecture 1.6 (1) for $B=L^\infty(M^n)$ immediately implies Gromov's systolic inequality 
with the constant that depends on the dimension $n$ in the optimal possible way: $sys_1(M^n)\leq const\sqrt{n}\ vol^{\frac{1}{n}}(M^n)$ for each essential Riemannian
manifold $M^n$. (Here $sys_1(M^n)$ denotes the length of the shortest non-contractible closed curve on $M^n$.)

The boxing inequality for compact sets holds for a more general class of geodesic spaces $B$ than Banach spaces. Recall that a metric space is called a geodesic space if
each pair of points $x$,$y$ can be connected by a path so that its length is equal to the distance between $x$ and $y$.

\begin{definition} Let $B$ be a geodesic space,
$\mu, \Lambda\geq 1$ be real numbers, $r_0$  a positive real number or $\infty$.
\par\noindent
We say that $B$ is $(\Lambda,\infty,r_0)$-linearly
contractible if for each $r< r_0$ each metric ball
$\beta(x,r)$ of radius $r$ centered at a point $x\in B$ is contractible to a point within the
concentric ball $\beta(x,\Lambda\cdot r)$
of radius $\Lambda\cdot r$. We call $B$ $(\Lambda,\mu, r_0)$-linearly contractible if it is $(\Lambda, \infty,r_0)$-linearly contractible, and also, for each metric ball $\beta(x,r)$ of radius $r\leq r_0$ one can find a contracting homotopy $H$ with the image in $\beta(x,\Lambda\cdot r)$ such that the length of the trajectory
of each point under $H$ does not exceed $\mu\Lambda r$. 
\par\noindent
\end{definition}

Observe that all Banach spaces
are $(1,1,\infty)$-linearly contractible, and the same is true for all complete CAT(0) spaces.

An interesting example of locally $(\Lambda,\mu, r_0)$-contractible spaces is provided by the class of all closed Riemannian manifolds of fixed dimension $n$, sectional curvature bounded below by a real constant $x$, volume bounded below by a constant $v>0$,
and diameter $\leq d$. K. Grove and P. Petersen proved in [KP] that all such manifolds
are $(\Lambda,\mu,r_0)$ linearly contractible
for some $\Lambda$, $\mu$, $r_0$ given by explicit functions of $n, x, v, d$.
Also, for a positive $i$, consider the class all closed Riemannian manifolds with 
injectivity radius greater than or equal to $i$.
All manifolds in this class are obviously  $(1,1,inj)$-linearly contractible.
As it is not natural to expect that a general metric space would contain embedded polyhedra, we are going to slightly modify the definition of $m$-filling: Now, it is not assumed that $(\lceil m\rceil-1)$-dimensional complex $K$ is embedded in $B$. Instead, we assume that $K$ is  mapped to $B$ by means of a continuous map $j:K\longrightarrow B$, and $H$ is the homotopy between the inclusion of $X$ in $B$ and $j\circ\phi$, where, as before, $\phi$ is a continuous map $X\longrightarrow K$. So, now a $m$-filling of $X$ is a quadruple $(K, j, \phi, H)$.

In the last section, we prove that

\begin{theorem}{\bf (The boxing inequality for compact sets in linearly contractible metric spaces)}

\medskip\noindent
    {\bf A.} For each positive $m$, $r_0\in (0, \infty]$, and $\Lambda,\mu \geq 1$ 
    there exist constants $\epsilon(m), \tilde{c_1}(m), \tilde{c_2}(m), a_1(m), a_2(m), a_3(m)$
    with the following property. 
    Let $M$ be a $(\Lambda,\infty ,r_0)$-linearly contractible Riemannian manifold, and $X$ be a compact subset of $M$ such that $\HC_m(X)\not =0$. If $r_0\not =\infty$, assume also that $\HC_m(X)\leq \epsilon(m)\Lambda^{-a_3(m)}r_0^m$. Then there exists
a $m$-filling $(K,\phi, j, H)$ of $X$ such that:
\medskip
\par\noindent
(1) For each $x\in X$ $$dist_B(\phi(x), x)\leq \tilde{c_1}(m)\Lambda^{a_1(m)}
\HC_m^{\frac{1}{m}}(X).$$ Moreover,  if $M$ is $(\Lambda,\mu,r_0)$-linearly contractible for $\mu\not=\infty$, then the length of the trajectory $H(\{x\}\times [0,1])$ of $x$ is bounded by $2\mu\tilde{c_1}(m)\Lambda^{a_1(m)}
\HC_m^{\frac{1}{m}}(X);$
\medskip
\par\noindent
(2) $\HC_m(H(X\times[0,1]))\leq \tilde{c_2}(m)\Lambda^{a_2(m)}
\HC_m(X)$.
\medskip
\par\noindent
If $\HC_m(X)=0$, then for each $\delta>0$  there exists an $m$-filling of $X$ such that for each $x\in X$ $dist_B(\phi(x),x)\leq \delta$, and $\HC_m(H(X\times [0,1]))\leq\delta.$

\medskip\noindent
{\bf B.} This assertion remains true, when $M$ is an arbitrary $(\Lambda, \mu, r_0)$-linearly contractible metric space and $X$
is 
a finite-dimensional polyhedron in $M$. 

\end{theorem}

\par\noindent
{\bf Remark.} Part A of the above theorem easily follows from Part B. Indeed, assume that $M$ is a Riemannian manifold.
As each (smooth) Riemannian metric can be approximated by
an analytic metric, it is sufficient to prove this theorem in the case of an analytic Riemannian manifold. The distance function on an analytic Riemannian manifold is subanalytic, and the metric balls will be triangulable. Now, given a compact subset set of $M$ $X$ and a finite collection of metric balls $\beta_i$ with radii $r_i$ covering $X$ such that
$\HC_m(X)\geq \Sigma_i r_i+\epsilon$, we can replace $X$ by the polyhedron$\bigcup_i\beta_i$. Now, the result of part A will follow from part B. Nevertheless, in the last section we will give
different proofs of parts A and B.

As an immediate corollary of this theorem and the theorem by Grove and Petersen,
we obtain the following theorem.

\begin{theorem}{\bf (The boxing inequality for subsets of closed Riemannian manifolds with sectional curvature bounded below).}
Let $M^n$ be a closed Riemannian manifold such that its sectional curvature is bounded below by a real number $x$, the volume is bounded below by $v>0$, and the diameter does not exceed $D$. Let $m$ 
be a positive number, and $X$ be a compact subset of $M^n$
such that $\HC_m(X)\not =0$. Then there exists positive functions $c_1(n,x,v,d), c_2(n,x,v,D), c_3(n,x,v,d)$ (that can be presented by explicit formulas) such that if $\HC_m(X)<c_3(n,x,v,D)$, then there exists
a $m$-filling $(K,\phi, H)$ of $X$ in $M^n$ with $K\subset M$ such that:
\medskip
\par\noindent
(1) For each $x\in X$ $$dist_{M^n}(x,\phi(x))\leq c_1(n,x,v,D)\HC_m^{\frac{1}{m}}(X)$$ and, moreover, the length of the trajectory $H(\{x\}\times [0,1])$ of $x$ is bounded by $c_1(n,x,v,D)\HC_m^{\frac{1}{m}}(X)$;
\medskip
\par\noindent
(2) $\HC_m(H(X\times[0,1]))\leq c_2(n,x,v,D)\HC_m(X)$.
\medskip
\par\noindent
If $\HC_m(X)=0$, then for each $\delta>0$  there exists an $m$-filling of $X$ such that for each $x\in X$ $dist_{M^n}( \phi(x),x)\leq \delta$, and $\HC_m(H(X\times [0,1]))\leq\delta.$
\end{theorem}

Similarly,

\begin{theorem}{\bf (The boxing inequality for subsets of closed Riemannian manifolds with injectivity radius bounded below).}
Let $M^n$ be a complete Riemannian manifold such that its injectivity radius  is bounded below by a positive number $inj$, and
$m\leq n-1$ 
be a positive real number. 

Then there exist positive constants $\epsilon(m), c_1(m), c_2(m)$ with the following property. Let $M^n$ be a complete Riemannian manifold such that its injectivity radius is bounded below by a positive number $inj$, and $X$ be a compact subset of $M^n$
such that $\HC_m(X)<\epsilon(m) inj^m$. If $\HC_m(X)$ is non-zero, then
there exists
a $m$-filling $(K,\phi, H)$ of $X$ in $M^n$ with $K\subset M$ such that:
\medskip
\par\noindent
(1) For each $x\in X$ $$dist_{M^n}(x,\phi(x))\leq c_1(m)\HC_m^{\frac{1}{m}}(X)$$ and, moreover, the length of the trajectory $H(\{x\}\times [0,1])$ of $x$ is bounded by $c_1(m)\HC_m^{\frac{1}{m}}(X)$;
\medskip
\par\noindent
(2) $\HC_m(H(X\times[0,1]))\leq c_2(m)\HC_m(X)$.
\medskip
\par\noindent
If $\HC_m(X)=0$, then for each $\delta>0$  there exists an $m$-filling of $X$ such that for each $x\in X$ $dist_{M^n}( \phi(x),x)\leq \delta$, and $\HC_m(H(X\times [0,1]))\leq\delta.$
\end{theorem}

\begin{corollary}{\bf (Meta-corollary.)}
Assume that $M$, $X$ are as in Theorem 1.9, or Theorem 1.10. Assume that $\HC_m(X)$ is ``controllably small", that is, it satisfies the assumptions of the corresponding theorem. Then $X$ can be homotoped to a subset of a $(\lceil m\rceil-1)$-dimensional subpolyhedron of $M$.
\end{corollary}

\forget
If we approximate the metric on $M$ by an analytic metric, the distance function will
be subanalytic, and our proof will produce a subanalytic homotopy $H$ that will have
a polyhedral image in the case when $X$ is a polyhedron. If $X$ is a cycle, then it will be the boundary of the image of $H$ regarded as a singular chain. Therefore, $X$
is null homologous. Also, recall that $\HC_m(X)\leq c(m)vol_m(X)$, if $X$ is a $m$-dimensional polyhedron. Combining this facts together we obtain another corollary"

\begin{corollary} There exist positive functions $c(n,x,v,d)$ and $const(n)$ with the following
property. Let $M$ be a closed $n$-dimensional Riemannian manifold with sectional curvature $\geq x$,
volume $\geq v>0$, and diameter $\leq d$. Then the $m$-dimensional volume
of each $m$-dimensional polyhedral cycle $X$ in $M$ representing a non-trivial $m$-dimensional 
homology class of $M$ is not less than $c(n,x,v,d)$. The same will be true 
for $\HC_m(X)$. If $M$ is a Riemannian manifold with injectivity radius $\geq i>0$,
then the volume of each $m$-dimensional polyhedral cycle $X$ in $M$ representing a non-zero class in $H_m(X)$ is at least $const(n)i^m$. The same is true for $\HC_m(X)$.
    
\end{corollary}

\forgotten
\forget
Gustin's boxing inequality was applied in analysis via 
Choquet integrals (see [PS], [CS], and the references there for some recent results in this direction).
Choquet integral $\int_\Omega fd\HC_m$ of a nonnegative function $f$ on a domain $\Omega\subset\R^n$ is defined by the formula
$\int_\Omega fd\HC_m=\int_0^\infty \HC_m(\{x\in\Omega\vert f(x)\geq t\}) dt$. One can also define {\it dyadic Hausdorff content}
$\HC_m^d$ defined for each set $A$ as the infimum of $\Sigma_i (\frac{l_i}{2})^m$, where $l_i$ denote side lengths of dyadic cubes in a covering of $A$ by dyadic cubes. The infimum is taken over all coverings of $A$ by dyadic cubes. Recall that dyadic cubes in $\R^n$ are cubes with
sides parallel to the coordinate axes, center of the form $2^{m-1}(2N+1)$ and side length $2^m$, where $m$ can be an arbitrary integer number (of any sign), and $N$ is an arbitrary $n$-dimensional vector with integer entries. It is well known and easy to see that $\HC_m(A)\leq \HC_m^d(A)\leq c(n)\HC_m(A)$ (cf. Lemma 3.3 below). Therefore, $\int_\Omega f d\HC^d_m$ differs from $\int_\Omega f d\HC_m$ by at most a constant factor $c(n)$
depending only on the dimension of the ambient space. If one is interested in inequalities with an arbitrary constant of this form,
then one can interchange Choquet integrals with respect to $d\HC_m$ and $d\HC^d_m$. The reason why one might want to do that is rather surprising: Choquet integrals with respect to $d\HC^d_m$ are countably subadditive, but Choquet integrals with respect to
$d\HC_m$ are not. However, in the present paper we are interested in constants that depend only on $m$ but not on $n$. Therefore,
for us $\HC^d_m$ and $\HC_m$ are not roughly the same. In this connection, we would like to state the following theorem:

\begin{theorem}
Let $B=l^n_\infty$. Then Theorems 1.2, 1.3, 1.5 remain true for $\HC^d_m$ instead of $\HC_m$.
For example, if $\Omega$ is a bounded domain in $\R^n$, then for each $m\leq n-1$
$\HC^d_m(\Omega)\leq c(m)\HC^d_m(\partial\Omega)$.
\end{theorem}

Theorems 1.2, 1.3 are almost immediate corollaries of Theorem 1.5. Dyadic cubes are metric balls in $l^\infty$ metric. Moreover, the part of our proof of Theorem 1.5 presented in Section 2 is similar to the argument in [LLNR], and we observed in [LLNR] that this argument remains valid when one considers only metric balls with centers in a prescribed subset $S$ of $B$ to define (a new version of) $\HC_m$. Moreover, the new part of our proof of Theorem 1.5 
presented in Section 3 immediately works for $\HC^d_m$.
Nevertheless, our proof
of Theorem 1.5 is not applicable as is to establish Theorem 1.6.
One complication is that when $Q$ is a dyadic cube, none of the larger cubes with the same center as $Q$ are dyadic. As a result, the argument presented in Section 2 needs to be amended.
We will present the details in Section 4.
\forgotten

\subsection{Some historical remarks and connections with analysis} 
Theorem 1.1 was proven by Gustin answering a question posed by W. Fleming. The main result of [G] was a (qualitatively weaker) inequality
$\HC_n(\Omega)\leq c(n)vol(\partial\Omega)$, and this is the form that is usually quoted as the boxing inequality. Soon after, a simpler proof of this inequality was found by H. Federer in [F], where it was used to prove a foundational result in geometric measure theory. However, the boxing inequality has been having numerous applications in analysis after V. Maz'ya (cf. [MaS]) and D.R. Adams (cf. [Ad]) discovered its functional forms, for example:

$$ \int_{\mathbb{R}^n}\vert\phi\vert d\HC_m\leq c(n,m)\Vert D^{n-m}\phi\Vert_{L^1(\mathbb{R}^n)},\ \ \ \ \ \eqno(1.1)$$

\noindent where $m<n$ is any integer, $\phi$ is an arbitrary function from $C_0^\infty(\mathbb{R}^n)$, and
$\int_{\mathbb{R}^n}\vert\phi\vert d\HC_m$ is the Choquet integral defined as $\int_0^\infty \HC_m(\{\vert \phi\vert>t\})dt.$
This inequality can also be extended for non-integer $m$; see [PS].
Since for $m=n-1$ the boxing inequality is finer than the non-sharp isoperimetric inequality, this is a qualitatively better inequality than the corresponding Sobolev inequality. When $m$ is lower, this inequality can be used to control the integrals of functions of interest over sets of lower dimensions. Also, note the fact (also proven by Maz'ya and Adams) that Hausdorff contents differ by not more than a constant factor from the corresponding capacities:
$\HC_m(A)\sim Cap_{n-m,1}(A)$.
Recall that for an integer $d$ and a compact $A$ $Cap_{d,1}(A)$ is defined as $\inf_{\phi\in C_c^\infty(\mathbb{R}^n); \phi\geq 0; \phi\vert_A=1}\Vert D^d\phi\Vert_{L^1(\mathbb{R}^n)}$. If $A$ is open, $Cap_{d,1}(A)$ is defined as the supremum of $Cap_{d,1}(K)$ over all compact $K\subset A$.
Therefore, for $m<n$

$$ \int_{\mathbb{R}^n}\vert\phi\vert d Cap_{n-m,1}\leq C(n,m)\Vert D^{n-m}\phi\Vert_{L^1(\mathbb{R}^n)},\ \ \ \ \ \eqno(1.2)$$

\noindent ([PS]).
As it was noted in [PS], capacities arise in studies of PDE when one needs to quantify the size of the set on which certain fine properties hold.
Thus, the inequalities (1.1), (1.2) and their modifications have numerous applications in analysis.

\par\noindent
{\bf Questions:} 
\par\noindent
A. Since Theorem 1.2 has a constant that depends only on $m$ but not on the dimension of the ambient space,
one can ask if it has functional analogs (such as inequalities (1.1) or (1.2)) where the constant also
does not depend on the dimension of the ambient space. If true, this might be helpful in situations where the dimension $n$ is uncontrollably large, yet one is interested in integrals over subsets of a small dimension.
\par\noindent
B. What are the functional analogs of our Theorem 1.3 that would generalize inequalities (1.1), (1.2) to higher codimensions similarly to how the Michel-Simon Sobolev inequality on submanifolds of Euclidean spaces generalizes the classical Sobolev inequality ([MS])? In particular, can one always fill any given $M^n\subset R^N$, ($N>n+1$), by some $W^{n+1}$ so that some analogues of (1.1), (1.2) would hold on $W^{n+1}$, preferably with constants that do not depend on $N$? 

\subsection{Isoperimetric inequality in Banach spaces and related inequalities}

The classical isoperimetric inequality holds in codimensions greater than one: F. Almgren proved in [A] that for each $n<N$
any $n$-dimensional submanifold $M^n$ in $\mathbb{R}^N$ can be filled by an
$(n+1)$-dimensional minimal submanifold $W$ possibly with singularities
of codimension $\geq 2$ so that $vol_{n+1}(W)\leq \gamma(n)vol_n(M^n)^{\frac{n+1}{n}}$ with equality only in the case
when $M^n$ is a round $n$-dimensional sphere $S^n$ and $W$ is a flat $(n+1)$-dimensional disc $D^{n+1}$. Here, $\gamma(n)$ is defined by considering the equality case.

If one does not care about the optimal value of the constant $\gamma(n)$, then the same
inequality was proven earlier by J.Michel and L. Simon ([MS]) and later generalized
by M. Gromov ([Gr]) to the case when the ambient space is an arbitrary finite- or infinite-dimensional Banach space. We would like to note that long before [MS] H. Federer and W. Fleming proved a weaker version of the high-codimension isoperimetric inequality in Euclidean spaces ([FF]). In their version the constant $\gamma$ is allowed to depend not on $m$ but a potentially larger dimension $n$ of the ambient Euclidean space.

Gromov's result also contained the upper bound $c(n)vol^{\frac{1}{n}}(M^n)$ for the maximal distance between a point of a filling and the boundary $M^n$. In fact, this bound was the main goal of Gromov, as it almost immediately implies the upper bound $c(n)vol^{\frac{1}{n}}(M^n)$ for the length $sys_1(M^n)$ of the shortest noncontractible closed curve on $M^n$ provided that $M^n$ is essential. (Essentiality is a constraint on homotopy type of closed nonsimply-connected manifolds. All nonsimply-connected closed surfaces, $\mathbb{R}P^n$, tori, or, more generally, aspherical manifolds are essential.) In [Gr] the most important ambient Banach spaces were $L^\infty(M^n)$ as each compact metric space $X$ can be isometrically embedded in $L^\infty(X)$ via the Kuratowski embedding.

Gromov's proof was later simplified by S. Wenger ([W]). Later, the proof scheme was adopted in the paper [LLNR] to prove a similar isoperimetric inequality for $m$-dimensional Hausdorff contents, where $m$ is an arbitrary positive parameter not related to the dimension $n$ of $M^n$ (although the inequality is non-trivial only when $m\in (1,n]$). A main result of [LLNR] was a proof of a conjecture by L. Guth ([Gu17]) asserting that for each integer $m$ $(m-1)$-dimensional Urysohn width of a compact metric space $X$ can be majorized by $const(m)\HC_m^{1\over m}(X)$ . Here the $(m-1)$-dimensional Urysohn width $UW_{m-1}(X)$ is, by definition, the infimum over all $(m-1)$-dimensional complexes $K$ and continuous maps $\phi:X\longrightarrow K$ of $sup_{k\in K}\diam(\phi^{-1}(k))$. As each $X$ can be regarded as a subset
of $L^\infty(X)$, the Hausdorff content - Urysohn width inequality conjectured by Guth immediately follows from Theorem 1.5 (1) and from a similarly looking isoperimetric inequality in [LLNR] that differs from Theorem 1.5 only in part (2).

Later, P. Papasoglu found a much simpler proof of Guth's conjecture ([P]; see also [N] for an even simpler version of the proof with better constants. 
However, the approach in [P] and [N] does not seem to lead to a proof of isoperimetric (or boxing) inequalities. 
An immediate reason behind this failure is that the map $\phi$ in the definition of $UW_{m-1}$ can have an uncontrollably large
Lipschitz constant that depends on the fine properties of $X$. Easy $1$-dimensional examples show that this is
an intrinsic unavoidable feature of the Hausdorff content - Urysohn width inequality. Therefore, the mere existence of $\phi$ that satisfies Theorem 1.5(1) cannot provide any information about $\HC_m$ of the image of $H$.


Even before [LLNR] Guth ([Gu 13]) and independently, Young ([Y]) observed that the Federer-Fleming
approach can be used for Hausdorff contents. The central idea of the Federer-Fleming approach was to project 
from a random point inside a simplex (or a cube or a ball) to the boundary and to observe that the 
volume cannot increase by more than a constant; according to [Gu 13] and [Y] the same is true for Hausdorff contents. Then one can iterate and project on the skeleta
of progressively lower dimension until the skeleton of the desired dimension is reached.

This result was used in [LLNR] and is used in the present paper. However, this approach is not
sufficient to establish a weak version of the boxing inequality with constant depending on the dimension of the ambient Euclidean space (while being sufficient to establish a similar weak version of the corresponding isoperimetric inequality for Hausdorff contents).

The reason is that while one can bound the $(m+1)$-dimensional Hausdorff
measure (or Hausdorff content) of the cone $CX$ over an $m$-dimensional set $X$
in terms of the product of the $m$-dimensional Hausdorff measure (or Hausdorff content) of $X$ and
the maximum distance between the tip of the cone and a point in $X$,
one cannot bound the $m$-dimensional Hausdorff content of $CX$ using the same quantity. To see the difference,
consider $X$ with a small $\HC_m(X)$ inside a $n$-dimensional cube with side length $1$, ($n>m$).
Cover $X$ by a collection of balls with side lengths $r_i$ so that $\Sigma_i r_i^m< \HC_m(X)+\epsilon$ for a very small $\epsilon$. 
A projection from a random point might project a ball $B_i$ of radius $r_i$ in the covering
to a ball $B'$ with the radius $cost(n)r_i$ on the $(n-1)$-dimensional boundary of the cube. (Iterating, we will eventually project everything to the $(m-1)$-skeleton of the cube.)
However, the trace of homotopy between $B$ and $B'$ will look like the union of straight line segments of length $\sim 1$ connecting points of $B_i$ with their projections. One needs $\sim {1\over r_i}$ rescaled copies
of $B_i$ with radii between $r_i$ and $const(n)r_i$ to cover the trace of homotopy. If we want to prove the isoperimetric inequality, then we need to sum the radii of these balls to the power of $(m+1)$, obtaining the total contribution of $B_i$ $\sim {1\over r_i}r_i^{m+1}\sim r_i^m$.
After summing over $i$, we obtain a quantity that is controlled in terms of $\HC_m(X)$. However, here
we are concerned with the boxing inequality , that is, about the sum of radii of these balls to the power of $m$.
We obtain $\Sigma_i const(n){1\over r_i}r_i^m\sim \HC_{m-1}(X)$, and $\HC_{m-1}(X)$ cannot be controlled in terms of
$\HC_m(X)$. So, the idea of proof seems to fail. Note that the reason why it failed
is that the cubes used to cover $X$ were too spacious and contained too little
"material" from $X$. If the side length of a cube $C$ was less than or equal to $c(m) \HC_m^{1\over m}(X\bigcap C)$, then the above problem with the argument would not have arisen: The image
of homotopy is contained in $C$, and by our assumption $\HC_m(C)\leq const(m) \HC_m(C\bigcap X)$.

It should be noted that there exists an alternative approach by Y. Eliashberg and M. Gromov ([GE]) towards proving the isoperimetric inequality with constant depending on the dimension of the ambient space.
This approach does not involve radial projections from a random point. Instead, the authors of [GE] used orthogonal projections to carefully chosen hyperplanes. If applied to prove the boxing inequality, this approach would fail essentially for the
same reason as the Federer-Fleming style approach: One would need $~\sim {1\over r_i}$ copies of a ball of radius $r_i$ to cover the trace of a homotopy that would lead to an upper bound for $\HC_m$ of the trace of the homotopy in terms of $\HC_{m-1}(X)$ instead of $\HC_m(X)$ that we need.

\subsection{Some remarks about our proof.} 

Generally speaking, our proof of Theorem 1.5 follows the same scheme as
Wenger's version ([W])
of Gromov's proof of isoperimetric inequality in Banach spaces (which was also previously used in [LLNR] to prove isoperimetric inequalities for the Hausdorff content). 
The idea is to run an iterative process. On step $i$, one wants
to improve $f_i:X\longrightarrow B$, where $f_1$ is the inclusion of $X$ into $B$, to lower its volume (or, in later versions, $\HC_m$). In addition, $f_{i+1}$ and $f_i$ should be sufficiently close, so that one could find a nicely controlled homotopy
between them. These homotopies become pieces of the homotopy $H$
that we are constructing. At each step, the volume decreases by at least the same factor $>1$. Once the volume becomes lower than some very small $\epsilon$ that is allowed to depend on the dimension of the ambient space, one finishes the proof using
the Federer-Fleming argument that provides the last piece of $H$.
(This last piece requires $B$ to be finite-dimensional, so 
one would need to reduce the general case of the theorem to
the situation when the ambient Banach space is finite-dimensional albeit with an uncontrollable large dimension.)

At each step (with the exception of the last step),
one first finds a family of disjoint good metric balls,  (``thin fingers") in $X$ and cuts them off and fills the boundaries of the cuts with sets with
a nicely controlled volume or $\HC_m$. This means that one alters
$f_i$ in the interior of each of these balls so that (i) the volume (or $\HC_m$) of the new part is smaller than the volume (or $\HC_m$) of the removed part, and (ii) the new map $f_{i+1}$
is homotopic to the old map by means of a homotopy such that the volume (or a Hausdorff content) of its image is nicely controlled.

\forget
As the result, one obtains a new set (or, more precisely, a map of a set) with a considerably smaller value of $\HC_m$. Then we construct the homotopies
between $f_i$ and the newly constructed $f_{i+1}$ for the considered good balls. These homotopies 
will form pieces of
the filling of $X$ with a nicely controlled $\HC_m$ that we are constructing. Then we iterate this step again, and again, and again,
until we obtain (a map into) a set of very small Hausdorff content. Now we finish
the proof using a Federer-Fleming style argument.

The same scheme was also
used in [LLNR].
\forgotten

But, as we already mentioned, our main difficulty in comparison with the proof of
the isoperimetric inequalities for Hausdorff content in [LLNR] is that we cannot
use the cone inequality for Hausdorff content proven in [LLNR]. From the technical point of view the main novelty is our observation that the cone inequality is not really necessary to prove the boxing inequalities (as well as the isoperimetric inequalities for the Hausdorff content).

First, we observe that the cone inequality is not really needed to estimate $\hc_m$ of subsets of the image of the homotopy that connect the old version of $X$
before cutting a ``thin finger" and after cutting it and filling the boundaries. 

Moreover, our construction
of homotopies between ``thin fingers" that are being cut off and minimal fillings
of the boundaries are more straightforward and geometric in nature than in [LLNR]. 
(In [LLNR]
we were proving an effective extension theorem simultaneously with the isoperimetric
inequality for $\HC_m$ using induction with respect to $m$. In each induction step,
we use this extension theorem for $m-1$ to
construct these homotopies and complete the proof of the isoperimetric inequality for $\HC_m$. Nothing like that is needed here.)

Second, we also demonstrate that a control over $\HC_m$ of cones is not needed in the course of proving a weaker version of Theorem 1.5, where we assume that the ambient Banach space is finite-dimensional
and the constants are allowed to depend on its dimension, and not only on $m$ (Theorem 3.1). 

As we mentioned, we need to be able to reduce Theorem 1.5 to its particular case, where $B$ is finite-dimensional (as in [LLNR]).
An important difference between the present paper and [LLNR] is in this reduction.
Unlike [LLNR] we do not use the Kadec-Snobar projection, and do not need to prove versions of auxiliary lemmas for the version of Hausdorff contents, where all centers of balls must be in a prescribed linear subspace. Instead, we use a very simple geometric argument; see Section 1.6 below.

On the other hand, in this paper we still use the crucial idea from [LLNR] that was used there to counter the non-additivity of $\HC_m$. Namely, when we determine the ``thin fingers" that must be ``cut",
we use $\widetilde{\HC_m}$ instead of $\HC_m$, where $\widetilde{\HC_m}$ is a version
of $\HC_m$ for subsets of the ambient space $X$, where one is allowed to use only metric balls from a fixed almost optimal covering of $X$ by metric balls. In this way,
it is ensured that the cutting off of the thin fingers significantly reduces the Hausdorff content of $X$.

As in [Gr] or [LLNR], we need to first prove a weaker version of Theorem 1.5, where we assume that the ambient Banach space is finite-dimensional,
and the constants are allowed to depend on its dimension, and not only on $m$ (Theorem 3.1).

As in [G],  we use a covering by dyadic cubes $Q$ such that $\hc_m(X\cap Q)$ is neither too small nor too large compared to $\HC_m(Q)$. As in [Gu13] and [Y], we inductively 
project from each cube to its skeleton of codimension one, until we reach the
skeleton of dimension $\lceil m\rceil-1$. 
But now we need to deal with a cubic complex that consists of cubic faces of different dimensions.
As the cubes are of different sizes, the union of their $\lceil m\rceil-1)$-skeleta will, in general, be disconnected.
Therefore, the obvious problem is that we cannot hope to construct a meaningful continuous map of $X$ to this disconnected
set. In the first version of the paper ([AN]) we used the idea of ``insertions" connecting
the $(\lceil m\rceil -1)$-dimensional skeleta of the cubes, obtaining a connected set that can be used
to obtain a continuous map with the desired properties ((see Figure 1 and the explanation of this idea in Section 1.5 of [AN]). 

However, we discovered that this approach works only if the adjacent cubes have comparable sizes. This can be achieved using the idea of ``collars" of dyadic cubes explained in Section 3 below.
An anonymous referee of the first version of this paper pointed out 
that, since we use cubes of comparable sizes, one can also use the idea
of QC simplicial complexes from [Y], that is, metric simplicial complexes where all simplices are $c$-bilipschitz homeomorhic to rescalings of standard simplices for some constant $c$.
Indeed, we triangulate cubic complexes where all the cubes are dyadic cubes, and the side lengths of adjoint cube differ by at most the factor of $2$ that were
used in the previous version of the paper into ``fat" simplices. Then we are able to apply a more conventional version
of the Federer-Fleming argument to the resulting QC triangulations. We adopted this approach in the present version of our paper. (Adjoint simplices must have comparable sizes since otherwise a radial projection of a set $X$ in a simplex $\sigma^m_1$ to its faces can completely fill a $(m-1)$-dimensional face $F\subset\partial\sigma^m_1$ of a much smaller adjoint simplex $\sigma^m_2$.)
Also note that the idea
of using a covering by dyadic cubes of comparable sizes was previously known in analysis and appeared in the classical Whitney covering lemma. Yet its application in the present paper in combination with the Federer-Fleming projection towards an isoperimetric inequality style result seems to be new.

Finally, we observe that, as we do not use the cone inequality, our proof adapts to the situation where the ambient space is not a Banach space, but a metric space with a linear contractibility function. The proof is somewhat easier when the ambient space is a Riemannian manifold, but it can be carried out with some modification of our basic construction in the case of more general ambient spaces. The proof of Theorem 3.1 cannot be adapted to this situation when the ambient space is not linear, but we circumvent this difficulty by embedding the ambient space $M$ into a Banach space, applying Theorem 3.1 in this Banach space, and then projecting $K$ and the image of the homotopy back to $M$. 
\forget
Therefore,
a more significant new idea is that we replace our original cover with dyadic cubes
with a coarser cover using the idea of ``collars". This idea is explained in Section 3.3. (See Figures~\ref{figure:paths} and ~\ref{figure:collar} that illustrate this idea). It
ensures that sizes of dyadic cubes (=balls) in the covering change ``slowly", which turns out to be sufficient to complete the proof. (This approach works even despite the fact that we lose a positive lower bound for ${\HC_m(X\cap Q)\over \HC_m(Q)}$ for some of the balls $Q$ of the new cover.) 

\begin{figure}[ht]
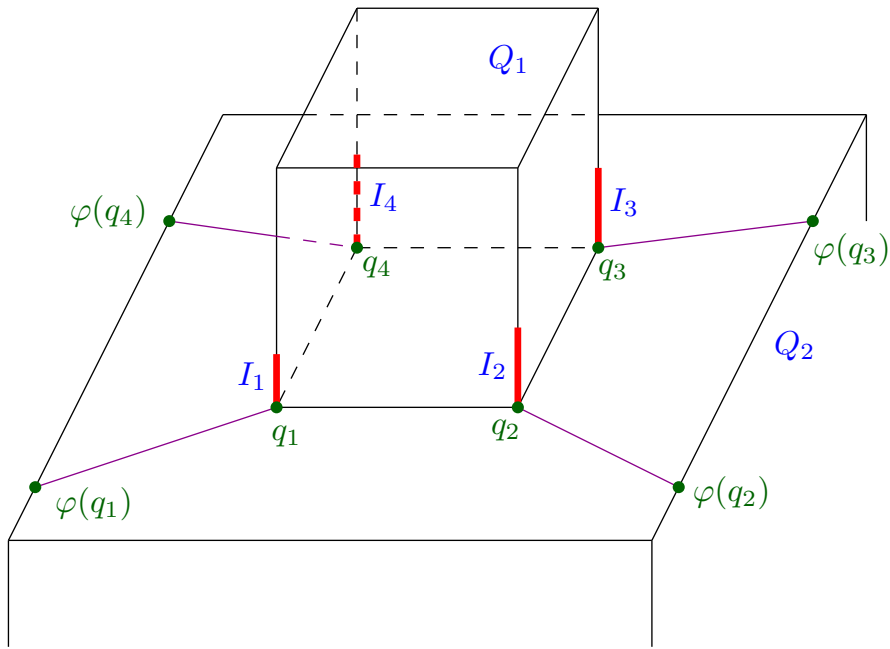

\center
\include{Cubes1.pdf}
\caption{Using insertions.}
\label{figure:twocubes}
\end{figure}

Finally, we use the idea of ``insertions". To explain this idea imagine a situation when we want to map a connected set $X\subset B$, where $B$ is the $3$-dimensional space with $l^\infty$ metric, into a $1$-dimensional complex at the distance $\leq c_1(2) \HC_2(X).$ Assume that $X$ is covered by two dyadic cubes $Q_1$ and $Q_2$ in $\mathbb{R}^3$ that are 
of different sizes, so that the smaller cube $Q_1$ is on top of $Q_2$ and the bottom $2$-dimensional face of $Q_1$
is a subset of the interior of the top $2$-dimensional face of $Q_2$ (see Figure  ~\ref{figure:twocubes} ). Assume that the diameters
of $Q_1$ and $Q_2$ are less than the upper bound desired in Theorem 1.5 (1), and that the sum
of squares of the side lengths of $Q_1$ and $Q_2$ does not exceed $c_2(m)\HC_2(X)$.
It is natural to want to attempt to map $X$ in the union of the $1$-skeletons of $Q_1$ and $Q_2$, but, of course, this would be impossible because this union is not connected.

Let us explore where this attempt will surely fail.
Assume that we know how to map $X\cap Q_1$
into the $1$-skeleton of $Q_1$ so that the intersections of $X$ with all closed faces of $Q_1$ remain in these faces. Denote this map by $\phi_1$. In addition, assume that we managed to project the intersection of $X$ with the interior of $Q_2$ radially to the boundary of $Q_2$, and, moreover, know how to project the images of $X$ in all $2$-dimensional faces  of $Q_2$ under this
map as well as under $\phi_1$ to the boundaries of these faces.
Once we perform all these projections, the result will be in the union
of the $1$-skeletons of $Q_1$ and $Q_2$, as we wanted. However, after performing
the projection of $X$ in the top face of $Q_2$, we will create discontinuity at the four
points $q_1,q_2, q_2, q_4$ at the bottom of the vertical edges of $Q_1$: While points in $\phi_1(X)$ arbitrarily close to $q_i$, $i=1,2,3,4$, will be unchanged by the projection $\phi$ 
in the top face of $Q_2$, the points $q_i$ will be mapped to different points $\phi(q_i)$
in the boundary of the top face of $Q_2$. Our idea now is to eliminate this discontinuity
as follows: (a) We define the $1$-dimensional complex $K$ from Theorem 1.5 as the union of four straight line intervals $[q_i\phi(q_i)]$ in the upper face of $Q_2$, and the $1$-skeleta of $Q_1$ and $Q_2$; (b) We will map some short intervals $I_i$ in the vertical edges of $Q_1$ incident to $q_i$ into $I_i\cup [q_i\phi(q_i)]$ regarded as a longer interval. Here $q_i\in I_i$
is mapped to $\phi(q_i)$ and the other endpoint of $I_i$ is mapped into itself.
It turns out that adding four intervals $[q_i\phi(q_i)]$ in this way works, and that this
idea can be generalized to all higher dimensional situations.
\forgotten


\forget
is an adaptation of the proof of the isoperimetric
inequality for Hausdorff contents in [LLNR] to the boxing inequality. It will be presented in section 2.
It follows the same general scheme as [Gr], [W], [Gu 17], and [LLNR]. 
We use the ideas from [LLNR] on how to overcome the non-additivity of Hausdorff contents, most notably using $\widetilde{\HC}_m$ (that will be explained in section 2.2 below).
Yet we made a number of modifications and simplifications to the proof of the isoperimetric inequality for Hausdorff contents that was presented in [LLNR]. For example, here we do not use the cone inequality for Hausdorff contents proven in [LLNR]. Also, in the present work we designed a much simpler reduction to a finite-dimensional case than in [LLNR] that we are going to explain in the next subsection.
\forgotten
\subsection{Reduction to a finite-dimensional case}
Theorem 1.5 can be easily reduced to the case where $B$ is a finite-dimensional space. For a sufficiently small $\epsilon>0$ cover $X$ by a finite collection of open metric balls $B_i=B(x_i,r_i)$ so that
$\Sigma_i r_i^m< \HC_m(X)+\epsilon$. Without loss of generality one can assume that $X$ is already covered
by the union of $B(x_i, r_i-\delta)$ for a very small positive $\delta$.
Let $L$ denote the finite-dimensional space spanned by the centers $x_i$ of all metric balls of the covering. Let $\phi_i$ be a partition of unity subordinate to the covering $\{B_i\}$. As usual, $\phi_i$ are non-negative continuous real-valued functions on the ambient Banach space, vanishing outside of $B_i$, such that for each $x\in X\subset \bigcup_i B(x_i,r_i-\delta)$ $\Sigma_i\phi_i(x)=1$.
Define a map $\psi:X\longrightarrow L$ using the formula $\psi(x)=\Sigma_i \phi_i(x)x_i$. Let $X_1$ denote the image $\psi(X)$. For each $x\in X$ let $I(x)$ denote the (finite) set of indices $i$ such that $\phi_i(X)\not= 0$.
Assume that $B_k$ is a ball with the largest radius in the collection of balls $B_i$, $i\in I(x)$, that is, $r_k=\max_{i\in I(x)}r_i$. We claim that (1) All centers $x_i$, $i\in I(x)$, are contained in $2B_k=B(x_k, 2r_k)$; (2) the straight line segment between $x$ and $\psi(x)$ is contained in $2B_k$.
Indeed, as the intersection of all balls $B_i$, $i\in I(x)$, is non-empty (because it contains $x$), these balls
pairwise intersect. Therefore, the distance between $x_i$ and $x_k$ does not exceed $r_i+r_k\leq 2r_k$,
and (1) follows. As all $x_i$ are in $2B_k$, and $2B_k$ is convex, $\psi(x)$ is in $2B_k$.
Therefore, the segment $[x\psi(x)]\in 2B_k$. As this is true for all $x$, we see that (1) $X_1\subset \bigcup 2B_i$; and (2) the image of the homotopy $H_1:X\times [0,1]\longrightarrow B$ between the identity map and $\psi$, defined as $H_1(x,t)= (1-t)x+t\psi(x)$ is also contained in $\bigcup_i 2B_i$. However,
$\HC_m(\bigcup 2B_i)\leq \sum_i (2r_i)^m\leq 2^m\HC_m(X)+2^m\epsilon.$

Now observe that if Theorem 1.5 is true for all finite-dimensional ambient spaces, we can apply it to $X_1\subset L$ and obtain a $(\lceil m\rceil-1)$-dimensional complex $K\subset L$, a map $\phi_0:X_1\longrightarrow K$, and a homotopy $H_2$ between the identity map and $\phi_0$ as in Theorem 1.5.
Now we can construct the desired $\phi$ as the composition of $\psi$ and $\phi_0$ and the homotopy $H$ as a concatenation of $H_1$ and $H_2$.

Therefore, in the following, we can assume that $B$ in Theorem 1.5 is finite-dimensional.

\section{Reduction of Theorem 1.5 to its weak version, where the constants in the right hand sides of (1), (2) depend on the dimension of $B$}

In this section, we will deduce Theorem 1.5 from its weaker version, where the constants
$c_1(m)$, $c_2(m)$ are allowed to depend on the dimension of the (finite-dimensional)
Banach space $B$ (and so are $c_1({\rm dim} B), c_2({\rm dim} B)$). This weaker version will be proven in Section 3.

Here we follow Gromov's ``cutting-off thin fingers" scheme of the proof of the main result in [LLNR] with some improvements indicated below. In turn, the scheme of the proof in [LLNR] is similar
to the scheme of the proof of the isoperimetric inequality in [Gr], and its later simplified version from [W]. However, [LLNR] contains several new ideas, in particular, tools to counter the nonadditivity of Hausdorff content, and these ideas are also used here.

The main
difference between our proof and proofs [W] and [LLNR] is that we cannot use the cone inequality (since the $m$-dimensional Hausdorff content of the cone $CX$ over $X$ cannot be majorized
in terms of $HC_m(X)$). So, we do not use it. Both [W] and [LLNR] use gradual "improvements" of the underlying manifold/metric space, where one cuts off "thin fingers" inside ``good" balls $B_i$ and replaces them by ``good filling" of the boundaries. Then one uses the coning to construct a homotopy between the removed and inserted pieces. These homotopies become the parts of the homotopy that we construct, and we need to majorize
$\HC_m$ of their images.
We made an observation that one does not need to evaluate the volume/Hausdorff content of the cones
over the union of the old and new parts of $X$ (like it was done in [W] and [LLNR]), but instead can majorize it by $HC_m(B_i)$ (equal to the $m$th power of its radius) - see Lemma 2.2.

Unlike [W], here and in [LLNR] we consider not only geometric objects $X_i$, but also their
maps to $K$. In [LLNR] we were also proving an effective version of the Dugundji theorem by obtaining an upper bound for the minimal $HC_m$ of the extension image.
This theorem was proven simultaneously with the isoperimetric inequality for Hausdorff contents and was used in the proof of the induction step for the isoperimetric inequality.  In contrast with [LLNR], here we are not proving a quantitative extension theorem, and give a direct construction of the required maps.
This helps to significantly simplify and shorten the proof. Also, we hope that in the future
this might be helpful for a dimensionality reduction algorithm based on the ideas
of this proof.

\subsection{ Plan of the proof.}

%

First, we use induction with respect to $m$. The base of induction is the case $m\in (0,1]$. In this case, Theorem 1.5 has already been proven in [LLNR]. The basic
idea is that, when $m\leq 1$, a replacement of two intersecting balls by a metric ball of the minimal radius that contains both of these balls does not increase the sum of the $m$th powers of their radii as this minimal radius does not exceed the sum of the radii of these two balls.
In fact, for $m\leq 1$ $r_1^m+r_2^m\geq (r_1+r_2)^m$ and, more generally, $\Sigma_i r_i^m\geq (\Sigma_i r_i)^m$. Proceeding by induction, we see that the same will be true for any finite collection of metric balls in $B$ such that their union is connected.

As a result, one can assume that the cover of
$X$ by closed metric balls that yields the minimal $\HC_m$ consists of pairwise disjoint metric balls $\hat B_j$ with radii that do not exceed $\Sigma_{i\in C_j}r_i$, where $\{C_j\}$ is a partition of the set of indices $i$. 
As $m\leq 1$, $\HC_m$ of the union of these disjoint balls does not exceed $\Sigma_j(\Sigma_{i\in C_j}r_i)^m\leq \Sigma_j\Sigma_{i\in C_j}r_i^m= \Sigma_i r_i^m$.
The $0$-dimensional complex $K$ is the union of the centers of these balls, the map $\phi$ sends a part of $X$ inside each ball to its center, and the image of the (obvious) homotopy $H$ is also contained in the union of these balls, resulting in
the inequality $\HC_m(H(X\times [0,1]))\leq \HC_m(X)$.

To prove the induction step, we assume that the theorem is true for $m-1$ and prove it for $m>1$. Theorem 3.1 can be immediately used to deduce Theorem 1.5 in the case when $\HC_m(X)$ does not exceed
a small positive constant that depends on the dimension of the ambient space.
If $\HC_m(X)$ does not satisfy this condition, we will inductively construct $X_0=X$, $X_1, \ldots $ and maps $\phi_i:X_i\longrightarrow X_{i+1}$ such that: (1) $\HC_m(X_i)$ exponentially decreases with $i$, where the exponent depends on $m$ but not on $n$; (2) For every $i$ and $x\in X_i$ $\Vert \phi_i(x)-x\Vert\leq const(m)\HC_m^{1\over m}(X_i)$. This construction will be explained in the next two subsections.

At the end of our construction, we will proceed as follows.
As $\HC_m(X_i)$ exponentially decreases with $i$, eventually it becomes arbitrarily small.
In particular, there exists
$N$ such that $\HC_m(X_N)\leq \min\{c_1(n)^m, c_2(n)\}^{-1}\HC_m(X)$, where $n$ is the dimension of
the ambient Banach space. (Recall that in light of the previous section we can assume that the ambient Banach space is finite-dimensional). Here $c_1(n), c_2(n)$ are the constants
in the weaker version of Theorem 1.5 stated as Theorem 3.1 below, where the constants are allowed to depend on $n$.
Then we can apply Theorem 3.1 for $X_N$ and obtain
a map $\phi_N:X_N\longrightarrow K$ and a homotopy $H_N:X_N\times [0,1]\longrightarrow K$
so that $K$ is a $(\lceil m\rceil-1)$-complex, and $\Vert \phi_N-x\Vert\leq c_1(n)\HC_m^{1\over m}(X_N)\leq \HC_m^{1\over m}(X)$, the lengths of the trajectories of points of $X_N$ during the homotopy $H_N$ similarly do not exceed $\HC_m^{1\over m}(X)$ and $\HC_m(H_N(X_N\times [0,1]))\leq \HC_m(X)$. Now we can define the desired map of $X$ to $K$ as the composition
of $\phi_1,\ldots , \phi_N$, and homotopy $H$ as the concatenation of homotopies $H_i$,
where for each $i<N$, homotopy $H_i$ moves all points $x$ of $X_i$ through the straight line
segment between $x$ and $\phi_i(x)$ in $B$. The quantity $\Vert x-\phi(x)\Vert$ does not exceed $\Sigma_{i=0}^N \Vert \phi_i(x)-\phi_{i+1}(x)\Vert\leq \Sigma_{i=1}^N const(m)\ \HC_m^{1\over m} (X_i)$, where $\phi_0(x)=x$. The last sum can be majorized by the sum of an infinite geometric series equal to $Const(m)\HC_m^{1\over m}(X)$ (as $\HC_m(X_i)$ exponentially decreases). Exactly the same 
argument implies that $\HC_m$ of the image of $H$ can be bounded
by $Const_2(m)\HC_m(X)$.

One may ask if the last step (for $X_N$) based on Theorem 3.1 
is really needed here. An alternative idea
might be to continue to perform the inductive steps for all $i$, and hope to get the desired $\phi$ as $\lim_{i\longrightarrow\infty}\phi_i\circ\phi_{i-1}\circ\ldots\circ\phi_1$.  In fact, for each $x\in X$ the distances 
$\Vert \phi_i\circ\phi_{i-1}\circ\ldots\circ\phi_1(x)-\phi_{i-1}\circ\ldots\circ\phi_1(x)\Vert$ exponentially decrease with $i$. Therefore, the sequence
$\phi_i\circ\phi_{i-1}\circ\ldots\circ\phi_1(x)$ converges, and one can define $\phi(x)$ as the limit. The problem is that
the Lipschitz constants of the maps $\phi_i\circ\phi_{i-1}\circ\ldots\circ\phi_1(x)$
are unbounded, and we cannot guarantee that the convergence is uniform, that so defined $\phi$ is continuous, and that $\phi(X)$ has zero $m$-dimensional Hausdorff content.
As we were unable to overcome these difficulties, we abandoned this approach.

Thus, modulo the last step the proof reduces to the induction step that will
be discussed in the next two subsections.

\subsection {Induction step: Good balls}

The proof of the classical boxing inequality in [G] as well as the
proofs of the high codimension isoperimetric inequality in [Gr] and [W] use coverings
of a manifold $M^n$ by ``good" balls. Depending on the context, the exact requirements for good
balls are different, but typically a good ball $B$ of radius $r$ centered at $x$ must satisfy
the requirement that the ratio of the measure of $M^n\cap B$ to $r^n$ is in an interval
$[const_1(n), const_2(n)]\subset (0,1)$ that depends only on the dimension of $M^n$.
The simplest way to construct such a ball centered at a prescribed point $x\in M^n$ is
to start from a ball $B=B(x,R)$ of a very large radius $R$, so that ${\vert M^n\cap B\vert\over R^n}$ is very small due to the compactness of $M^n$, and then gradually lower the value of
$R$ to the maximal value $r_0$ such that this ratio is equal to $const_2(n)$. 
This would imply that for $r_1=2r_0$, $\vert M^n\bigcap B(x, r_1)\vert < const_2(n) r_1^n$, and the measure of the annulus $A(r_0,r_1)=B(x, r_1)\setminus B(x, r_0)$ does not exceed
$2^nconst_2(n)r_0^n$.

Another requirement (in [W]) is that the boundary of $B$ has a ``small" $(n-1)$-dimensional measure. This can be ensured by applying the coarea formula to the annulus
$A(r_0,r_1)$, where $r_1$ is chosen as a multiple of $r_0$, say $2r_0$. The coarea formula implies the existence of some $r_*\in (r_0, r_1)$
such that $\vert M^n\cap\partial B(x,r_*)\vert \leq 2^n const_2(n)r_0^{n-1}\leq 2^nconst_2(n)r_*^{n-1}$.
On the other hand $\vert M^n\cap B(x,r_*)\vert< const_2(n)r_*^n$, and as $B(x,r_*)$ contains $B(x,r_0)$, $\vert M^n\cap B(x,r_*)\vert\geq const_2(n)r_0^n\geq {const_2(n)\over 2^n}r_*^n$.
Therefore, $B(x, r_*)$ satisfies 
all three conditions if $const_1(n)={const_2(n)\over 2^n}$. 
Now we can use the Vitali construction to obtain a finite family of disjoint ``good" balls with radii $r_i$ so that concentric balls with radii $3r_i$ cover $M^n$.
The construction of ``good" balls implies that $\vert M^n\cap B(x,3r_i)\vert\leq c(n)\vert M^n\cap B(x_i, r_0)\vert$
for each of these "good" balls.

It seems natural that if one is interested in $\HC_m$ rather than volumes, one can
try to repeat the same argument for Hausdorff contents.
One of the new ideas of [LLNR] was to use this argument with the functional $\HC_m$ on the set of all subsets
of $X$ replaced by another functional denoted $\widetilde{\HC}_m$. 
To define $\widetilde{\HC}_m$ consider an almost optimal covering of $X$ by a finite collection of metric balls $\{\beta_i\}$ with radii $r_i$ so that $\Sigma_i r_i^m<\HC_m(X)+\epsilon$ for a very small $\epsilon$. We also assume that no ball $\beta_j$ is contained in another ball $\beta_k$. Now for each subset $A$ of $X$ we define $\widetilde{\HC}_m(A)$ as the minimum
of $\Sigma_i r_i^m$, where the minimum is taken over all subcollections of $\{\beta_i\}$ that cover $A$. Observe, that 1) $\widetilde{\HC}_m(A)\geq \HC_m(A)$, and $\widetilde{\HC}_m(X)$ and $\HC_m(X)$ are almost equal (differ my at most $\epsilon$); 2) $\widetilde{\HC}_m$ of any ball of a sufficiently small
radius $r$ centered at a point $x$ of $X$ is greater than or equal to $r_*^m$, where $r_*>>r$ is the smallest radius of a ball from the collection $\{\beta_i\}$ containing $x$; 3) $\widetilde{\HC}_m$ (as well as $\HC_m$) satisfy the following coarea inequality (Lemma 2.3 in [LLNR]; see also the discussion before Lemma 2.2 and Example after Lemma 2.3 in [LLNR]):
For all $x\in X$. and positive $r_1$, $r_2>r_1$ there exists $r\in (r_1,r_2)$ such that
$\widetilde{\HC}_m(X\bigcap S(x,r))\leq {2\over r_2-r_1}\widetilde{\HC}_m(A(x, r_1, r_2))$. Here $S(x,r)$ denotes the metric sphere in $X$ centered at $x$ of radius $r$, and $A(x,r_1,r_2)$ is the metric annulus in $X$ centered at $x$ bounded by metric spheres with radii $r_1$ and $r_2$. (The same inequality holds for $\HC_m$.) As $\widetilde{\HC}_m$ is not additive, we cannot estimate $\widetilde{\HC}_m(A(x,r_1,R_2))$ by subtracting $\widetilde{\HC}_m(B(x,r_1))$ from $\widetilde{\HC}_m(B(x,r_2))$. Instead, we use $\widetilde{\HC}_m(B(x,r_2))$ as an upper bound for $\widetilde{HC}_m(A(x,r_1,r_2))$.

Now, we can define ``good" balls using $\widetilde{\HC}_m$ (and $\widetilde{\HC}_{m-1}$) and proceeding exactly as described above. (Note that $\widetilde{\HC}_{m-1}$ here and below is calculated using the collection of balls $\{\beta_j\}$ used in a definition
of $\widetilde{\HC}_m$, and {\it not} an analogous collection that almost minimizes $\Sigma_i r_i^{m-1}$ among covers of $X$ by metric balls.)
A small change made in [LLNR] is that one
looked for $r_*$ not in the $(r_0,\ 2r_0)$ interval but $((1+{1\over m})r_0, (1+{1\over m})^2r_0)$, where the choice of the constant $(1+{1\over m})^2$ instead of $2$ was made
to make some constant smaller, and the fact that the left end of the interval was chosen to be strictly greater than $r_0$ will be explained later. The following theorem was proven in [LLNR] (see the proof of Proposition 3.3) for $C_1(m)=1+{1\over m}$, $C_2(m)=(1+{1\over m})^2, A(m)>m$. However, the proof for arbitrary constants $C_1(m),C_2(m)$ is identical.

\begin{theorem}([LLNR]) Assume that $C_1(m), C_2(m), A(m)$ 
are not necessarily continuous positive functions defined for all $m>1$ such that $C_2(m)>C_1(m)>1$, $A(m)>{2\over C_1(m)-1}$. 
Then for each Banach space $B$ and its compact subset $X$
there exist a positive integer $N$, a finite set $\{x_i\}_{i=1}^N$ of points of $X$, a set of $N$ positive real numbers $r_i$, $i=1,\ldots, N$, and a set of positive real
numbers $R_i\in [C_1(m)r_i, C_2(m)r_i]$ such that:
\par\noindent 
(1) Metric balls $B(x_i, R_i)$ are disjoint, but balls $B(x_i, 3R_i)$ cover all $X$;
\par\noindent
(2) For each $i$ $\widetilde{\HC}_m(B(x_i,r_i)\cap X)=({r_i\over A(m)})^m$, and 
for each $r>r_i$
$\widetilde{\HC}_m(B(x_i,r)\cap X)< ({r\over A(m)})^m$;
\par\noindent
(3) $\widetilde{\HC}_{m-1}(\partial B(x_i,R_i)\cap X)\leq {2\over C_2(m)-C_1(m)}C_2(m)^{m-1} {r_i^{m-1}\over A(m)^m}$;
\par\noindent
(4) $\widetilde{\HC}_m(X\setminus \cup_{i=1}^N B(x_i, R_i))\leq \widetilde{\HC}_m(X)-\Sigma_{i=1}^N \widetilde{\HC}_m(B(x_i,r_i)\cap X)$;
\par\noindent
(5) $\Sigma_{i=1}^N\widetilde{\HC}_m(B(x_i,r_i)\cap X)\geq {1\over (3C_2(m))^m}\widetilde{\HC}_m(X)$. This inequality and (4) imply that $\widetilde{\HC}_m(X\setminus \cup_{i=1}^N B(x_i,R_i))\leq (1-{1\over (3C_2(m))^m})\widetilde{\HC}_m(X)$.
    
\end{theorem}
Before the text of this theorem we explained how to choose $r_i$ and $R_i$ to satisfy (1)-(3). An outstanding issue is that,
as $\lambda(r)={\frac{\widetilde{\HC}_m(B(x_i,r)\cap X)}{r^m}}$
is, in general, not continuous, it may not
be immediately clear that for $r_i$ defined as $\sup\{r| \lambda(r)\geq \frac{1}{A^m(m)}\}$,
$\lambda(r_i)=\frac{1}{A^m(m)}$, which is equivalent to (2).  (Note that as $x_i\in X$,
$\lim_{r\longrightarrow 0^+}\lambda(r)=\infty$. On the other hand, $\lim_{r\longrightarrow\infty}\lambda(r)=0$. Therefore, $r_i$ exists and is finite.)
Observe that the
monotonicity of $\widetilde{\HC_m}$, and the continuity of $r^m$ imply that for each $\rho$
$\lim \inf_{r\longrightarrow \rho^+}\lambda(r)\geq \lambda(\rho)\geq \lim \sup_{r\longrightarrow \rho^-}\lambda(r)$. Combining the definition of $r_i$ with the first of these two inequalities for $\rho=r_i$, we see that $\lambda(r_i)\leq \frac{1}{A^m(m)}$.
The second inequality and the definition of $r_i$ similarly imply that $\lambda(r_i)\geq \frac{1}{A^m(m)}$, which completes the proof of (2).
%

On the other hand, we did not explain (4). Here we need to use $\widetilde{\HC}_m$ instead of $\HC_m$ and
to use the assumption that $C_1(m)>1$. (This assumption was not needed when we dealt with volumes). The underlying difficulty here is the lack of additivity of the Hausdorff contents.

The idea behind the proof of (4) is the following one. Observe that the balls of the collection $\{\beta_i\}$ that are used in an optimal covering of $B(x_i,r_i)\cap X$ have radii
$\leq \widetilde{\HC}^{\frac{1}{m}}_m(B(x_i,r_i)\cap X)= {r_i\over A(m)}.$ Therefore, all these balls are contained in $B(x_i,R_i)$. (Here, we use the assumption that ${2\over A(m)}<C_1(m)-1$).
Hence, these balls cannot be used in an optimal covering of $X\setminus\cup_{i=1}^N B(x_i,R_i)$ by balls $\beta_i$. Even if all the remaining balls are used to cover $X\setminus\cup_{i=1}^N B(x_i,R_i)$, the sum of their $m$th powers does not exceed the right-hand side in inequality (4).

In order to prove (5) first note that $X\subset\cup_{i=1}^N B(x_i, 3R_i)$, $\widetilde{\HC}_m(B(x_i, 3R_i)\cap X)\leq ({3R_i\over A(m)})^m$, and
$3R_i\leq 3c_2(m)r_i$. Also, recall that $(\frac{r_i}{A(m)})^m=\widetilde{\HC}_m(B(x_i,r_i)\cap X)$. Therefore, $\widetilde{\HC}_m(X)\leq (3c_2(m))^m\Sigma_{i=1}^N\widetilde{\HC}_m(B(x_i,r_i)\cap X)$.

We would like to finish this section with an easy-to-prove inequality that could have been stated and proven in [LLNR] (but was not).

\begin{lemma} In the notation of the previous theorem
$$\Sigma_{i=1}^N ({r_i\over A(m)})^m\leq \widetilde{\HC}_m(X).$$
\end{lemma}
\begin{proof}
Consider coverings of balls $B(x_i,r_i)\cap X$ by balls from the collection $\{\beta_j\}$ such that the sums of the $m$th powers of the radii are minimal. All balls from collection $\{\beta_j\}$ used in such a covering for $B(x_i,r_i)\cap X$ have radii $\leq {r_i\over A(m)}$. Therefore, no pair of balls from the collection $\{\beta_j\}$ used in two of these coverings for different values of $i$ can intersect. Indeed, these two balls are contained in the interiors of the disjoint balls $B(x_i,R_i)$ for two different values of $i$. So, a ball from the collection $\{\beta_j\}$ used to cover one of $B(x_i,r_i)\cap X$ does not appear in the covering
of $B(x_k, r_k)\cap X$ for any other value of $k$. In other words, when one combines the optimal subcollections of balls from $\{\beta_j\}$
used to cover $B(x_i,r_i)\cap X$ for all values of $i$, the resulting collection contains no more than one copy of each ball from $\{\beta_j\}$. Denote the resulting collection $\{\beta_{j_k}\}$. We see that the sum $\Sigma_k r_{j_k}^m$  of the $m$th powers of the radii of all balls $\{\beta_{j_k}\}$ is equal to $\Sigma_i \widetilde{\HC}_m(B(x_i,r_i)\cap X)$. On the other hand, this collection is a subcollection
of the original collection $\{\beta_j\}$, and $\Sigma_k r_{j_k}^m\leq \Sigma_j r_j^m=\widetilde{\HC}_m(X)$. Combining this inequality with the equality $\Sigma_kr_{j_k}^m=\Sigma_i\widetilde{\HC}_m(B(x_i,r_i)\cap X)$ we obtain the inequality
in the lemma.
\end{proof}

Note that we are going to use $\widetilde{\HC}_m$ below not only for (subsets of ) $X=X_0$
but for all $X_k$. The construction of $\phi_k$ and $X_{k+1}$ 
will start with a choice of a nearly optimal covering of $X_k$ by metric balls
that will be used to define the functional $\widetilde{\HC}_m$ defined for all subsets of $X_k$. 

\subsection{Induction step: the construction}

Here we construct maps $\phi_k:X_k\longrightarrow X_{k+1}$ that lower $\HC_m$ by at least a constant factor $>1$ and do not move points too far in $B$. Fix an optimal
covering $\{\beta_j\}$ of $X_k$ by metric balls and find the disjoint balls $B(x_i, R_i)$
as in Theorem 2.1 of the previous section (for $X_k$ instead of $X$). Here, the balls $\{\beta_j\}$ were (and will be) used to determine $\widetilde{\HC}_m$ and $\widetilde{\HC}_{m-1}$ of subsets of $X_k$, yet $\widetilde{\HC}_m(X_k)$ exceeds $\HC_m(X_k)$ by not more than an arbitrarily small summand $\tau$.
We are going to perform the same procedure for all values of $i$. As a result, we will remove the intersection of $X_k$ with the interior of $B(x_i,R_i)$ and replace it by something else so that $\widetilde \HC_m$ decreases. After performing this procedure for all values of $i$,
we will obtain $X_{k+1}$ and $\phi_k$ is going to be constructed while we construct $X_{k+1}$.

\begin{figure}
    \centering
    \includegraphics[scale=0.7]{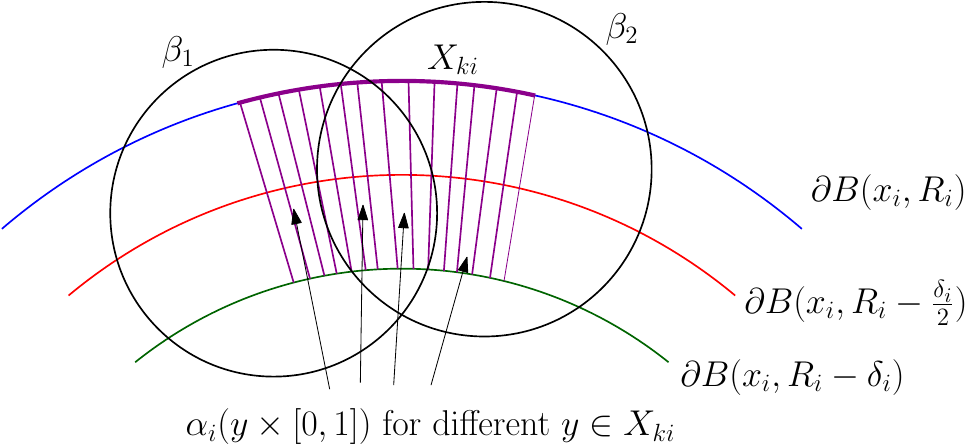}
    \caption{We added the union of the metric balls $\beta$ covering $X_k\cap\partial B(x_i, R_i)$ to $X_k$
    and denoted the result by $\tilde X_k$.
    These balls $\beta$ cover all points of $X_k$ in a $\delta$-thick metric annulus bounded by $\partial B(x_i,R_i)$ on the outside. We identify the intersection of this annulus with $\tilde{X}_k$ and
    $(\tilde{X}_k\cap\partial B(x_i,R_i))\times [0,1]$ by means of homeomorphisms $\alpha_i$.}
    \label{figure:homeomorphism}
\end{figure}

\begin{figure}
    \centering
    \includegraphics[scale=0.7]{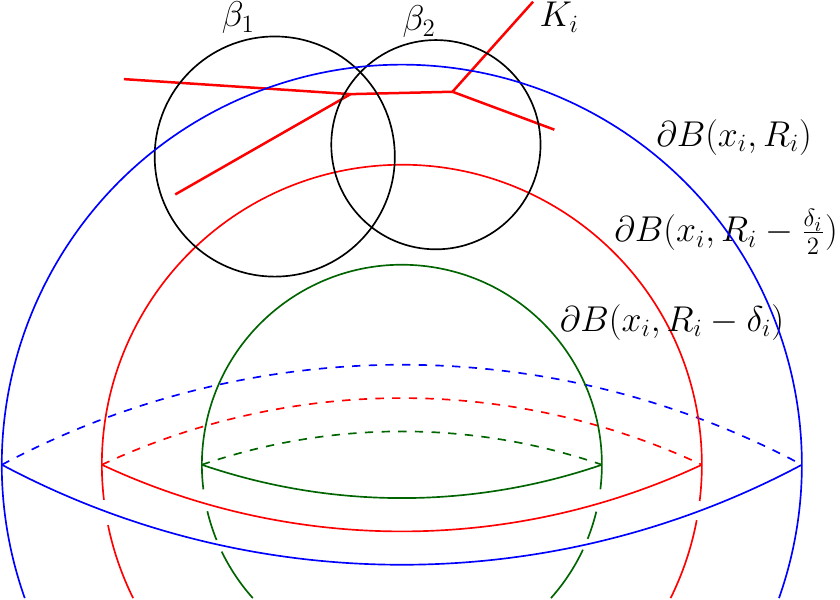}
    \caption{The part of $X_k$ inside the closed ball $B(x_i,R_i-\delta_i)$ is mapped to $x_i$. Points of $X_k$ outside of the interior of the ball $B(x_i,R_i)$ are mapped to themselves.
    Points of $X_k$ on the middle sphere $\partial B(x_i,R_i-\frac{\delta_i}{2})$ are mapped to $K_i$ using $\theta_i$. Two annuli between spheres are sliced into radial arcs, and the arcs that are both in $X_k$ and the inner annulus are mapped into the straight line segments between the images of their endpoints. The arcs that are both in $X_k$ and the outer annulus are mapped using the homotopy $H_i$.}
    \label{figure:emptyingball}
\end{figure}

\begin{figure}
    \centering
    \includegraphics[scale=0.7]{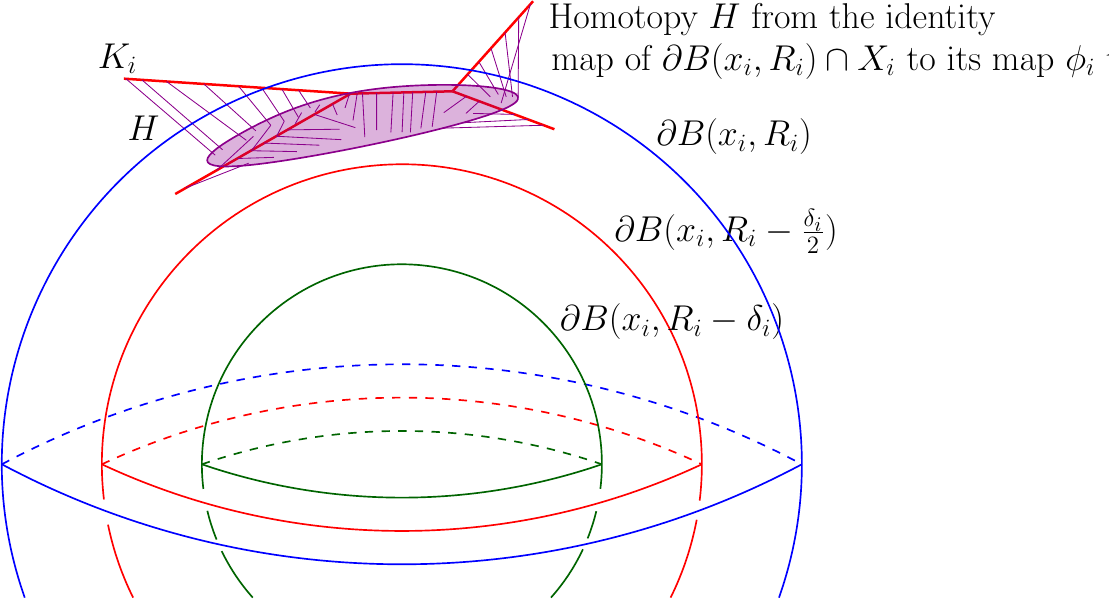}
    \caption{Points of $X_k$ in the outer annulus $A(x_i, R_i-\frac{\delta_i}{2}, R_i)$ are mapped using the homotopy $H=H_i$.}
    \label{figure:boundaryhomotopy}
\end{figure}

\begin{figure}
    \centering
    \includegraphics[scale=0.7]{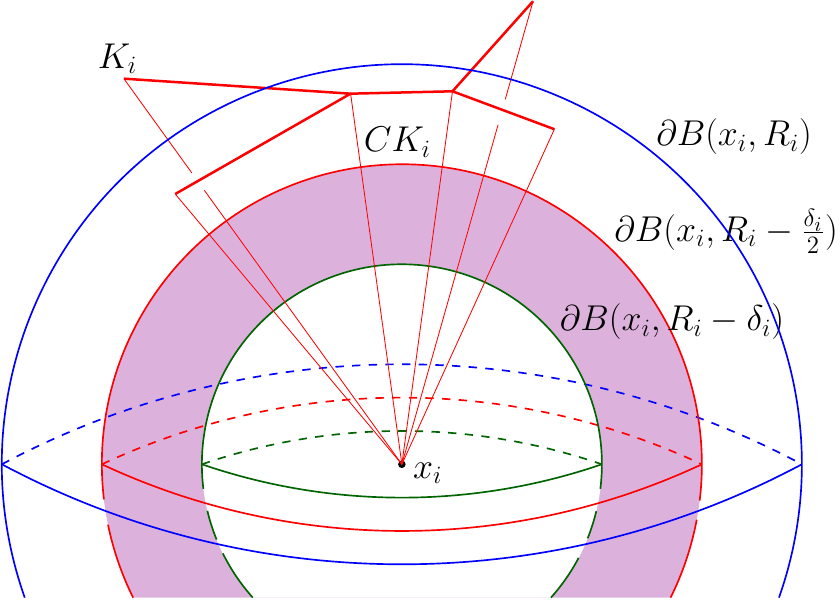}
    \caption{Points of $X_k$ in the inner annulus $A(x_i, R_i-\delta_i, R_i-\frac{\delta_i}{2})$ are mapped to the cone $CK_i$.}
    \label{figure:maptocone}
\end{figure}
Consider a covering of $\partial B(x_i,R_i)$ by a 
subcollection $C_{ki}$ of $\{\beta_j\}$ 
that realizes $\widetilde{\HC}_{m-1}(\partial B(x_i,R_i)\cap X_k)$.
Now we enlarge $X_k$ by adding to it the union of all balls in the collection $C_{ki}$ to $X_k$. Denote the result by $\tilde{X}_k$. (Obviously, this small addition will not change $\widetilde{\HC}_m$ with respect to the collection $\{\beta_j\}$.) By slightly enlarging the newly added balls (in $C_{ki}$) and choosing an arbitrarily small positive $\delta_i$, we can ensure
that there exists a homeomorphism 
$\alpha_i: (\tilde{X}_k\cap \partial B(x_i,R_i))\times [0,\delta_i ]\longrightarrow \tilde{X}_k\cap A(x_i, R_i-\delta, R_i)$ 
such that $\alpha_i(*,0)$ is the identity map of
$\tilde{X}_k\cap \partial B(x_i,R_i)$ and $\alpha_i(*,1)$ has its image in the geodesic sphere $\partial B(x_i,R_i-\delta_i)$ (see Figure~\ref{figure:homeomorphism}).

We are going to apply the induction assumption to $X_{ki}=\tilde{X}_k\cap \partial B(x_i,R_i)$. There
exists a $(\lceil m\rceil -2)$-dimensional simplicial complex $K_i$, a map $\theta_i:X_{ki}\longrightarrow K_i$ and a homotopy $H_i$ between the identity map and $\theta_i$ that satisfies the condition of Theorem 1.5. (See Figure~\ref{figure:emptyingball}).
Denote the cone over $K_i$ with the vertex at $x_i$ by $CK_i$. Observe that $K_i$ and, therefore,  $CK_i$ are contained in $B(x_i, R_i+ c_1(m-1)(\frac{2}{C_2(m)-C_1(m)})^{\frac{1}{m-1}}C_2(m){\frac{r_i}{A(m)^{\frac{m}{m-1}}}})$,
where $c_1(m-1)$ is the constant of Theorem 1.5 , and 
$C_1(m), C_2(m), A(m)$ are the constants of Theorem 2.1. (Here, we have combined Theorem 1.5(1) with the inequality in Theorem 2.1(3) for $m-1$ to estimate the distance
between any point of $K_i$ and $\partial B(x_i,R_i)$.)

Now we define $X_{k+1}$ first by removing from $X_k$ its intersections with the interiors
of all balls $B(x_i,R_i)$, but then adding 1) all cones $CK_i$ and 2) the images of all homotopies $H_i$. To construct $\phi_i$, we map each point $x$ of $X_k$ that is not inside
one of the balls $B(x_i,R_i)$ to itself now regarded a point of $X_{k+1}$). Then we map all points in each of the closed balls 
$B(x_i, R_i-\delta_i)$ to the center $x_i$.
We map
$\tilde X_k\cap \partial B(x_i, R_i-{\delta_i\over 2})$ to $K_i\subset CK_i\subset X_{k+1}$ first identifying it with $X_{ki}$ by applying $\alpha_i^{-1}$ and projecting it onto the first coordinate, and then applying $\theta_i$ (see Figure~\ref{figure:maptocone}). Now we rescale each interval $[0, {\delta_i\over 2}]$ to $[0,1]$, and map $\tilde X_k\cap A(x_i,R_i-{\delta_i\over 2}, R_i)$ first identifying it with $X_{ki}\times [0,1]$ using $\alpha_i^{-1}$ and then the rescaling of the second coordinate, and then applying the homotopy $H_i$ (see Figure~\ref{figure:boundaryhomotopy}).

It remains to map $\tilde X_k\cap A(x_i, r_i-\delta_i, r_i-{\delta_i\over 2})$ to $K_{k+1}$. This time we are going to map to $CK_i$. We will again apply $\alpha_i^{-1}$, and for each $x\in X_{ik}$ map 
$x\times [{\delta_i\over 2},\delta_i]$ to the straight line segment in $B$ that connects 
$\theta_i(x)\in K_i$ with the center $x_i$.

\subsection{End of the proof. Estimates for $c_1(m)$ and $c_2(m)$.}

Theorem 2.1(5) implies that the deletion of the part of $X_k$
in the union of all balls $B(x_i,R_i)$ decreases 
$\widetilde{\HC}_m(X_k)$ by at least 
$\Sigma_i({r_i\over A(m)})^m\geq {1\over (3C_2(m))^m}\widetilde{\HC}_m(X_k)$. (As before, $\widetilde{\HC_m}$ is calculated with respect to the nearly optimal covering of $X_k$ by 
balls $\{\beta_j\}$.)
As $\HC_m(X_k\setminus\bigcup_i B(x_i,r_i))\leq\widetilde{\HC}_m(X^k\setminus\bigcup_i B(x_i,r_i))$ and $\HC_m(X_k)\geq \widetilde{\HC}_m(X_k)-\tau$, we conclude that $\HC_m(X_k)$ would decrease by at least the same amount minus $\tau$.
The cones $CK_i$ have dimension $\lceil m\rceil -1$ and, therefore, 
their $m$-dimensional Hausdorff content is zero. The $(m-1)$-dimensional Hausdorff
content of the images of each
homotopy $H_i$ is bounded by $c_2(m-1)\HC_{m-1}(X_k\cap \partial B(x_i,R_i))\leq 
C_3(m)({r_i\over A(m)})^{m-1}$, where $C_3(m)=c_2(m-1){2\over C_2(m)-C_1(m)}{C_2(m)^{m-1}\over A(m)}$. Here $c_2(m-1)$ is the constant from Theorem 1.5(2), and we cannot choose it, but the other constants are from the text of Theorem 2.1 and can be chosen as we please, as long as they satisfy the constraints stated at the beginning of the text of Theorem 2.1. As for each set $A$ $\HC_m(A)\leq \HC_{m-1}^{m\over m-1}(A)$, the total $m$-dimensional Hausdorff content of $X_{k+1}$ does not exceed
$\HC_m(X_k)+\tau-(1-C_3(m)^{\frac{m}{m-1}})
\Sigma_i ({r_i\over A(m)})^m$. 
For our strategy to work, we need to ensure that the coefficient at $\Sigma_i ({r_i\over A(m)})^m$ is negative. If so, we can replace $\Sigma_i ({r_i\over A(m)})^m$ by its
lower bound ${1\over (3C_2(m))^m}\HC_m(X_k)$ provided by Theorem 2.1(5).
The condition that the coefficient at $\Sigma_i({r_i\over A(m)})^m$ is negative is equivalent to $C_3(m)<1$, and
can be rewritten as another lower bound for $A(m)$, namely,
$$ A(m)> c_2(m-1)C_2(m)^{m-1}{2\over C_2(m)-C_1(m)}.\ \ \ \ (*)$$
We would like the coefficient $1-C_3(m)^{\frac{m}{m-1}}$ to
be not merely negative but less than $-{1\over 2}$, which leads to a stronger inequality where $2$ in the right hand side of the previous lower bound for $A(m)$ is replaced by $2^{2-{1\over m}}$. We are going to drop $-{1\over m}$ which takes care of arbitrary small quantities such as $\tau$. If $A(m)$ would satisfy the inequality with the factor $4$ instead of $2$ on the right-hand side of inequality (*), then $\HC_m(X_{k+1})$ decreases by at least the factor $q=1-{1\over 2\cdot 3^mC_2(m)^m}$ compared to $\HC_m(X_k)$. As we observed above,
the sum of the geometric series $\Sigma_k q^k={1\over 1-q}= 2\times 3^mC_2(m)^m$ will enter the
estimate for $\HC_m(H(X\times [0,1])$, that is, it will be a factor in the formula
for $c_2(m)$, while a similar geometric series $\Sigma_k (q^{1\over m})^k= {1\over 1-q^{1\over m}}<2m\ 3^mC_2(m)^m$ will be a factor in the formula for $c_1(m)$.
We can choose $C_1(m)=1+{1\over m},\ C_2(m)=1+{2\over m}$, and 
$A(m)=30mc_2(m-1)$.
Note that $C_2(m)^m<\exp(2)<7.5$. As $c_2(m-1)\geq 1$, this $A(m)$ will satisfy the constraint in the text of Theorem 2.1
as well as the above inequality. 

Now we will evaluate the distance between $x\in X_0=X$ and $\phi_0(x)$ in the first step; the formula for $c_1(m)$ will then be obtained by multiplying
by the upper bound $20m\cdot 3^m$ for the sum 
${\Sigma_k \HC_m(X_k)^{1\over m}\over \HC_m(X)^{1\over m}}$ and adding $1$ to account for the very last step.

During the first step, the points do not move more than \par\noindent $R_i+c_1(m-1)(\frac{2}{C_2(m)-C_1(m)})^{\frac{1}{m-1}}C_2(m) \frac{1}{A(m)^{\frac{m}{m-1}}}r_i$ which is less than $[(1+\frac{2}{m})A(m)+c_1(m-1)\ (2m)^{\frac{1}{m-1}}\frac{1+\frac{2}{m}}{A(m)^{\frac{1}{m-1}}}]\frac{r_i}{A(m)}$. Substituting the formula
for $A(m)$, we see that the last expression does not exceed
$[90mc_2(m-1)+\frac{3c_1(m-1)}{(15\ c_2(m-1))^{\frac{1}{m-1}}}]\frac{r_i}{A(m)}$.
On the other hand, 
$r_i\leq A(m)\HC_m(X)^{1\over m}$. 
Combining the last two inequalities and using the argument in the previous paragraph, we see
that we can choose any
$$c_1(m)\geq 1+20m\cdot 3^m\cdot (90m\ c_2(m-1)+\frac{3c_1(m-1)}{(15\ c_2(m-1))^{\frac{1}{m-1}}}). \ \ \ \ \ \eqno(2.1)$$
Note that $c_1(m)\HC_m^{1\over m}(X)$
will also be an upper bound for the length of the trajectory of $H$ that will start as a broken line with exponentially decreasing segments followed by another broken line
obtained during the last stage and explained in the section below. Thus, our choice of $c_1(m)$ ensures that Theorem 1.5(1) holds.

Similarly, we can estimate $m$-dimensional Hausdorff content of the image
of the homotopy between $X_0$ and $\phi_0(X_0)$. This homotopy $H_0$ will consist of straight-line segments connecting points in $X_0\cap B(x_i,R_i)$ with their images. All these segments and, therefore, the image of the homotopy will be contained in $\bigcup_i B(x_i, 
R_i+c_1(m-1)(\frac{2}{C_2(m)-C_1(m)})^{\frac{1}{m-1}}C_2(m)\frac{1}{A(m)^{\frac{m}{m-1}}}r_i)$. Proceeding as above and using Lemma 2.2, we can majorize the sum over $i$ of the $m$th powers of these radii, $\tilde R_i$, by $(90mc_2(m-1)+\frac{3c_1(m-1)}{(15c_2(m-1))^{\frac{1}{m-1}}})^m\HC_m(X).$ Thus,
$\HC_m(H_0(X_0\times [0,1]))\leq \Sigma_i \tilde R_i^m\leq (90mc_2(m-1)+\frac{3c_1(m-1)}{(15c_2(m-1))^{\frac{1}{m-1}}})^m\HC_m(X)$.

Multiplying by $\Sigma q^k< 20\cdot 3^m$ and adding $1$ for the very last step, we see
that any
$$c_2(m)\geq 1+20\cdot 3^m (90mc_2(m-1)+\frac{3c_1(m-1)}{(15c_2(m-1))^{\frac{1}{m-1}}})^m \eqno(2.2)$$
would satisfy Theorem 1.5(2).
Recall that if $m\leq 1$, then $c_1(m)=c_2(m)=1$. Thus, we can take $c_1(m)$ to be equal to the right-hand side of (2.1), and $c_2(m)$ to be equal to the right-hand side of (2.2) and obtain the recurrence relations that can be used to determine $c_1(m)$ and $c_2(m)$ for all values of $m$. This completes the proof of Theorem 1.5 modulo Theorem 3.1 that will
be proven in the next section.

\par\noindent

\section{Boxing inequality with the constant that depends on the dimension of the ambient space}


\begin{theorem}
\label{thm:boxing_n}
For any positive $m$ and integer $n \ge m$ there are positive constants $c_{m,n},\ c'_{m,n}$ with the following property.
Let $X\subset B^n$ be a closed set in a $n$-dimensional Banach space $B$. Let $\delta > 0$ be a number\footnote{The summand $\delta$ is only there to account for the case $\hc_m(X)=0$ and can be ignored by the reader.}.
Then there exists a homotopy $H:X\times[0,1]\to\R^n$ such that
\begin{itemize}
\item $H_0=H(*,0)$ is the inclusion $X\subset B^n$,
\item $H_1(X)=H(X,1)$ is covered by an $(m-1)$-dimensional polyhedron in $\R^n$,
\item $\hc_m(H(X)) < c_{m,n}\hc_m(X)+\delta$,
\item $\Vert x-H_1(x)\Vert< c'_{m,n}(\hc_m(X))^{\frac{1}{m}}+\delta$. Moreover, the length of the trajectory of $x$ under $H$ does not exceed $c'_{m,n}(\hc_m(X))^{\frac{1}{m}}+\delta$ for every $x\in X$.
\end{itemize}
\end{theorem}

{\bf Remarks.} (1) Each trajectory of $H$ will consist of finitely many straight line segments connecting points $x_ix_{i+1}$ for $i=1,2\ldots, N(x)-1$ with $x_1=x=H(x,0)$
and $x_{N(x)}=H_1(x)$. The length of the trajectory is equal to $\Sigma_i\Vert x_{i+1}-x_i\Vert$.
\par\noindent
(2) Given a pair of $n$-dimensional Banach spaces $B_1, B_2$ there exists a linear isomorphism $A:B_1\longrightarrow B_2$ such that $\Vert A\Vert\Vert A^{-1}\Vert\leq n$.
(This is a well-known corollary of F. John's ellipsoid theorem.)
Therefore, Theorem 3.1 follows immediately from its particular case, where $B^n=\ell^n_\infty$. Indeed, given any Banach space $B^n$ consider a linear isomorphism
$A:B^n\longrightarrow \ell^n_\infty$ such that $\Vert A\Vert\Vert A^{-1}\Vert\leq n$. Observe that the image under $A$ of a metric ball of radius $r$ in $B^n$ is contained in a metric ball of radius $\Vert A\Vert r$ in $\ell^n_\infty$. So, $\HC_m(A(X))\leq \Vert A\Vert ^m\HC_m(X)$. Apply Theorem 3.1 to $A(X)$
in $\ell^n_\infty$, and then map the filling back to $B^n$ using $A^{-1}$. As a result, we will get
the claim of Theorem 3.1 for $X$ in $B^n$ but with extra factors $\Vert A\Vert\Vert A^{-1}\Vert\leq n$ in $c'_{m,n}$ and $\Vert A\Vert^m\Vert A^{-1}\Vert^m\leq n^m\leq n^n$ in
$c_{m,n}$.

Therefore, in the following, we are going to present the proof of the theorem
only for this particular case. (Recall, that metric balls in $\ell^n_\infty$ are cubes
with sides parallel to coordinate axes.)
\par\noindent
(3) If $B^n=\ell^n_\infty$ and $X$ is piecewise smooth or PL, the homotopy $H$ that
will be constructed in our proof will also be piecewise smooth or, correspondingly, PL.
\par\noindent
(4) Our upper bounds for constants that depend on $m$ and $n$ will increase
with $n$ and $m$. As $m\leq n$, this automatically implies that they can be replaced by $c_n$ and $c'_n$ that depend only on $n$. 
However, we prefer to keep the notations as in the theorem to emphasize that the optimal values of these constants probably also depend on $m$.

\subsection{Plan of the proof.} 
Let us first give an informal explanation of how the proof of the theorem works. We can assume that $B^n=\ell^n_\infty$.

We first cover $X$ by dyadic cubes so that the total $\hc_m$ of the covering is comparable to $\hc_m(X)$, each cube contains a relatively small ``amount'' (measured in $\hc_m$) of $X$, and the neighboring cubes of the covering are comparable in size. The homotopy $H$ will be constructed within this covering.

We then triangulate the dyadic covering, to get a triangulation where each simplex is bilipschitz close to a round ball (or an equilateral simplex), the complexity of neighborhoods and sizes of adjacent simplices are bounded, and each simplex still contains a relatively small ``amount'' of $X$. This is very similar to the concept of QC complexes in [Y].

Finally, we use H. Federer and W. Fleming approach ([FF]) to iteratively ``push'' $X$ first from the top-dimensional simplices of the triangulation, then its codimension $1$ skeleton, and so on, until we arrive at its $(\lceil m\rceil-1)$th skeleton. The properties of the triangulation ensure that at each step there is not too much $X$ in every simplex.


\subsection{On the dyadic cubes.}
A \emph{dyadic face} is a face (of any dimension) of a dyadic cube.
When we talk of \emph{larger} or \emph{smaller} dyadic cubes (or faces) or refer to their \emph{size} we mean their diameters in $\ell^\infty$ metric.
For any two dyadic cubes, either their interiors are disjoint, or one of the cubes is contained in the other.

All dyadic cubes and their faces are assumed to be closed unless otherwise specified.

\begin{lemma}
\label{lemma:dyadic_hausdorff}
An $\ell^n_\infty$ ball of diameter (=size) $s$ can be covered by $4^n$ dyadic cubes of size at most $2s$.
\end{lemma}
\begin{proof}
Let $S_i$ be the projection of the $\ell^\infty$ ball onto the $i$th coordinate axis.

If there are at least two vertices of dyadic line segments of length $d$ in $S_i$, then at least one of them is also a vertex of a dyadic line segment of length $2d$. So, there is the smallest number $d_i$ such that there is exactly one vertex $x_i\in S_i$ of a dyadic line segment of length $d_i$ in $S_i$. Then there are at least two vertices of dyadic line segments of length $d_i/2$ in $S_i$, which means that $d_i/2 \leq |S_i| = s$ and $d_i \leq 2s$. On the other hand, $S_i\subset (x_i-d_i, x_i+d_i)$ and so $d_i > s/2$ and $S_i$ can be covered either by $2$ dyadic line segments of length $d_i$ or by $4$ of length $d_i/2$.

Suppose all $d_i$ are equal. Then the $\ell^\infty$ ball can be covered by $2^n$ dyadic cubes of size $d_1\leq 2s$.

If $d_i>d_j$ for some $i$ and $j$, then from $s/2 < d_i,d_j \leq 2s$ we get $d_j=d_i/2\leq s$.
So, each of the projections of the ball can be covered by at most $4$ dyadic line segments of length $d_j\leq s$.
Then the ball itself can be covered by $4^n$ dyadic cubes of size $d_j\leq s$
\end{proof}

\subsection{Construction and properties of the dyadic covering $X\subset \bigcup Q_\ell$.}

Denote by $\hc_m^{d}$ the Hausdorff content taken with the restriction that every covering cube must be dyadic.
By Lemma~\ref{lemma:dyadic_hausdorff}, for any bounded subset $U\subset \R^n$ we have $\hc_m(U)\leq \hc_m^{d}(U) \leq 4^n\cdot 2^m\cdot\hc_m(U)$.

For a dyadic cube (or any bounded set) $Q$ define its \emph{density} as $\frac{\hc_m^{d}(X\cap Q)}{\hc_m^{d}(Q)}$. Let $0<\epsilon < 1$ be a number that we will choose later, and $\delta>0$ the same as in Theorem \ref{thm:boxing_n}.
We will also need an arbitrarily small positive parameter $\delta'<<\delta$ that
can be easily chosen later to demonstrate Theorem \ref{thm:boxing_n} when $\hc_m(X)=0$.

In this subsection, we construct a \emph{good} covering $X\subset \bigcup Q_\ell$ by dyadic cubes $\{Q_\ell\}$ with the following properties:

\begin{itemize}

\item[(i)] The density of every $Q_\ell$ is less than $\epsilon$.

\item[(ii)] If $\hc_m(X)\not= 0$, then $\hc_m(\bigcup Q_\ell)\le c_{m,n,\epsilon} \hc_m(X)$ for some constant $c_{m,n,\epsilon}$. If $\hc_m(X)=0$, then
$\hc_m(\bigcup Q_\ell)<\delta'$.

\item[(iii)] If $\hc_m(X)\not= 0$, then the size of each $Q_\ell$ is at most $c_{m,n,\epsilon}(\hc_m(X))^{\frac{1}{m}}$. If $\hc_m(X)=0$, then the size of
each $Q_\ell<\delta'$.

\item[(iv)] For any cube $Q_\ell$ of size $s$ the $\ell^\infty$-distance between $Q_\ell$ and any point $P$ of a smaller cube in $\{Q_\ell\}$ of size at most $s/4$ is at least $s/2$.
\end{itemize}

Let $\bigcup Q''_i$ be a finite covering of $X$ by dyadic cubes $Q''_i$ such that if $\HC_m(X)\not= 0$, then $\sum_i\hc_m^{d}(Q''_i)< 1.1 \cdot \hc_m^{d}(X)$. If $\HC_m(X)=0$, we choose $\bigcup Q''_i$ so that $\sum_i\hc_m^{d}(Q''_i)<\delta'$ for an appropriate positive $\delta'<<\delta$, and sizes of all cubes $Q_\ell<\delta'$.
By incorporating $4^n\cdot 2^m$ and $1.1$ into the constant $c_{m,n}$ in the statement of the theorem we may now assume that $X=\bigcup Q''_i$.

For any $Q''_i$, let $Q'_i$ be the smallest dyadic cube containing $Q''_i$ and such that both $Q'_i$ and all larger dyadic cubes containing $Q''_i$ have a density lower than $\epsilon$. A cube $Q'_i$ that satisfies this condition exists because, as the size of dyadic cubes containing $Q''_i$ goes to infinity, their density goes to $0$, and the density of $Q''_i$ is exactly $1$.

Remove from the set of cubes $\{Q'_i\}$ the ones contained in a larger cube from the set. Renumber the remaining cubes $Q'_1, Q'_2,\ldots $. They cover the whole $X$.
For this covering, the following inequalities hold:

\begin{multline}
\label{eq:cover_size}
\sum_j \hc_m(Q'_j) = \sum_j \hc_m^{d}(Q'_j) \overset{(1)}{<} \sum_j \frac{2^m}{\epsilon} \cdot \hc_m^{d}(X\cap Q'_j) \overset{(2)}{<} \\
\overset{(2)}{<} 6^n\cdot \frac{2^m}{\epsilon}\cdot \sum_i \hc_m^{d}(Q''_i) < 1.1\cdot 6^n\cdot \frac{2^m}{\epsilon} \cdot \hc_m^{d}(X) \leq 1.1\cdot 6^n\cdot \frac{4^{m+n}}{\epsilon} \cdot \hc_m(X).
\end{multline}

Let us explain some of these inequalities:
\begin{itemize}
\item[(1)] By the definition of $Q'_j$, there is a dyadic cube $Q\subset Q'_j$ such that its size is equal to half the size of $Q'_j$, and its density is greater than $\epsilon$. So, $\frac{\hc_m^{d}(X\cap Q)}{\hc_m^{d}(Q)}>\epsilon$. This implies that

\[
\hc_m^{d}(Q'_j) = 2^m\hc_m^{d}(Q) < \frac{2^m}{\epsilon}\hc_m^{d}(X\cap Q)\leq \frac{2^m}{\epsilon}\hc_m^{d}(X\cap Q'_j).
\]

\item[(2)] Recall that $X=\bigcup Q''_i$. So, $\hc_m^{d}(X\cap Q'_j)$ is less than or equal to $\sum \hc_m^{d}(Q''_i)$ where the sum is taken over those $Q''_i$ that intersect $Q'_j$ by a subset of dimension at least $\lceil m\rceil$. If this happens, then $Q''_i$ is smaller than $Q'_j$, as otherwise the density of $Q'_j$ would be $1$. Therefore, $Q'_j$ contains a face of $Q''_i$ of dimension at least $\lceil m\rceil$. Each face of $Q''_i$ of dimension $k$ can be contained in at most $2^{n-k}$ of the cubes $Q'_j$ because the interiors of those cubes are disjoint. 
We see that each $Q''_i$ contributes to $\hc_m^{d}(X\cap Q'_j)$
for at most $C$ different $Q'_j$, where $C$ is the sum of $2^{n-k}$ over all $k\ge \lceil m\rceil$ and all faces of $Q''_i$ of dimension $k$. We can conclude that $C<6^n$ because the cube has less than $3^n$ faces of positive dimension and $2^{n-k}\leq 2^n$. 

\end{itemize}

Now we modify the covering $X\subset\bigcup Q'_j$ further. For a cube $Q'_j$ its \emph{collar} consists of several \emph{layers} defined inductively; see Fig.~\ref{figure:collar}.
The first layer is the cube $Q'_j$ itself. The $(k+1)$th layer consists of all dyadic cubes that are $2^k$ times smaller than $Q'_j$ and that touch the $k$th layer from the outside.
The last layer of the collar of $Q'_j$ consists of cubes of size equal to the size of the smallest cube in the set $\{Q'_j\}$.

\begin{figure}[ht]
\center
\includegraphics{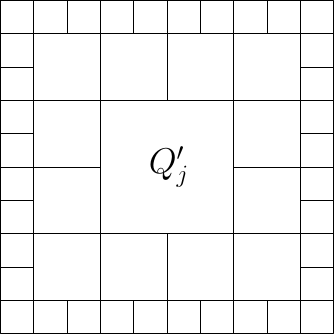}
\caption{The first $3$ layers of the collar of $Q'_j$.}
\label{figure:collar}
\end{figure}

Take the set of all the cubes in the collars of all $Q'_j$. Remove from this set any cube covered by another.
Remove any remaining cube if it contains no \emph{interior} points of any of $Q''_i$.
What remains is a covering $\bigcup Q_\ell\supset X$.
Let us check that this covering satisfies the properties of a good covering:

\begin{itemize}

\item[(i)] The density of every $Q_\ell$ is less than $\epsilon$. Indeed, $Q_\ell$ contains an interior point of $Q''_i$ for some $i$. Let $Q'_j$ be the unique cube in $\{Q'_j\}$ containing $Q''_i$. Then $Q'_j\subseteq Q_\ell$. This implies that the density of $Q_\ell$ is less than $\epsilon$, because all dyadic cubes containing $Q'_j$
have density less than $\epsilon$ by the construction of the covering $\bigcup Q'_j$.

\item[(ii)] $\hc_m(\bigcup Q_\ell)\le 3^m \sum_j \hc_m(Q'_j)< 3^m\cdot 1.1\cdot 6^n\cdot \frac{4^{m+n}}{\epsilon} \cdot \hc_m(X)$. Indeed, for any $j$ the collar of $Q'_j$ can be covered by a copy of $Q'_j$ scaled by the factor of $3$.

\item[(iii)] The size of each $Q_\ell$ is at most $c_{m,n,\epsilon}(\hc_m(X))^{\frac{1}{m}}$. Indeed, it is at most the size of the largest cube in $\{Q'_j\}$ that can be bounded in terms of $\varepsilon$ and $(\hc_m(X))^{\frac{1}{m}}$ because the density of each $Q'_j$ is at least $2^{-m}\varepsilon$ by definition.

\item[(iv)] For any cube $Q_\ell$ of size $s$ the $\ell^\infty$-distance between $Q_\ell$ and any point $P$ of a smaller cube in $\{Q_\ell\}$ of size at most $s/4$ is at least $s/2$. 
\end{itemize}

\begin{lemma}
\label{lemma:collar_distance}
Property (iv) holds.
\end{lemma}
\begin{proof}
The cube $Q_\ell$ comes from the collar of $Q'_j$ for some $j$. If it is in the last layer of the collar, then there are no smaller cubes in the covering and the property holds automatically.

If it is not in the last layer, then its $s/2$-neighborhood is covered by the cubes of the collar of size at least $s/2$. If one of these cubes is not in $\{Q_\ell\}$ then either it was replaced by a larger cube that contained it, or it did not contain interior points of $Q''_i$ and was removed. In the latter case, observe that none of smaller cubes inside of the removed cube can contain an interior point of $Q''_i$. Therefore, if such a smaller cube were present in a collar of another cube, it would have been removed, too, and  cannot be present in the covering $\{Q_\ell\}$.
\end{proof}

\subsection{Triangulating the covering $X\subset \bigcup Q_\ell$.}

It is intuitively clear that the covering $\bigcup Q_\ell$ can be triangulated so that for some constant $C_n$ independent of $X$ and $\{Q_\ell\}$ we have the following:

\begin{itemize}
\item[(i)] Every simplex in the triangulation is $C_n$-bilipschitz close to an equilateral simplex.
\item[(ii)] $\frac{1}{C_n} < \frac{\hc_m(\sigma_1)}{\hc_m(\sigma_2)} < C_n$ for any two simplices $\sigma_1\cap\sigma_2\neq\emptyset$.
\item[(iii)] Every simplex is contained in less than $C_n$ of simplices of higher dimension.
\item[(iv)] The density of every simplex is less than $\epsilon C_n$. 
\end{itemize}

For the reader's convenience, let us prove that such a triangulation exists. (This triangulation is in no way unique and there might be many different constructions).

For a closed face (of any dimension) $F$ of size $s$ of a cube in $\{Q_\ell\}$ denote by $N(F)$ the set of \emph{open} faces of size at most $s$ of all cubes $Q_\ell$ intersecting $F$. Clearly, every face in $N(F)$ is contained in $F$ and, in particular, its dimension is at most $\dim\ F$. If $F_1\in N(F)$, then either $F_1$ is a face of $F$, or the size of $F_1$ is $s/2$. In the latter case $N(F_1)$ is just the set of subfaces of $F_1$, which we call $N(F_1)$ \emph{simple}.

Let us say that $N(F_1)$ is \emph{isomorphic} to $N(F_2)$ if there is a translation and scaling which induces a bijection between the sets. Clearly, there are a finite number of possible isomorphism classes and these classes are independent of $\{Q_\ell\}$ or $X$.

We now define a partial order on the isomorphism classes in such a way that $N(F_1)<N(F_2)$ if $\dim F_1< \dim F_2$. If $\dim F_1 = \dim F_2$, we put $N(F_1)<N(F_2)$ whenever $N(F_1)$ is simple and $N(F_2)$ is not. From what was said above, we see that if $F_1\in N(F)$ then $N(F_1) < N(F)$.

Let us now triangulate each face of $\{Q_\ell\}$ following the partial order, that is, triangulating $F_1$ before $F_2$ whenever $N(F_1)<N(F_2)$. By the time we triangulate $F$, we have already triangulated every face in $N(F)$. We can extend this triangulation to a triangulation of $F$, possibly adding vertices in the interior of $F$. We triangulate every face $F'$ with $N(F')$ isomorphic to $N(F)$ \emph{in the same way} - meaning that the translation and scaling inducing the isomorphism between $N(F')$ and $N(F)$ also induces a bijection between the triangulations of $F'$ and $F$. It is possible to do so because $F\notin N(F')$ and $F'\notin N(F)$, and so there are no ``circular dependencies''.

Property (i) holds because there are only a finite number of congruence classes of possible simplices. Property (ii) holds because there is only a finite number of congruence classes of possible pairs of intersecting simplices. Property (iii) holds because there are only a finite number of possible combinatorial classes of neighborhoods of simplices. Property (iv) holds because the ratio between $\hc_m$ of a cubical face and a simplex in its triangulation can only take a finite possible number of values.

\subsection{Construction of the homotopy $H$.}

The homotopy $H$ is going to be a series of elementary projections from the simplices of the triangulation using the following lemma:

\begin{lemma}
\label{lemma:projection}
There exist positive constants $C_{scale}$ and $C_{max density}$ depending only on $m$ and $n$, such that for every $m < k \leq n$, a $k$-dimensional simplex $S\subset \R^n$ which is $C_n$-bilipschitz close to an equilateral simplex, and $X\subset \inte S$ with
$\hc_m(X) < C_{max density}\cdot\hc_m(S)$ there is a point $O\in \inte S\setminus X$ such that $\hc_m(\phi_O(X)) \le C_{scale}\cdot\hc_m(X)$, where $\phi_O:S\setminus O\to\partial S$ is the radial projection in $S$ with the center $O$.
\end{lemma}

See a proof in [Gu13], Lemma 7.2. (See also Lemma 2.5 in [Y] for a similar result). The difference from our statement is that our simplex is not equilateral, but only $C_n$-bilipschitz close to an equilateral, which worsens the constants $C_{scale}$ and $C_{max density}$ but does not affect the lemma otherwise.

We now use the lemma to project $X$ from each simplex of dimension at least $\lceil m\rceil$ of the triangulation constructed in the previous subsection. After projecting from the $\lceil m\rceil$-dimensional simplices $X$ will lie in the $(\lceil m\rceil-1)$th skeleton of the triangulation.

We first project from the top-dimensional simplices, then from the simplices of codimension $1$, then $2$, and so on. To apply the lemma, we only need to ensure that at each step the density of $X$ in the simplex where we project is not greater than $C_{max density}$. We can achieve this by choosing $\varepsilon$ sufficiently small, because the original density in each simplex is at most $\varepsilon C_n$ and because every simplex is contained in less than $C_n$ of simplices of higher dimension and hence we can control how much of $X$ arrives in a simplex (from the simplices of higher dimension) after each round of projections.

Because the image of the whole homotopy $H$ lies in $\bigcup Q_\ell$, its Hausdorff content is bounded by $\hc_m(\bigcup Q_\ell)$ which itself is bounded in terms of $\hc_m(X)$. The length of the trajectory of each point $x\in X$ is bounded because $x$ is only projected $(n-\lceil m\rceil+1)$ times, each time moving by not more than the size of the largest cube in $\{ Q_\ell\}$ which in turn is bounded in terms of $\varepsilon$ and $(\hc_m(X))^{\frac{1}{m}}$.

\forget
First, we cover $X$ by a set of dyadic cubes  $\bigcup Q''_i$ so that $\sum_i \hc_m(Q''_i)$ is not too large compared to $\hc_m(X)$. This allows us to identify $X$ with $\bigcup Q''_i$.

Then we construct a bigger cover $\bigcup Q'_j$ of $X=\bigcup Q''_i$, again by dyadic cubes.
This time we choose each cube $Q'_j$ in such a way that $\hc_m(X\cap Q'_j)$ is neither too large nor too small compared to $\hc_m(Q'_j)$. The second property ensures that if we manage to construct the homotopy $H$ inside $\bigcup Q'_j$, then the size of image of this homotopy $\hc_m(H)$ will not be too large compared to $\hc_m(X)$. On the other hand, the first property gives us enough space to actually construct the homotopy.

At this point we would like to start radially projecting $X$ from the faces of the cover $\bigcup Q'_j$. 
We start with the smallest cube in $\bigcup Q'_j$ and project from a suitable point in its interior, choosing the point in such a way that the projection is not increasing $\hc_m(X)$ too much, using the fact that there is not too much of $X$ in the cube. Then we project from the boundary faces of the cube to their respective boundaries and so on, until we push $X$ to the $(m-1)$-dimensional skeleton of the smallest cube. In this manner we go through all the cubes in $\bigcup Q'_j$ in order of the increasing size and project first from the cubes and then from their faces of lower and lower dimension. We skip the faces which are covered by strictly larger cubes and do not project from them.

\begin{figure}[ht]
\center
\includegraphics{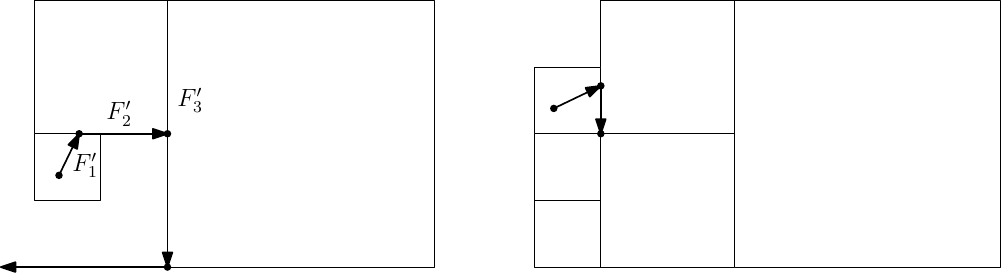}
\caption{On the left the point is projected from a $2$-dimensional face $F'_1$ to a $1$-dimensional face $F'_2$ which is twice as big. Then it is projected to a $1$-dimensional face $F'_3$ which is again, twice as big as $F'_2$. In principle, this spiral can continue for an uncontrollable number of steps.
On the right the point reaches a face which is twice as big as the one it starts at, but cannot reach an even bigger face on the right.}
\label{figure:paths}
\end{figure}

The issue with this approach is that a point in $X$ can potentially be  projected an uncontrollable number of times. For instance, it can be projected from some face $F'_1$ of a cube $Q'_1$ to a point in $\partial F'_1$ which also belongs to a face $F'_2$ of a bigger cube $Q'_2$. Then it is projected to even a bigger face $F'_3$ of some cube $Q'_3$, and so on. The length of the (finite) sequence $F'_1, F'_2, \ldots$  cannot be bounded in terms of $m$ and $n$ since the dimensions of the faces in it are not necessarily decreasing and might even increase, see Fig.~\ref{figure:paths} (left).

As the result, by the moment we decide to process a face $F$ of size $s$, it could have
absorbed ``material" from uncontrollably many cubes with sizes $\frac{s}{2^k}$, where
$k$ can be arbitrarily large. These cubes will reach $F$ via many spirals as on Fig. \ref{figure:paths} (left). On the first glance the problem might seem not too serious,
as (a) we immediately realize that all these small cubes will be at the distance
at most $s$ from $F$; and (b) the $m$-dimensional content of the intersection of $X$
with each small cube of size $\frac{s}{2^{k}}$ is at most $\frac{s^m}{2^{km}}\longrightarrow 0$, as $k\longrightarrow\infty$. However, the number of dyadic cubes of size $\frac{s}{2^k}$ that can potentially reach $F$ behaves as $\sim 2^{kn}$,
which grows much faster than $2^{km}$. In addition, the Hausdorff content might increase
by a constant factor during a radial projection. When the number of projections cannot be bounded in terms of the quantities that we control, we cannot claim that
the contribution of even an individual cube of an uncontrollably small size is $o(1)$.

To combat this issue we construct the third and final cover $\bigcup Q_\ell$. To do so we add a ``collar'', see Fig.~\ref{figure:collar}, to each cube in $\bigcup Q'_j$ and then get rid of all the cubes in the union of the collars which are covered by a bigger cube. This new cover $\bigcup Q_\ell$ has the following property that the cover $\bigcup Q'$ lacked: Any two cubes of size $s$ and $4s$ of the covering $\bigcup Q_\ell$ are always separated by least one cube of size $\geq 2s$. 
 This allows us to ensure that each point of $X$ is going to be projected a number of times that we can majorize in terms of $m$ and $n$. Indeed, a point will not be able to reach a cube more than $2$ times larger than the cube it started in, see Fig.~\ref{figure:paths} (right).

Finally, note that, a radial projection in a face can cause a discontinuity at a subset of this face of codimension at least $2$ (see the figure at the end of section 1.4. Recall, that the main reason for this phenomenon is that we use cubes with different side lengths.) To remedy this problem we use ``insertions" as sketched in section 1.4 and explained in details below (see Figures~\ref{figure:homfix} and ~\ref{figure:squished}).

\subsection{On the dyadic cubes.}
A \emph{dyadic face} is a face (of any dimension) of a dyadic cube.
When we talk of \emph{larger} or \emph{smaller} dyadic cubes (or faces) or refer to their \emph{size} we mean their diameters in $\ell^\infty$ metric.
For any two dyadic cubes either their interiors are disjoint or one of the cubes is contained in the other.

All the dyadic cubes and their faces are assumed to be closed.

\begin{lemma}
\label{lemma:dyadic}
Let $Q$ be a dyadic face. Suppose some other dyadic face $Q'$ contains an interior point of $Q$ but does not cover $Q$.
Then $Q'$ is strictly smaller than $Q$.
\end{lemma}
\begin{proof}
Let $\widetilde{Q}'$ be a dyadic cube whose face $Q'$ is. Suppose that $Q\subset \widetilde{Q}'$. Then $Q\subset Q'$, because $Q'$ contains an interior point of $Q$, contradiction. So, $Q$ is not covered by $\widetilde{Q}'$.

Because $\widetilde{Q}'$ contains an interior point of $Q$, there is a dyadic cube $\widetilde{Q}$ having $Q$ as its face and intersecting the interior of $\widetilde{Q}'$.
So, either $\widetilde{Q}\subseteq\widetilde{Q}'$, or $\widetilde{Q}'\subset \widetilde{Q}$. But the former was already proved to be impossible, and the latter case means that $Q'$ is strictly smaller than $Q$.

\end{proof}

\begin{lemma}
\label{lemma:dyadic_hausdorff}
An $\ell^n_\infty$ ball of diameter (=size) $s$ can be covered by $4^n$ dyadic cubes of size at most $2s$.
\end{lemma}
\begin{proof}
Let $S_i$ be the projection of the $\ell^\infty$ ball  onto the $i$th coordinate axis.

If there are at least two vertices of dyadic line segments of length $d$ in $S_i$, then at least one of them is also a vertex of a dyadic line segment of length $2d$. So, there is the smallest number $d_i$ such that there is exactly one vertex $x_i\in S_i$ of a dyadic line segment of length $d_i$ in $S_i$. Then there are at least two vertices of dyadic line segments of length $d_i/2$ in $S_i$ meaning that $d_i/2 \leq |S_i| = s$ and $d_i \leq 2s$. On the other hand, $S_i\subset (x_i-d_i, x_i+d_i)$ and so $d_i > s/2$ and $S_i$ can be covered either by $2$ dyadic line segments of length $d_i$ or by $4$ of length $d_i/2$.

Suppose all $d_i$ are equal. Then the $\ell^\infty$ ball can be covered by $2^n$ dyadic cubes of size $d_1\leq 2s$.

If $d_i>d_j$ for some $i$ and $j$, then from $s/2 < d_i,d_j \leq 2s$ we get that $d_j=d_i/2\leq s$.
So, each of the projections of the ball can be covered by at most $4$ dyadic line segments of length $d_j\leq s$.
Then the ball itself can be covered by $4^n$ dyadic cubes of size $d_j\leq s$
\end{proof}

\subsection{Construction and properties of the cover $X\subset \{Q_\ell\}$.}

Denote by $\hc_m^{d}$ the Hausdorff content taken with the restriction that every covering cube has to be dyadic.
By Lemma~\ref{lemma:dyadic_hausdorff}, for any bounded subset $U\subset \R^n$ we have $\hc_m(U)\leq \hc_m^{d}(U) \leq 4^n\cdot 2^m\cdot\hc_m(U)$.

Let $\bigcup Q''_i$ be a finite cover of $X$ by dyadic cubes $Q''_i$ such that if $\HC_m(X)\not= 0$, then $\sum_i\hc_m^{d}(Q''_i)< 1.1 \cdot \hc_m^{d}(X)$. If $\HC_m(X)=0$, we choose $\bigcup Q''_i$ so that $\sum_i\hc_m^{d}(Q''_i)<\delta'$ for an appropriate positive $\delta'<<\delta$.
By incorporating $4^n\cdot 2^m$ and $1.1$ into the constants $c_{m,n}$ and $c'_{m,n}$ we may now assume that $X=\bigcup Q''_i$.

For a dyadic cube $Q$ define its \emph{density} as $\frac{\hc_m^{d}(X\cap Q)}{\hc_m^{d}(Q)}$.

Let $0<\epsilon < 1$ be a number which we will choose later. For any $Q''_i$, let $Q'_i$ be the smallest dyadic cube containing $Q''_i$ and such that both $Q'_i$ and all larger dyadic cubes containing $Q''_i$ have density less than $\epsilon$. Such a cube $Q'_i$ exists because, as the size of dyadic cubes containing $Q''_i$ goes to infinity, their density goes to $0$, and the density of $Q''_i$ is exactly $1$.

Remove from the set of cubes $Q'_i$ the ones contained in a larger cube from the set. Renumber the remaining cubes $Q'_1, Q'_2,\ldots $. They cover the whole $X$.
For this cover the following inequalities hold:

\begin{multline}
\label{eq:cover_size}
\sum_j \hc_m(Q'_j) = \sum_j \hc_m^{d}(Q'_j) \overset{(1)}{<} \sum_j \frac{2^m}{\epsilon} \cdot \hc_m^{d}(X\cap Q'_j) \overset{(2)}{<} \\
\overset{(2)}{<} 6^n\cdot \frac{2^m}{\epsilon}\cdot \sum_i \hc_m^{d}(Q''_i) < 1.1\cdot 6^n\cdot \frac{2^m}{\epsilon} \cdot \hc_m^{d}(X) \leq 1.1\cdot 6^n\cdot \frac{4^{m+n}}{\epsilon} \cdot \hc_m(X).
\end{multline}

Let us explain some of these inequalities:
\begin{itemize}
\item[(1)] By the definition of $Q'_j$, there is a dyadic cube $Q\subset Q'_j$ such that its size is equal to the half the size of $Q'_j$, and its density is greater than $\epsilon$. So, $\frac{\hc_m^{d}(X\cap Q)}{\hc_m^{d}(Q)}>\epsilon$. This implies that

\[
\hc_m^{d}(Q'_j) = 2^m\hc_m^{d}(Q) < \frac{2^m}{\epsilon}\hc_m^{d}(X\cap Q)\leq \frac{2^m}{\epsilon}\hc_m^{d}(X\cap Q'_j).
\]

\item[(2)] Recall that $X=\bigcup Q''_i$. So, $\hc_m^{d}(X\cap Q'_j)$ is less than or equal to $\sum \hc_m^{d}(Q''_i)$ where the sum is taken over those $Q''_i$ which intersect $Q'_j$ by a subset of dimension at least $\lceil m\rceil$. If this happens, then $Q''_i$ is smaller than $Q'_j$, as otherwise the density of $Q'_j$ would be $1$. Therefore, $Q'_j$ contains a face of $Q''_i$ of dimension at least $\lceil m\rceil$. Each face of $Q''_i$ of dimension $k$ can be contained in at most $2^{n-k}$ of the cubes $Q'_j$ because the interiors of those cubes are disjoint. 
We see that each $Q''_i$ contributes to $\hc_m^{d}(X\cap Q'_j)$
for at most $C$ different $Q'_j$, where $C$ is the sum of $2^{n-k}$ over all $k\ge \lceil m\rceil$ and all faces of $Q''_i$ of dimension $k$. We can conclude that $C<6^n$ because the cube has less than $3^n$ faces of all dimensions $\geq 1$ 
and $2^{n-k}\leq 2^n$. 

\end{itemize}

Let us now modify the cover $\bigcup Q'_j$ further. For a cube $Q'_j$ its \emph{collar} consists of several \emph{layers} defined inductively, see Fig.~\ref{figure:collar}.
The first layer is the cube $Q'_j$ itself. The $(k+1)$th layer consist of all dyadic cubes which are $2^k$ times smaller than $Q'_j$ and which touch the $k$th layer from the outside.
The last layer in the collar of $Q'_j$ consists of the cubes with the size equal to the size of the smallest cube in the cover $\bigcup Q'_j$.

\begin{figure}[ht]
\center
\includegraphics{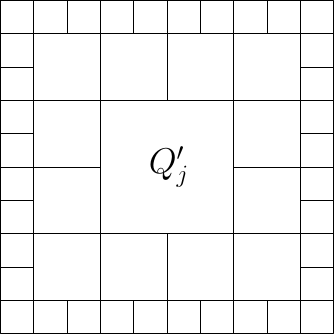}
\caption{The first $3$ layers of the collar of $Q'_j$.}
\label{figure:collar}
\end{figure}

Take the set of all the cubes in the collars of all $Q'_j$. Remove from this set any cube covered by another.
Remove any remaining cube if it contains no interior points of any of $Q''_i$.
What remains is a cover $\bigcup Q_\ell\supset X$.
This cover has the following properties:

\begin{itemize}

\item[(i)] The density of every $Q_\ell$ is less than $\epsilon$. Indeed, $Q_\ell$ contains an interior point of $Q''_i$ for some $i$. Let $Q'_j$ be the unique cube in $\{Q'_j\}$ containing $Q''_i$. Then $Q'_j\subseteq Q_\ell$. This implies that the density of $Q_\ell$ is less than $\epsilon$, because all dyadic cubes containing $Q'_j$
have density less than $\epsilon$ by the construction of the cover $\bigcup Q'_j$.

\item[(ii)] $\hc_m(\bigcup Q_\ell)\le 3^m \sum_j \hc_m(Q'_j)$. Indeed, for any $j$ the collar of $Q'_j$ can be covered by a copy of $Q'_j$ scaled by the factor of $3$.

\item[(iii)] For any cube $Q_\ell$ of size $s$ the $\ell^\infty$-distance between $Q_\ell$ and any point $P$ of a smaller cube in $\{Q_\ell\}$ of size at most $s/4$ is at least $s/2$.
\end{itemize}

\begin{lemma}
\label{lemma:collar_distance}
Property (iii) holds.
\end{lemma}
\begin{proof}
The cube $Q_\ell$ comes from the collar of $Q'_j$ for some $j$. If it is in the last layer of the collar, then there are no smaller cubes in the cover and the property holds automatically.

If it is not in the last layer, then its $s/2$-neighborhood is covered by the cubes of the collar of size at least $s/2$. If one of these cubes is not in $\{Q_\ell\}$ then it was either replaced by a bigger cube containing it, or it contained no interior points of $Q''_i$ and was removed. In the latter case, observe that none of smaller cubes inside of the removed cube can contain an interior point of $Q''_i$. Therefore, if such a smaller cube were present in a collar of another cube, it would have been removed, too, and  cannot be present in the covering $\{Q_\ell\}$.
\end{proof}

\subsection{Construction of the homotopy $H$.}

Let us now construct the homotopy $H$ with its image in $\bigcup Q_\ell$. By (ii) and because of the inequality (\ref{eq:cover_size}) we will have that
\[
\hc_m(H)\leq \hc_m(\bigcup Q_\ell) < 1.1\cdot 3^m\cdot 6^n\cdot \frac{4^{m+n}}{\epsilon} \cdot \hc_m(X).
\]
So, $H$ will satisfy the restriction on $\hc_m(H)$ in Theorem~\ref{thm:boxing_n}.

The homotopy $H$ will be a concatenation of several homotopies $H^k:X^k\times [0,1]\to \R^n$, where $X^0:=X$ and for $k>0$ the closed set $X^k$ will contain the image of $X$ under $H^{k-1}\circ \ldots\circ H^2 \circ H^1$ and will be fully defined later. (Here and below in this subsection we prefer to imagine that the ambient Banach space $l^n_\infty$ is $\R^n$. The difference between their norms will play no role in our arguments.)

At the $k$th step the set $X^k$ will be covered by a union of a set of dyadic faces $\{Q^k_\ell\}$ with an $(\lceil m\rceil-1)$-dimensional polyhedron $Y^k$. The set of faces $\{Q^k_\ell\}$ is a subset of the faces of the cubes in the original cover $\{Q_\ell\}$. The homotopy $H^k$ will ``push'' $X^k$ out of one of the faces $\{Q^k_\ell\}$, leaving behind a ``trace'' (added to preserve continuity; cf. section 1.4) of dimension at most $(\lceil m\rceil-1)$ which will be later covered by $Y^{k+1}$. After the ``push'' the face will be removed from $\{Q^k_\ell\}$. Its uncovered facets will be added to either $\{Q^{k+1}_\ell\}$ or $Y^{k+1}$. If the facet's dimension is at least $\lceil m\rceil$ (and we will want to push out of it later), it will be added to $\{Q^{k+1}_\ell\}$. Otherwise, it will be added  to $Y^{k+1}$. (``Uncovered" means here that the facet is not a subset of a larger facet
of another cube.)

Formally, at the $k$th step we will have a closed set $X^k$, a set of faces $\{Q^k_\ell\}$, a polyhedron $Y^k$, and a homotopy $H^k:X^k\times [0,1]\to \R^n$. None of the dyadic faces in $\{Q^k_\ell\}$ will be contained in another and their dimensions will be at least $\lceil m\rceil$. We initialize these objects by setting:

\begin{itemize}
\item $X^0:=X$,
\item $\{Q^0_\ell\}:=\{Q_\ell\}$,
\item $Y^0:=\emptyset$.
\end{itemize}

The start $H^k_0$ of the homotopy $H^k$ will be the inclusion $X^k\subset \R^n$. (Here and below we use the notation $H^k_t$ for the map $H^k(t, *)$ from $X^k$ to $\R^n$.) After the $k$th step we will have the following properties:

\begin{itemize}
\item $\dim Y^{k+1}\leq \lceil m\rceil-1$ and $\dim(Y^{k+1}\cap(\bigcup Q^{k+1}_\ell))\leq \lceil m\rceil-2$,
\item $X^{k+1}:=H^k_1(X^k)\cup Y^{k+1}$,
\item $X^{k+1}\subset \bigcup Q^{k+1}_\ell\cup Y^{k+1}$.
\end{itemize}

\begin{lemma}
\label{lemma:projection}
There exist positive constants $C_{scale}$ and $C_{max density}$ depending only on $m$ and $n$, such that for every $m < k \leq n$, cube $Q\subset \R^k$, and $X\subset \inte Q$ with
$\hc_m(X) < C_{max density}\cdot\hc_m(Q)$ there is a point $O\in \inte Q\setminus X$ such that $\hc_m(\phi_O(X)) \le C_{scale}\cdot\hc_m(X)$, where $\phi_O:Q\setminus O\to\partial Q$ is the radial projection in $Q$ with the center $O$.
\end{lemma}
\begin{proof} 
A variant of this lemma is proved in [Gu13], Lemma 7.2. (See also Lemma 2.5 in [Y] for a similar result.) 
The difference between the claim above and Lemma 7.2 in [Gu13] is that Lemma 7.2 is stated for an equilateral simplex $\Delta$ instead of a cube $Q$, and that $m$ ($d$ in the notations of Lemma 7.2) is an integer. The latter assumption is never used in the proof of Lemma 7.2, and, therefore, we can drop it. Let us now explain  why the proof of lemma 7.2 works for cubes instead of simplices. 

Lemma 7.2 guarantees that a suitable center of projection can be found in $\Delta_{1/2}\subset\Delta$ - the scaled by the factor of $1/2$ concentric copy of $\Delta$. 

It is easy to see that the scale factor of $1/2$ in the lemma can be replaced by any positive constant $s<1$. Let us choose $s$ so that $\Delta_s\subset Q\subset \Delta$. Then we can apply Lemma 7.2 directly, find a suitable center of projection $O\in \Delta_s\subset Q$ and then project back from $\partial\Delta$ to $\partial Q$ with the same center $O$. This final projection can increase $\hc_m$ only by a factor bounded by a constant depending on $k$.
\end{proof}

We are ready to describe $H^k$. Pick a smallest dyadic face $Q\in \{Q^k_\ell\}$. We call it the \emph{active} face of the $k$th step.
By our choice of $\epsilon$, we will have that $\hc_m(X^k\cap \inte Q)< C_{max density}\cdot\hc_m(Q)$. We will verify this property after finishing describing the construction. Let $O\in \inte Q$ be the point whose existence is guaranteed by Lemma~\ref{lemma:projection}. Let $\widetilde{H}^k:X^k\times [0,1]\to \R^n$ be the identity outside of $Q$ and the linear interpolation between the identity and the map $\phi_O$ inside.

\begin{figure}[ht]
\center
\includegraphics{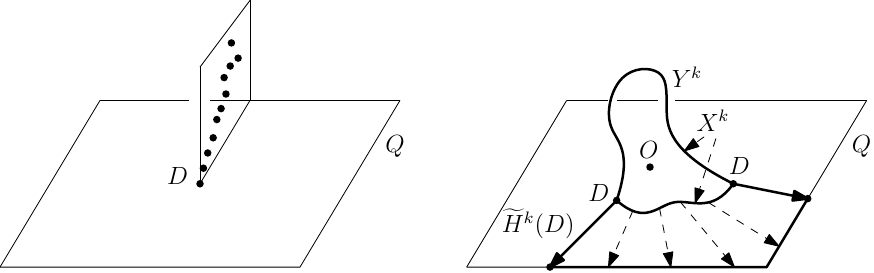}
\caption{Homotopy $H^k$.}
\label{figure:homotopy}
\end{figure}

The homotopy $\widetilde{H}^k$ can be discontinuous at some point in the interior of $Q$. This can only happen if there is a sequence of points of $X^k$ outside of $Q$ converging to it. By Lemma~\ref{lemma:dyadic}, if such sequence of points is covered by a dyadic face, then the face either covers $Q$, or is strictly smaller than $Q$, see Fig.~\ref{figure:homotopy} (left). In both cases such a face cannot belong to $\{Q^{k}_\ell\}$. (Indeed, $Q$ is a smallest face of dimension $\ge \lceil m\rceil$ that still remains.) So, since this sequence belongs to $X^k$ and $X^k\subset \bigcup Q^k_\ell\cup Y^k$, then it is covered by $Y^k$. This means that all the points where $\widetilde{H}^k$ is discontinuous belong to $\inte Q\cap Y^k$.

Let $D$ be the closure of the set of points where $\widetilde{H}^k$ is discontinuous, $D$ is a subpolyhedron of $Y^k$, and it has codimension $\geq 1$ in $Y^k$. Therefore, $\dim D\leq \lceil m\rceil-2$.

After a barycentric subdivision we may assume that $D$ is induced, i.e., that any simplex of $Y^k$ with its vertices in $D$ belongs to $D$. Let $U$ be the simplicial neighborhood of $D$ in $Y^k$, i.e., the union of closed simplices in $Y^k$ which intersect $D$. After several more barycentric subdivisions we can make $U$ arbitrarily small. Any simplex $\sigma\in U$ is a join $\sigma_D*\sigma_Y$ of two of its faces $\sigma_D\in D$ and $\sigma_Y$, $\sigma_Y\cap D=\emptyset$. The complement $\sigma_D*\sigma_Y\setminus \sigma_Y$ retracts to $\sigma_D$ in the standard way: any point of $\sigma_D*\sigma_Y\setminus \sigma_Y$ has the form $a_D\cdot x_D * a_Y\cdot x_Y$, where $0< a_D \leq 1$, $0\leq a_Y < 1$, $x_D\in\sigma_D$, $x_Y\in\sigma_Y$, and is mapped by the retraction to $x_D$. Under this retraction the preimage of every point in $\sigma_D$ is an open cone with its apex being this point and the base being $\sigma_Y$. The retractions agree on intersections of simplices in $U$ and so define a retraction $U\setminus\partial U\to D$, where $\partial U$ is the union of simplices of $U$ disjoint with $D$.

Every simplex in $D$ is a facet of a simplex in $U$, otherwise for its interior points there would not be a sequence in $U\setminus D$ converging to them and these interior points would not be points of discontinuity of $\widetilde{H}^k$.
So, for every $P\in D$ its preimage $U_P$ under the $U\setminus\partial U\to D$ is an open cone with nonempty base, i.e, this cone is not just the point $P$ itself. In each cone $U_P$ we take a smaller cone $U'_P\subset U_P$ with the same apex $P$. Cone  $U'_P$ is obtained from $U_P$ by scaling $U_P$ by the factor $0\leq \frac{\dist(P,\partial Q)}{\diam Q} < 1$, which is $0$ only when $P\in\partial Q$.

\begin{figure}[ht]
\center
\includegraphics{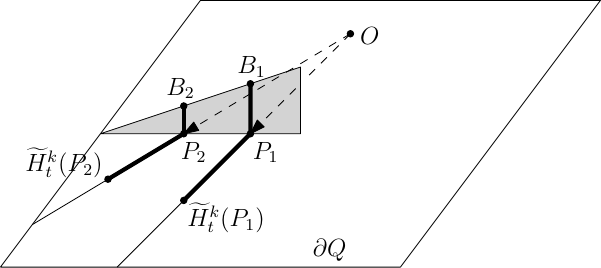}
\caption{Homotopy $H^k$ stretches the line segment $B_1P_1\subset U'_{P_1}$ to the concatenation of $B_1P_1$ with the line segment $P_1\widetilde{H}^k_t(P_1)$. Same for $P_2$.}
\label{figure:homfix}
\end{figure}

Cone $U'_P$ is a union of line segments connecting its apex $P$ with points in its base.
Let us describe $H^k$ on each such segment. Let $B$ be a point in the base of $U'_P$. We have that $B\neq P$, unless $P\in \partial Q$. So, we define $H^k$ on the segment $BP$ as the homotopy fixing $B$ and stretching $BP$ to the concatenation of $BP$ with the line segment $P\widetilde{H}^k_t(P)$, see Fig.~\ref{figure:homfix}. This can be done even if $P\in \partial Q$, because in this case we have $B=P=\widetilde{H}^k_t(P)$.

\begin{figure}[ht]
\center
\includegraphics{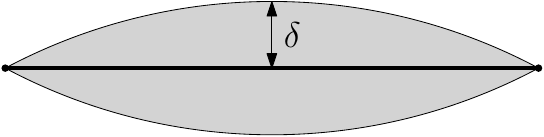}
\caption{A squished $\delta$-neighborhood of a line segment.}
\label{figure:squished}
\end{figure}

Note, that $H^k$ moves only the points in a small \emph{squished $\delta$-neighborhood} of $Q$, which is the union of $\delta\cdot\dist(P,\partial Q)$-neighborhoods of points $P\in Q$, see Fig.~\ref{figure:squished}. The $H^k$-images of these points are also in the squished $\delta$-neighborhood of $Q$. Moreover, we can make $\delta$ arbitrary small by subdividing $D$ several times before constructing $U$. We are going to use this option later, when we will be estimating the distance points move under $H$.

It is clear that $X^{k+1}=H^k_1(X^k)\cup Y^{k+1}$ is covered by $(\bigcup Q^k_\ell\setminus\inte(Q)) \cup Y^k \cup \widetilde{H}^k(D)$.
Let us now describe $\{Q^{k+1}_\ell\}$ and $Y^{k+1}$. We start from
$\{Q^{k+1}_\ell\}=\{Q^{k}_\ell\}$ and $Y^{k+1}=Y^k$ and then modify them. First, we remove $Q$ from $\{Q^{k+1}_\ell\}$.
Then for every facet of $Q$ we do the following. If the facet is covered by some of the faces in $\{Q^{k+1}_\ell\}$ we disregard it.
If it is not covered and its dimension is greater than $(\lceil m\rceil-1)$ (hence the dimension of $Q$ is greater than $\lceil m\rceil$) then we add it to $\{Q^{k+1}_\ell\}$.
If it is not covered and its dimension is exactly $(\lceil m\rceil-1)$ then we add it to $Y^{k+1}$.
Finally, we add $\widetilde{H}^k(D)$ to $Y^{k+1}$.

Let us check that $\{Q^{k+1}_\ell\}$ and $Y^{k+1}$ satisfy the required properties.
We only added facets of $Q$ to $\{Q^{k+1}_\ell\}$ if they were not already contained in a face from $\{Q^{k}_\ell\}\setminus Q$.
The polyhedron $\widetilde{H}^k(D)$ which we added to $Y^{k+1}$ was of dimension at most $(\lceil m\rceil-1)$, and the facets of $Q$ we added to $Y^{k+1}$ only if they were of dimension $(\lceil m\rceil-1)$. The points of $Y^{k+1}\cap(\bigcup Q^{k+1}_\ell)$ which were not already in $Y^{k}\cap(\bigcup Q^{k}_\ell)$ either belong to $\widetilde{H}^k(D) \cap\partial Q$ whose dimension is at most $(\lceil m\rceil-2)$; or they belong to the intersection of an $(\lceil m\rceil-1)$-dimensional facet $F$ of $Q$ (if $Q$ was $\lceil m\rceil$-dimensional) with $\bigcup Q^{k+1}_\ell$.
By Lemma~\ref{lemma:dyadic}, faces in $\{Q^{k+1}_\ell\}$ cannot contain an interior point of $F$ because they neither cover $F$ (otherwise $F$ would not be added to $Y^{k+1}$), nor can they be strictly smaller than $F$ by the choice of $Q$. So, $F\cap(\bigcup Q^{k+1}_\ell)$ belongs to $\partial F$ and has dimension at most $(\lceil m\rceil-2)$.

Clearly, after repeating these homotopies $H^k$ a finite number $k_{max}$ of times (bounded from above by the number of faces of dimension at least $\lceil m\rceil$ in $\{Q_\ell\}$) we will have that $\{Q^{k_{max}}_\ell\}=\emptyset$ and $X^{k_{max}}\subset Y^{k_{max}}$.

\subsection{Choosing $\epsilon$.}
We will now verify that we can choose $\epsilon$ so that at each step $k$ we will have that $\hc_m(X^k\cap \inte Q)< C_{max density}\cdot\hc_m(Q)$ for the active face $Q$ of this step.

First, we can completely disregard the ``insertions" (or ``traces") as $Y$ is $(\lceil m\rceil-1)$-dimensional, and, therefore, $HC_m(Y)=0$. At any step the only other source of new points
in $\inte Q$ for some face $Q$ is the arrival of points from the interior of the active face under the radial projection. Of course, this
can happen only if the boundary of the active face intersects $\inte Q$.

So, let us consider the set of all faces which can be active at some step and from the interior of which a point may arrive to the interior of $Q$ after performing
several steps of our construction.
Formally, a sequence $Q^{k_1}, Q^{k_2}, Q^{k_3}\ldots, Q^{k_i}=Q$ of faces from $\{Q_\ell\}$ is called \emph{admissible (for $Q$)} if:
\begin{itemize}
\item The dimension of $Q^{k_{j}}$ is at least $\lceil m\rceil$  for all $1\leq j\leq i$.
\item The sizes of the faces in the sequence are non-decreasing.
\item Any face from $\{Q_\ell\}$ covering  $Q^{k_{j}}$ is a coface\footnote{$A$ is a coface of $B$ iff $B$ is a face of $A$.} of $Q^{k_{j}}$  for all $1\leq j\leq i$.
\item $\partial Q^{k_j}\cap \inte Q^{k_{j+1}}\neq \emptyset$ for all $1\leq j < i$.
\end{itemize}
The first two conditions are equivalent to $Q^{k_{j}}$ being an active face of some step.
The last condition means that at this step a point of $X$ might project from $\inte Q^{k_{j}}$ to $\inte Q^{k_{j+1}}$.
Denote by $A(Q)\subset \{Q_\ell\}$ the set of all faces which belong to at least one admissible for $Q$ sequence.

\begin{lemma}
\label{lemma:face_sequence}
The number of elements in $A(Q)$ is bounded from above by some constant $C_{adm}$ depending only on $m$ and $n$.
All faces in $A(Q)$ are not larger than $Q$.
\end{lemma}

\begin{proof}
Let $Q^{k_1}, Q^{k_2}, Q^{k_3}\ldots, Q^{k_i}=Q$ be an admissible sequence.
Lemma~\ref{lemma:dyadic} implies that for each $1\leq j < i$ one of the following two options is realized:
\begin{itemize}
\item[(1)] Face  $Q^{k_j}$ covers $Q^{k_{j+1}}$. Then $Q^{k_j}$ is a coface of $Q^{k_{j+1}}$ and $Q^{k_{j+1}}$ is a face of $Q^{k_j}$.
\item[(2)] Face  $Q^{k_{j+1}}$  is strictly larger than $Q^{k_j}$.
\end{itemize}

Without any loss of generality we may assume that the size of $Q^{k_1}$ is $1$.

Suppose first that all faces in the sequence have size $1$. Then option $(2)$ is never realized and all the faces in the sequence are faces of $Q^{k_1}$. Let us call such sequence a \emph{primitive} sequence of $Q^{k_1}$. All the faces in  it are covered by $Q^{k_1}$ which itself is covered by a single dyadic cube of size $1$.

Suppose now the sequence has a face of size at least $2$ and let $Q^{k_j}$ be the first such face.
By Lemma~\ref{lemma:collar_distance}, the size of $Q^{k_j}$ is exactly $2$. Consider a point $x\in\partial Q^{k_{j-1}}\cap \inte Q^{k_j}$. It lies in the cube $Q^{k_{j-1}}$ of size $1$, so by Lemma~\ref{lemma:collar_distance}, the $\ell^\infty$-distance from $x$ to any of the cubes of size at least $4$ is at least $2$. On the other hand, $x\in \inte Q^{k_j}$ and the size of $Q_{k_j}$ is $2$, so the $\ell^\infty$-distance from $x$ to any point of $Q_{k_j}$ is strictly less than $2$. Therefore, $Q^{k_j}$ is disjoint with any of the cubes of size at least $4$.

So, $Q^{k_{j+1}}$ cannot be of size $4$ or more and so has size $2$, same as $Q^{k_j}$. Meaning that we are dealing with option $(1)$ and $Q^{k_{j+1}}$ is a face of $Q^{k_j}$. Then $Q^{k_{j+1}}\subset Q^{k_{j}}$ is also disjoint with any of the cubes of size at least $4$. This implies that $Q^{k_{j+2}}$ likewise has size $2$ and $Q^{k_{j+2}}\subset Q^{k_{j+1}}\subset Q^{k_{j}}$ is also disjoint with any of the cubes of size at least $4$. And so on.
So, the rest of the sequence $Q^{k_j}, Q^{k_{j+1}}, Q^{k_{j+2}}, \ldots, Q^{k_i}=Q$ is just a primitive sequence of $Q^{k_j}$.

We proved that our original admissible sequence  $Q^{k_1}, Q^{k_2}, Q^{k_3}\ldots, Q^{k_i}=Q$ is either primitive or is a concatenation of two primitive sequences. In the first case the size of $Q$ is $1$ and the whole sequence is covered by a dyadic cube of size $1$ of which $Q$ is a face. In the second case the size of $Q$ is $2$, the first primitive sequence is covered by a dyadic cube of size $1$, the second primitive sequence is covered by a dyadic cube of size $2$, these cubes intersect and $Q$ is a face of the bigger one.

In either case, the $\ell^\infty$-distance between any face in the sequence and $Q$ is at most $1.5$ times larger than the size of $Q$.
Also, every face in the sequence is at most $2$ times smaller than $Q$. The number of faces in $\{Q_\ell\}$ satisfying both of these conditions is clearly bounded by some constant depending only on $n$ and $m$.
\end{proof}

For brevity let us introduce notations $h^k(Q):=\hc_m(\inte Q\cap X^{k})$ and $h^k(A(Q))=\sum_{Q_1\in A(Q)}h^k(Q_1)$.
We are interested in the dependence  of  $h^k(A(Q))$ on the step number $k$.

Take some face $Q_1\in A(Q)$. Suppose that at the $k$th step with the active face $Q_a\neq Q_1$ some new points of $X$ arrived to $\inte Q_1$, meaning that there is a point in $\inte Q_1\cap X^{k+1}$ which was not in  $\inte Q_1\cap X^{k}$. At this step we only did a projection from the interior of $Q_a$ to its boundary, so $\partial Q_a\cap \inte Q_1\neq\emptyset$. So, $Q_a$ is in $A(Q)$ as well.

We see that $h^{k+1}(A(Q))$ differs from $h^{k}(A(Q))$ only if the active face $Q_a$ of the $k$th step is in $A(Q)$.
When we project from $Q_a$, only parts of $Y$ remain in $\inte Q_a$ and their Hausdorff content is zero.
(More precisely, $\inte Q_a\cap X^{k+1}=\inte Q_a\cap Y^{k+1}$ and $\hc_m(\inte Q_a\cap Y^{k+1})=0$ because $Y^{k+1}$ is at most $(\lceil m\rceil-1)$-dimensional.) The content which was in $\inte Q_a$, i.e., $h^{k}(Q_a)$, can be bounded from above by $h^{k}(A(Q))$. During the projection it is scaled by the factor of at most $C_{scale}$ and projected to the interior of some of the faces of $Q_a$, some of which can also be in $A(Q)$. The number of such faces is bounded by $3^n-1$, the total number of faces in an $n$-dimensional cube. In total we get that (assuming $C_{scale}>1$ which we can safely do):
\[
h^{k+1}(A(Q)) < h^{k}(A(Q)) + h^{k}(A(Q)) \cdot (3^n-1) \cdot C_{scale} < h^{k}(A(Q)) \cdot 3^n \cdot C_{scale}.
\]
We also have that
\[
h^0(A(Q)) < \epsilon\cdot\hc_m(Q)\cdot |A(Q)| \leq \epsilon\cdot\hc_m(Q)\cdot C_{adm}
\]
because every face has initial density less than $\epsilon$, and by Lemma~\ref{lemma:face_sequence}, every face in $A(Q)$ is not larger than $Q$ and $|A(Q)|\leq C_{adm}$.

Finally, we have that $h^{k}(A(Q))$ changes only at the steps with active face from $A(Q)$, and there are at most $|A(Q)|\leq C_{adm}$ such steps. Assembling it all together we get that
\[
\hc_m(\inte Q\cap X^{k})=:h^k(Q)\leq h^k(A(Q)) <  \epsilon\cdot\hc_m(Q)\cdot C_{adm} \cdot (3^n C_{scale})^{C_{adm}}.
\]
By choosing $\epsilon < \frac{C_{max density}}{C_{adm} \cdot (3^n C_{scale})^{C_{adm}}}$ we can guarantee that $\hc_m(\inte Q\cap X^k)< C_{max density}\cdot\hc_m(Q)$ as required.

 \subsection{Upper bound on the distance points move under $H$.}

Let $x\in X$ be a point. We will now estimate how far it can travel under $H$.

Let $x_i$ be the image of $x$ after the $(i-1)$th step.
Let $k_1, k_2, \ldots$ be the steps for which $x_{k_i+1} \neq x_{k_i}$.
Let $Q^{k_1}, Q^{k_2},  \ldots$ be the corresponding active faces.

Recall that by the construction of $H^k$ both $x_{k_i}$ and $x_{k_i+1}$ have to be in the squished $\delta$-neighborhood of $N_\delta (Q^{k_i})$, otherwise we would have $x_{k_i+1} = x_{k_i}$. Likewise, $x_{k_i+1}=x_{k_i+2}=\ldots=x_{k_{i+1}}$ has to be in $N_\delta (Q^{k_{i+1}})$. Moreover, $x_{k_{i+1}}$ cannot be in the boundary of $Q^{k_{i+1}}$, otherwise it would not be moved by $H^{k_{i+1}}$.

For $\delta$ small enough compared to the sizes of $Q^{k_i}$ and $Q^{k_{i+1}}$, we have that $N_\delta (Q^{k_i})\cap (N_\delta (Q^{k_{i+1}})\setminus \partial Q^{k_{i+1}})=\emptyset$, unless $Q^{k_i}\cap Q^{k_{i+1}}$ contains an interior point of either $Q^{k_i}$ or $Q^{k_{i+1}}$. By Lemma~\ref{lemma:dyadic}, either the sizes of $Q^{k_i}$, $Q^{k_{i+1}}$ are different, or they are the same size but one is covered by another.

In the first case we have that $Q^{k_{i+1}}$ is at least two times bigger than $Q^{k_{i}}$, because active faces cannot get smaller. In the second case we have that $Q^{k_{i+1}}$ is a face of $Q^{k_{i}}$, because if $Q^{k_{i}}$ was a face of $Q^{k_{i+1}}$ then $Q^{k_{i}}$ would be covered by some face in $\{Q^{k_{i}}_\ell\}$.

Now consider the sequence $Q^{k_1}, Q^{k_2},  \ldots$, and observe that each next face is either a face of the previous one or is at least two times larger than the previous one. This implies that the sum of sizes of $Q^{k_1}, Q^{k_2},  \ldots$ and hence the distance between $x$ and $H_1(x)$ is bounded from above by $n\cdot 2d_{max}\leq n\cdot 2(\hc_m(X))^{\frac{1}{m}}$, where $d_{max}$ is the size of the maximal cube in the cover $\{Q_\ell\}$. (Here $n$ majorizes the number of times that a point can be projected within the same closed $n$-dimensional cube.)
\forgotten

    
\forget
\subsetcion{Proof of Theorem 1.6}

The proof of Theorem 3.1 remains as is since in this proof we were considering covering by dyadic cubes.

The derivation of Theorem 1.5 from Theorem 3.1 remains mostly intact when we replace $\HC_m$ by $\HC^d_m$. 

The only insignificant difference is in the proof of the base of the
induction (the case of $m\in (0,1]$. We will still start from a nearly optimal finite covering
of $X$ by dyadic cubes. We choose an arbitrarily small $\epsilon>0$ and a finite collection of dyadic cubes $Q_i$ with radii $r_i$ such that $\HC_m(X)+\epsilon\geq\Sigma_i \HC_m(Q_i)$.

Each connected component of the union of these dyadic is contained in a metric ball (= a cube) of radius that does not exceed $\Sigma_i r_i\leq (\Sigma_i r_i^m)^{1\over m}\leq \HC_m^{1\over m}(X)+\epsilon.$ We map the connected component of $\cup_i Q_i$ to the center of this ball and homotop it along the radii of the ball.

The proof of the induction step is as in the proof of Theorem 1.5. Good balls are not necessarily dyadic cubes, as they are not used to cover anything. The covering used in
the definition of $\bigtilde{\HC_m}$ is a nearly optimal covering of $X$ by dyadic cubes.
Otherwise the proof is the same.
\forgotten

\section{Proof of Theorem 1.3 using Theorem 1.5}

If $B$ is infinite-dimensional, then the homotopy between $M^n$ and its image in a finite-dimensional subspace of $B$ that we constructed in section 1.5 and used to reduce to the case of a finite-dimensional Banach space can be made PL. Also, after a small perturbation in the infinite-dimensional space, its image will become injective for $t\in [0,1)$. We can also make it injective for all $t$ because we can always raise the dimension of the finite-dimensional target space at $t=1$  to at least $2n+3$ by
adding new ``dummy" dimensions and then use them for a small perturbation.
 
 So, it is sufficient to consider the case of a finite-dimensional
 $B$. In the PL- case we are going to find a $2$-bilipschitz approximation of the norm on $B$ with a norm such that the unit ball is a convex polygon. In the smooth case, we find a $2$-bilipschitz approximation by a norm with smooth metric spheres.

 Now we can ensure that the boundaries of all good balls constructed in section 2.2 intersect $X=M^n$ along a PL- (correspondingly, smooth) 
 submanifold. As a result, the image $C^{n+1}$ of the homotopy $H$
 constructed in our proof will be PL (correspondingly,
 a smooth polyhedron). After taking a subdivision, we can assume that the image of the homotopy $H$ is a (linear or smooth) polyhedron embedded in $B$, and its boundary is $M^n$.i

 First, consider the case where the codimension of $M^n$ in $B$ is one. Then $M^n$ is the boundary of a compact manifold with boundary $W^{n+1}\subset B$. Moreover, for each point $p$ in the interior of $W^{n+1}$, $M^n$ does not bound in $B\setminus\{p\}$.
 On the other hand, $M^n$ is the boundary of the image of $H$. 
 Therefore, each point of $W^{n+1}$ is contained in the image of $H$. 
 Therefore, $\HC_m(W^{n+1})\leq \HC_m(H(M^n\times [0,1]))$, and Theorem 1.3(2) follows from Theorem 1.5(2).
 Observe that each point of $W^{n+1}$ will be contained in a trajectory $H(\{x\}\times [0,1])$ of some point $x\in M^n$. Therefore, Theorem 1.3(3) follows from Theorem 1.5(1). 

 From now on, we are going to assume that $M^n$ has a finite codimension greater than one.
 Denote the image of homotopy $H$ by $C^{n+1}$.
 It remains to convert $C^{n+1}$ into a pseudomanifold. 
 Consider 
 a simplex $s^n$ in $C^{n+1}$ of codimension $1$. 
 Assume that $N>2$ of $(n+1)$-dimensional
 simplices of $C^{n+1}$ meet at $s^n$; we need to reduce this number to at most two.
 The general idea is to detach an even number of these $(n+1)$-dimensional simplices
 from the interior of $s^n$, pair them and attach their boundaries to each other, thereby reducing the
 number of cofaces of $S^n$ to $2$, if $s^n$ is not in $M^n$, and to $1$ if $s^n\subset M^n$. 
 To accomplish this plan, we choose a very small non-negative function $\epsilon(x)$ on $s^n$ that vanishes only on the boundary of $s^n$,
 and identify the finite-dimensional ambient Banach space with the Euclidean space of the 
 same dimension by means of a linear isomorphism.
 Let $L$ (correspondingly, $L(x)$)
 denote the orthogonal complement to the plane
 of $s^n$ in the PL case (correspondingly, the orthogonal complement to the tangent space to $s^n$ at $x$ is the smooth case). 
 For each $x\in s^n\setminus\partial s^n$
 consider the sphere $S(x,\epsilon(x))$ centered at $x$ of radius $\epsilon(x)$ in $L+x$, (correspondingly, $L(x)+x$).
 As the codimension of $M^n$ in $B$ is at least $2$, the dimension of this sphere is at least $1$. Suppose a general position of all $(n+1)$-simplices incident to $s^n$ (which can always be achieved by a small perturbation) and a smallness of $\epsilon(x)$. Then $S(x,\epsilon(x))$
 intersects each of these $N$ $(n+1)$-dimensional simplices at a point. We are going to choose an even number of these $(n+1)$-simplices equal to $N-2$ or $N-1$, group
these $N-1$ or $N-2$ points in pairs, connect these pairs by pairwise non-intersecting arcs in $S(x,\epsilon(x))$ and remove the parts of these $N-2$ or $N-1$ simplices inside
$S(x,\epsilon(x))$ in $L+x$ (or $L(x)+x$). The choice of simplices that we disconnect from the interior $s^n$ and their grouping in pairs should be the same for all points $x$, which is clearly possible.

Proceeding in this way, we obtain a polyhedron that we will denote by $\tilde{W}^{n+1}$.
By construction, each $n$-simplex of $\tilde{W}^{n+1}$ is incident to only two $(n+1)$-simplices, if this simplex is not in $M^n$, and to one $(n+1)$-simplex, if it is in 
$M^n$. Now we remove all simplices of dimension $\leq n$ in $\tilde{W}^{n+1}$ that are not faces of $(n+1)$-simplices $\tilde{W}^{n+1}$.  Denote the result by $W_1^{n+1}$. Finally, consider the equivalence classes of $(n+1)$-dimensional simplices in $W_1^{n+1}$, where two simplices are equivalent if they can be connected by a chain of ``flips" over $n$-dimensional simplices to an adjacent $(n+1)$-simplex. Consider the
unions of (closed) $(n+1)$-simplices in each of these equivalence
classes. The boundary of the union of simplices in one of these equivalence classes will be $M^n$, and this union will be a pseudomanifold $W^{n+1}$ that we want to construct. Indeed, the boundary of the union of all simplices in an equivalence class must be either $M^n$ or empty, because each proper subcollection of $n$-simplices in $M^n$ has a nonempty boundary and therefore cannot be a boundary of a $(n+1)$-dimensional chain.

\section{Boxing inequalities on Riemannian manifolds and linearly contractibele geodesic spaces}

\subsection{Introduction}
In this section, we will prove Theorems 1.8 (and, as an almost immediate corollary, Theorems 1.9 and 1.10). 

A tempting way to prove this theorem would be to embed the ambient space $M$ into $L^\infty(M)$ by means of the Kuratowski embedding, to construct the desired filling in $L^\infty(M)$ and then try to map it to $M$ using linear contractibility of $M$. This approach does not work, as the dimension of the image
of $H$ is, in general, greater than the dimension of $M$, and the dimension of $M$ can be uncontrollably greater than $m$. If one tries to map the image of $H$ to $M$ using an induction argument, where
one extends to skeleta of increasing dimension, then the radii of the balls will increase by a factor of $\Lambda$ at each step, eventually leading
to the factor $\Lambda^{dim\ M}$. However, this approach can be adapted to prove an analogue of Theorem 3.1. Meanwhile, we will be just modifying the approach explained in section 2 to deal with linearly
contractible $M$. When $M$ is a manifold, we will just need to verify that the approach works and change the
recurrent relations between $c_1(m)$ and $c_2(m)$ to incorporate a dependence on $\Lambda$. In the situation of Theorem 1.8 B the basic construction used in section 2 needs to be somewhat changed.

Let $M$ be a $(\Lambda,\mu,r_0)$-linearly contractible Riemannian manifold, where $\mu$ or $r_0$ (or both) can be equal to $\infty$.
Let $X$ be a compact subset
of $M$. In the case where $r_0\not=\infty$ we assume that $\HC_m(X)<\epsilon(m,\Lambda)r_0^m$ for a sufficiently small $\epsilon(m,\Lambda)={\epsilon(m)\over\Lambda^{a_3(m)}}$ that will be chosen later. We are going to review the proof of Theorem 1.5, and see what is needed to be changed to prove Theorem 1.8.A. We do not need to reduce to a finite-dimensional case (as was done in section 1.6). 


\subsection{Homotopies and relative homotopies in linearly contractible spaces}

Note that we cannot connect the images of points $x$ under two continuous maps $f_1,f_2$ of a set $A\subset X$ to a metric ball $\beta$ of radius $r$ by a family of minimizing geodesics in $M$ that continuously depends on $x$. But if the radius $r$ of $\beta$ is less than $r_0$ then $\beta$ can be contracted to a point $q$ within the concentric ball $\Lambda\beta$ of radius $\Lambda r$ by a homotopy $Q$, and now for each $a\in A$ we can connect $f_1(a)$ and $f_2(a)$ by a trajectory first connecting $f_1(a)$ and $p_0$ using the trajectory of $f_1(a)$ under $Q$, and 
then connecting $p_0$ and $f_2(a)$ using the trajectory of $f_2(a)$ under $Q$ traversed in the opposite direction. These trajectories will form a homotopy between $f_1$ and $f_2$ in $\Lambda\beta$, and if
$X$ is $(\Lambda,\mu,r_0)$-linearly contractible for $\mu\not=\infty$, the 
lengths of these trajectories will not exceed $2\mu\Lambda r$.

We will also need to be able to connect two maps $f_1,f_2$ from a polyhedral pair $(A,D)$ to a ball $B$ of radius $r$ that coincide on $D$ by means of a homotopy $H$ so that for all $d\in D$ $H(d,t)=f_1(d)=f_2(d)$ for all $t$. The homotopy $H_0$ described in the previous paragraph is not quite the homotopy we want, as the points of $D$ will move during this homotopy. More precisely,
each point $d$ first moves to $x$ along a path in $\Lambda B$ of
length $\leq \mu\Lambda r$, and then returns retracing the same path in the opposite direction. We can deform $H_0$ into another homotopy $H_1$ between $f_1$, $f_2$ which is as desired via intermediate homotopies as follows. Since $(A, D)$ is a polyhedral pair, there is an arbitrarily small open neighborhood $U$ of $D$ in $A$ that retracts to $D$ so that
$U\setminus D$ is homeomorphic to $N\times (0, 1)$ for some $N$,
and the restriction of retraction $\rho$ of $U$ to $D$ to $U\setminus D$ first maps $U\setminus D$ to $N$ and then to $D$.
(Here we assume that the points of $U$ with small second coordinate are close to $D$.) We can identify the closure $\bar U$
of $U$ with $N\times (0,1]\bigcup D$.
Therefore, for each point $c\in U\setminus D$ we can regard its
second coordinate $t(d)\in (0,1)$ as a ``distance" from $D$, and $\rho(d)$ as a ``near" point in $D$. 

Now we can define the homotopies $H_\tau$, $\tau\in [0,1]$, as follows. In the complement of $U$ $H_\tau=H_0$ 
for all $\tau$. If $d\in D$, let $\gamma_d$ denote the path from $f_1(d)=f_2(d)$ to $x$ during the first half of homotopy $H_0$.
In order to define $H_\tau$ on $D$ we just gradually cancel $\gamma_d$ 
by going along shorter and shorter subarcs of $\gamma_d$ before
turning back and following the same arc in the opposite direction. When $\tau=1$ we do not move $d$ at all.
As we now have a notion of closeness to $D$ in $U$, we can define
$H_\tau$ on $U\setminus D$ by interpolating between $H_0$ for
values of $t(d)$ closer to $1$, and the just constructed
$H_\tau$ on $D$ for values of $t$ closer to $0$. Specifically,
for $\tau\in (0,1)$, we will follow the first half of the
trajectory of $f_1(d)$ to a point $H_0(f_1(d), {1\over 2}-(1-t(d))\tau)$, then cross to $H_0(f_2(d), {1\over 2}+(1-t(d))\tau)$,
and then follow it along the trajectory of $(H_0(f_1(d),t)$ to $H_0(f_1(d),1)=f_2(d)$. In other words, we shorten both halves of the trajectory of $f_1(d)$ under $H_0$, and cross over 
from the first half to the second half. As $\tau$ grows from $0$
to $1$, we stop the trajectory closer to its
beginning (and the end). On the other hand, the further $d$ is from $D$, the further are the points where we cut (even for
$\tau=1)$. Once $d$ reaches the boundary of $D$ we do not shorten
the trajectories of $H_0$, and $H_1=H_\tau=H_0$, as we wanted.

In order to make the definition unambiguous, we need to be able to
connect a family of pairs of very close points in $M$ parametrized by a polyhedron by a continuous family of paths in $\Lambda B$ that are also very short if $\mu<\infty$.
Here, given a positive $\delta$, we can postulate
that the points that we connect are $\delta$-close, and, if $\mu$ is finite, we must ensure that the lengths of the paths are at most $\epsilon$. This is trivial in the case where $M$ is a Riemannian manifold that we are now considering. In fact, as long as $\delta$ is less than the injectivity radius of $M$, we can continuously connect these points with minimal geodesics (of length equal to the distance).

If $M$ is a $(\Lambda,\mu, r_0)$-contractible space (as in Theorem 1.8 B), we can still contract the desired continuous family of short curves as follows: Triangulate the polyhedron $P=\bar U\times [0,1]$
that parametrizes the family of pairs of points into very fine simplices so that $D$, $D\times [0,1]$, $N\times \{1\}\subset \bar U$, and $N\times\{1\}\times [0,1]$ are subcomplexes. 
The pairs of points corresponding to $D\times[0,1]\subset P$ that for each $(d,t)$, $d\in D, t\in [0,1]$, consist of two copies of the same point $\gamma_d({t\over 2})$ are connected by the trivial path. The pairs of points 
$(f_1(y), f_2(y))$ for $y\in N\times \{1\}$ are connected by paths
$H_0(x,t),\ t\in [0,1]$. Similarly, pairs of points $H_0(y,{t\over 2})$, $H_0(y,1-{t\over 2})$ for $y\in N\times\{1\},\ t\in [0,1]$ are connected by arcs of $H_0(y,\tau)$ corresponding to the interval of values of $\tau$ between ${t\over 2}$ and $1-{t\over 2}$.
We need to extend a given map of $P\times \{0, 1\}\bigcup (D\bigcup N\times\{1\})\times [0,1]$ to $B\subset M$ to a map of $P\times [0,1]$. Consider a cell subdivision
of $P\times [0,1]$ into prisms $\sigma\times [0,1]$ over simplices $\sigma$ of $P$. We are going to extend to $P\times [0,1]$ by induction with respect to the dimension
of cells in the cell subdivision of $P\times [0,1]$. The extension
to the $1$-simplices corresponds to the shortest paths between their endpoints.
In order to extend to each cell $c$ of the $(i+1)$-dimensional skeleton
after the map is defined on the $i$-skeleton, including the boundary of $c$, we use the linear contractibility of the spheres.
Eventually, we obtain the desired extension to a map $\psi:P\times [0,1]\longrightarrow B$. If $\mu$ is finite, the length of each path $\psi(p_0\times[0,1])$ will be bounded by
$2dim P\mu\delta\Lambda^{dim P}$, and if $\mu<+\infty$ can be made arbitrarily small by choosing a sufficiently small $\delta$. (In $\mu=\infty$, we do not care about the lengths of these paths,
only about the fact that they are contained in $B$, which will similarly hold once $\delta$ is sufficiently small.)

\forget
Here is a different way to construct a homotopy $H_1$
that is constant on $D$. First, we are going to homotop $f_1$
into a map $f'_1$ that coincides with $f_2$ on a very small neighborhood of $D$ and with $f_1$ outside of a somewhat larger
open neighborhood $V$ of $D$. The image of this homotopy will be in 
$(1+\delta)B$ for an arbitrarily small $\delta$. If $\mu<\infty$, then the lengths of the trajectories of this homotopy will have an arbitrarily small length. Finally, this homotopy will be constant on $D$. To construct this homotopy, we need to solve an extension problem $\bar U\times [0,1]$, where the map on $\bar V\times \{0\}$ is $f_1$, the map on $\bar U\times \{1\}$ is $f_2$, the map on $D\times [0,1]$ is identity, and the map on $\bar V\setminus V\times [0,1]$ is $f_1$. We can choose an arbitrarily fine triangulation of the domain, and do this extension using the induction with respect to the skeleta gaining in the process the factor $\Lambda^{dim A}$,
where we cannot control $dim A$. But this factor will be multiplied by the $\max_{x\in\bar V} dist (f_1(x), f_2(x))+\epsilon$, where $\epsilon$ characterizes the sizes of images of simplices under $f'_1$ and $f_2$ and becomes arbitrarily small for sufficiently fine triangulations.
Since the sizes of the neighborhoods $U$ and $V$ can depend 
on $dim\ A$, we can always ensure that the maximal distance above is very small, and its product with $\Lambda^{diam A}$ and, if $\mu\not= \infty$, with $\mu\Lambda^{dim P}$ is very small, too.

Now we construct the homotopy $H_0$ (that does not fix $D$) for
$f_1'$ and $f_2$ as above at the beginning of this section. Since these maps coincide on $U$, the trajectory of each point $u\in U$ under $H_0$ consists of two
arcs connecting $f'_1(u)=f_2(u))$ with $x$ and traversed in the opposite directions. As above, we can identify $U\setminus D$ with $N\times (0,1]$ for a polyhedron $N$, where points of $D$ are
limits of points $(n,t)$, as $t\longrightarrow 0$. We think
about $N\times \{1\}$ as the ``outer" boundary of $U$. Now, we 
continuously alter $H_0$ on $U$ keeping it as is outside of $U$. In other words, we will obtain a family of homotopies $H_\tau$,
$\tau\in [0,1]$, where $H_1$ is the desired homotopy relative to $D$. At each moment of time $\tau$ we will be shortening the trajectories of $H_0$ from $f'_1(u)=f_2(u)$ to the same point as follows: If $u=(n,t), n\in N, t\in [0,1]$, then we will be following the trajectory of $H_0(u,*)$ to $H_0(u,{1\over 2}(1-\tau(1-t)))$, and then will be returning back along the same trajectory. Note that $H_0=H_\tau$ for all $\tau$
and all points $u\in N\times\{1\}$, yet as $t\longrightarrow 0$, trajectories of $(n,t)$ for $\tau=1$ will tend to the constant trajectories, as desired.
\forgotten
Thus, we have analogs of homotopies that were constructed in
section 2 by simply connecting points by straight line segment.
We are going to verify that the argument in section 2 can be carried through (with some modifications).

\subsection{Proof of Theorem 1.8.A}

The reduction of Theorem 1.8.A to its weaker version is very similar to the same reduction for Theorem 1.5 discussed in 
section 2. As in section 2, we choose $C_1(m)=1+\frac{1}{m}$,
$C_2(m)=1+\frac{2}{m}$, so that $\frac{2}{C_2(m)-C_1(m)}=2m$.
But now $A_1$, $c_1$, $c_2$ depend not only on $m$, but also on $\Lambda$ (but not on $\mu$). In section 2 we have chosen $A(m)=30mc_2(m-1)$. Now we observe that the argument in section 2.4 to determine how to choose $A(m)$ involving the inequality (*) remains completely the same, and we can choose $A(m,\Lambda)=30mc_2(m-1,\Lambda)$.
The cones $CK_i$ are defined by fixing a contraction $c_i:K_i\times [0,1]\longrightarrow M$ in a ball of radius $\Lambda\cdot (R_i+c_1(m-1,\Lambda)({2\over C_2(m)-C_1(m)})^{1\over m-1}C_2(m)({r_i\over A(m)^{m\over m-1}})$.

Now the points during the first step do not move by more than
$\Lambda (R_i+c_1(m-1, \Lambda)(\frac{2}{C_2(m)-C_1(m)})^{\frac{1}{m-1}}C_2(m) \frac{1}{A(m)^{\frac{m}{m-1}}}r_i$,
and inequality (2.1) becomes
$$c_1(m, \Lambda)\geq 1+20m\cdot 3^m\Lambda \cdot (90m\ c_2(m-1, \Lambda)+\frac{3c_1(m-1,\Lambda)}{(15\ c_2(m-1,\Lambda))^{\frac{1}{m-1}}}). \ \ \ \ \ \eqno(5.1)$$

Assume that for each $m$ we will choose $c_2(m,\Lambda)$ so that $c_2(m,\Lambda)\geq (3\Lambda\ c_1(m,\Lambda))^{m-1}$. Then we can recursively define $c_1(m,\Lambda)$
as $c_1(m,\Lambda)$ for $m\leq 1$, and for $m>1$ as $c_1(m,\Lambda)=2000m^2\ 3^m\Lambda\ c_2(m-1,\Lambda)$.

However, $c_1(m, \Lambda)\HC_m(X)^{1\over m}$ will no longer control the lengths of the trajectories 
of homotopy $H$. However, if the ambient manifold is $(\Lambda,\mu,r_0)$-linearly contractible for a finite $\mu\geq 1$, this length does not exceed 
\par\noindent
$2\mu\ c_1(m, \Lambda)\HC_m(X)^{1\over m}$.

All this will be true only if the homotopies explained in section 2.3 can be done in our situation. As we explained, this just means that the maps that need to be connected by a homotopy have their images in a metric ball of radius $r_0$ in $X$. 
Looking at the description of our construction in section 2.3, we see these images are naturally contained in metric balls in $X$ and the maximal radii of these balls do not exceed
\begin{multline*}
\Lambda\cdot (R_i+ c_1(m-1,\Lambda)(2m)^{1\over m-1}{1+{2\over m}\over A(m)^{1\over m-1}}{r_i\over A(m)}\leq \\
\leq \Lambda\cdot (90mc_2(m-1,\Lambda)+{3c_1(m-1,\Lambda)\over (15c_2(m-1,\Lambda))^{1\over m-1}})\HC_m(X)^{1\over m}.
\end{multline*}
Using our assumption that for all $m$ $c_2(m,\Lambda)\geq (3\Lambda c_1(m,\Lambda))^{m-1}$, we see that this radius will be less than $r_0$, if $HC_m(X)\leq \epsilon(m,\Lambda)r_0^m$, where $\epsilon(m,\Lambda)={1\over (100m\Lambda c_2(m-1,\Lambda))^m}$.

Assuming $c_2(m,\Lambda)\geq (3\Lambda c_1(m,\Lambda))^{m-1}$, we see that, similarly to inequality (2.2), we can choose
$$c_2(m,\Lambda)=(6000m^2 3^m\Lambda^2 c_2(m-1,\Lambda))^m$$
for $m>1$ (and $c_1(m,\Lambda)=\Lambda^m$ for $m\leq 1$,
Thus, $c_2(m,\Lambda)=\Lambda^{a_2(m)}
\tilde c_2(m)$
for $a_2(m)={2m\over m-1}\cdot (m^{\lceil m\rceil}-1)$, and
$\tilde{c}_2(m)=1$ for $m\leq 1$, and $\tilde{c}_2(m)=(6000m)^m\tilde{c}_2(m-1)^m$ for $m>1$. Returning to our recurrent
relations for $c_1(m,\Lambda)$ we see that $c_1(m,\Lambda)=\Lambda^{a_1(m)}\tilde{c}_1(m)$ for $a_1(m)=1$, when $m\leq 1$, and $a_1(m)=a_2(m-1)+1$ for $m>1$, and for an appropriate $\tilde{c}_1(m)$. Now we see that
$\epsilon(m,\Lambda)={\epsilon(m)\over \Lambda^{a_3(m)}}$ for
$a_3(m)=ma_2(m-1)+1$ and an appropriate $\epsilon(m)$.

It remains to prove an analogue of Theorem 3.1. Here is the assertion that we are going to prove:

\begin{lemma}
    For each positive integer $n$ a positive $m\leq n-1$, a positive real $R>0$,
    a $n$-dimensional Riemannian manifold $M$, 
    and a point $x_0\in M$ there exist positive constants $c_1(M,R, x_0), c_2(M,R, x_0), \epsilon(M,R, x_0)$
    with the following property. Let $X$ be a compact subset of $M$ contained in the open ball $B(x_0, R)$ of radius $R$ centered at $x_0$ such that $\HC_m(X)<\epsilon(M,x_0,R)$.
    Then there exists
    a $m$-filling $(K,\phi, H)$ of $X$ in $M$ with $K\subset M$ such that:
\medskip
\par\noindent
(1) For each $x\in X$ 
$$dist_B(\phi(x), x)\leq c_1(M,x_0, R)
\HC_m^{\frac{1}{m}}(X).$$ 
Moreover, the length of the trajectory $H(\{x\}\times [0,1])$ of $x$ is bounded by 
$c_1(M,R, x_0)\HC_m^{\frac{1}{m}}(X);$
\medskip
\par\noindent
(2) $\HC_m(H(X\times[0,1]))\leq c_2(M,R,x_0)
\HC_m(X)$.
\medskip
\par\noindent
If $\HC_m(X)=0$, then for each $\delta>0$  there exists an $m$-filling of $X$ such that for each $x\in X$ $dist_B(\phi(x),x)\leq \delta$, and $\HC_m(H(X\times [0,1]))\leq\delta.$
\end{lemma}

The estimates in this lemma let us know when we can stop doing inductive steps in the proof discussed above and finish the proof by applying the lemma. We can also ensure
that by the time we apply this lemma, $\HC_m$
of the considered set is sufficiently small to ensure that the right-hand sides
of the inequalities in (1), (2) can be made arbitrarily small. Since $X$ is contained
in an open ball $B(x_0,R)$, we can replace $M$ with its submanifold with boundary $M_R$ that contains a small neighborhood of $M\bigcap \bar B(x_0, 2R)$.

Using Nash theorem embed $M_R$ isometrically in a high-dimensional Euclidean space
$R^a$. There exists a positive $\epsilon=\epsilon(M,R)$ such that for each point $x$ in the $\epsilon$-neighborhood of the embedded $M_R$ in $R^a$ there exists a unique $y\in M_R$ such that $\Vert x-y\Vert= dist(x,M_R)$, and the map $\Phi$ that sends $x$ to this nearest point $y$ in $M_R\subset\mathbb{R}^a$ is $2$-Lipschitz. (Here we could take any other constant $>1$ instead of $2$.)

\forget
Observe that $M_R$ is isometrically embedded in $L^\infty(M_R)$, and for each $\epsilon>0$ $M$ can be similarly $(1+\epsilon)$-bilipschitz embedded in a finite-dimensional subspace $l^{N}_\infty\subset L^\infty(M_R)$, where $N=N(M_R,\epsilon)$ depends on both $M_R$ and $\epsilon$. The idea here is to use distances to all points from a very large but finite subset of points of $M_R$. (See Lemma 2 in [Gu06] for a proof. The main idea here is that if $x,y\in M$ are very close to each other and $z$
is a third point in a very small metric ball centered at $x$ but much farther away from $x$ and $y$ than $dist(x,y)$, and $z$ is near the geodesic $(x,y)$, then 
$dist(x,y)(1-\epsilon)\leq \vert dist(x,z)-dist(y,z)\vert\leq dist(x,y)$.
We can choose here $\epsilon=2$
and define $a=a(M,R)$ 
as $N(M_R,2)$. As $l^{a(M)}_\infty$
is bi-lipschitz equivalent to $\mathbb{R}^a$ for a Lipschitz constant depending on $a(M)$, we can assume that the ambient space is $\mathbb{R}^a$. The inclusion of $M$ in $\mathbb{R}^a$ is not necessarily smooth, as the distance function from a point $x$ is not continuous at $x$. The image of inclusion in $\mathbb{R}^a$ will have singularities in the finite set of points used for the distance functions.
But we can modify the embedding to $L^\infty(M)$ by replacing $dist(x, *)$ for each $x$ by a smooth non-negative function $\tilde d_x$ that differs from $dist(x,*)$ only on a very small neighborhood of $x$ so that for each $y$ $d_x(y)\leq dist(x,y)$. As we use only "distant" $z$ the proof of Lemma 2 implies that for each $\epsilon>0$ choosing
sufficiently many points we can ensure that $(1-\epsilon)dist(x,y)\leq \vert d_z(x)-d_z(y)\vert$ (as $d_z(x)=dist(x,z)$, $d_z(y)=dist(y,z)$). 
But if $z$ is very close to $x$ and $y$ we cannot claim that $\vert d_z(x)-d_z(y)\vert\leq dist(x,y)$. In order to avoid this difficulty, we will make a more specific choice of function $d_x$. For each $\delta>0$ define $\psi_\delta:[0,\infty)$ as ${x^2\over\delta}$ when
$x\in [0,\delta]$, and $\psi_\delta(x)=x$ when $x\geq\delta$. Observe that for all $r,s>0$ $\vert\psi_\delta(r)-\psi_\delta(s)\vert\leq 2\vert r-s\vert$.
This function is not smooth
at $\delta$, but we can find an increasing smooth function $\phi_\delta\leq\psi_\delta$
that coincides with $\psi_\delta$ outside the interval $(0.99\delta,\delta)$ such that $\vert\phi_\delta(r)-\phi_\delta(s)\vert\leq 2\vert r-s\vert$. Now we can define 
$d_x(y)$ as $\phi_\delta(dist(x,y))$ for an appropriately small delta. We see that
the embedding of $M$ into $R^{a}$ defined using $d_z$ for a sufficiently dense finite set
of points $z\in M$ will be $2$-bi-Lipschitz to the embedding defined using functions $dist(z,*)$, and therefore bi-Lipschitz equivalent to the original manifold $M$.
\forgotten

As in the proof of Theorem 3.1 there exists a covering of $M_R$ by
a system of dyadic cubes in the ambient Euclidean space so that this system satisfies properties (i)-(iv) in section 3.3.
The size of maximal cubes is bounded by $const(a)\HC_m(X)^{1\over m}$, and becomes arbitrarily small when $HC_m(X)$ is small. 





Then we can triangulate these cubes in the ambient Euclidean space as in section 3.4. 

The proof of Theorem 3.1 implies the existence of a map $\phi$ of $X$ to the $(\lceil m\rceil -1)$ dimensional skeleton of the considered triangulation and a homotopy $H$ between the identity map of $X$ and $\phi$ with the image in the union of simplices of the triangulation. The lengths of trajectories of $H$ and the $m$-dimensional Hausdorff content of its image satisfy the inequalities in the text of Theorem 3.1. But now we can compose $\phi$ and $H$ with $\Phi$.
As a result, $X$ will be mapped by $\Phi\circ\phi$ to a $(\lceil m\rceil-1)$-dimensional complex in $M_R$, and the homotopy $\Phi\circ H$ will take place inside $M$. But $\Phi$ will not increase the lengths, distances, and $\HC_m$ by much. This completes the proof of Theorem 1.8.

This completes the proof of Theorem 1.8 A.  

\subsection{Proof of Theorem 1.8.B}
Next, we assume that $M$ is an arbitrary $(\Lambda,\mu, r_0)$-linearly contractible metric space, where $\mu$ and/or $r_0$ can be infinite, and $X\subset M$ is a finite-dimensional
polyhedron. First, we will explain how Theorem 1.8 reduces to an analog of Lemma 5.1 where we no longer assume that $M$ is a Riemannian manifold. In this part of the proof no assumptions about $X$ are required.

We observe that we did not need the finite-dimensionality of the ambient Banach space
in section 2. (It became important only in section 3.) So, even if $M$ is infinite-dimensional, we do not attempt to reduce the construction to a finite-dimensional case.

The only difference from the first part of the proof corresponding to section 2 will be in the construction described in section 2.3 and illustrated in Figures 1-4.
Recall that in section 2.3 we described how to modify $X_k$ in the interior of a ball $B(x_i, R_i)$. There we considered a finite collection of closed balls $C_{ki}$ realizing $\widetilde{\HC_m}(\partial B(x_i,R_i)\bigcap X_k)$ and added the union $U_i$ of all balls in this collection to $X_k$. The result was denoted $\tilde X_k$ and we noticed that $\widetilde{\HC_m}(\tilde X_k)=\widetilde{\HC_m}(X_k)$. In addition, $\widetilde{\HC_m}(U_i)=\widetilde{\HC_m}(\partial B(x_i, R_i)\bigcap X_k)\leq {2\over C_2(m)-C_1(m)}C_2(m)^{m-1}{r_i^{m-1}\over A(m)^m}=2m(1+{2\over m})^{m-1}{r_i^{m-1}\over A(m)^m}$.

Now we are going to proceed differently. We observe that for a sufficiently small $\delta_i>0$ $U_i$ will contain the intersection of $X_k$ with the annulus $A(x_i, R_i-\delta_i, R_i)$ centered at $x_i$ and contained between metric spheres with radii $R_i-\delta_i$ and $R_i$. This time we will apply the induction assumption to $U_i$ (instead of $\tilde X_k\bigcap \partial B(x_i,R_i)$). 
We obtain a $(\lceil m\rceil -2)$-dimensional complex $K_i$, 
and maps $\phi:U_i\longrightarrow K_i$ and $j:K_i\longrightarrow B$
as well as a homotopy $H_i$ between the inclusion of $U$ into $M$ and $j\circ\phi$
that satisfies the conditions of Theorem 1.8 for $\lceil m\rceil -1$.

Now, as in section 2.4,  we define $K_{i+1}$ from $K_i$ by removing from $X_k$ its intersections with the interiors of balls $B(x_i, R_i)$, then adding the images of homotopies $H_i$ and cones $CK_i$ over $K_i$ with the vertices at $x_i$.
(Recall that $CK_i$ is defined as the image of a continuous map $c_i$ contracting
$K_i$ to $x_i$ in a metric ball centered at $x_i$ with a radius that is $\Lambda$ times
greater than the metric ball centered at $x_i$ and containing $K_i$.)
Further, we define the map of $X_k$ to $X_{k+1}$ as the identity map outside the interiors of $B(x_i, R_i)$, and as the constant map to $x_i$ in the interiors of $B(x_i, R_i-\delta_i)$. Finally, we define this map $\phi_k$ at each point $x\in X_k$
in the interior of the annulus $A(x_i,R_i-\delta_i, R_i)$ by first calculating
$t(x)=1-{dist_M(x,x_i)-(R_i-\delta_i)\over \delta_i}$, and $\phi_k(x)=H_i(x, 2t(x))$ for
$t(x)\leq {1\over 2}$, and $\phi_k(x)=c_i(\phi(x), 2t(x)-1)\in CK_{i+1}$, if $t(x)\geq {1\over 2}$. The rest of the reduction of Theorem 1.8 to an analogue of Theorem 3.1 goes on as our proof in the case where $M$ is a Riemannian manifold described above. 

Now we will extend Lemma 5.1 (and, as a corollary, Theorem 1.8) to the case where $X$ is a finite-dimensional polyhedron.
In this version of Lemma 5.1 all constants will depend
on $M$ and the dimension of $X$. The dimension of $X$ does not change during the inductive ``cutting off" thin fingers, and therefore, the dependence on $dim\ X$ will be acceptable for our argument.
As before, we assume that $M$ is $(\Lambda,\mu,r_0)$-linearly contractible. 

Consider the Kuratowski embedding of $M$ in $L^\infty(M)$. Theorem 1.5 implies that when $\HC_m(X)$ is small, there exists a small $m$-filling of $X$ in $L^\infty(M)$. The image of $H$ will
be a polyhedron contained in a small neighborhood of $M\subset L^\infty(M)$. Now we combine $H$ with a projection of points in its image to $M\subset L^\infty(M)$. This projection will be built inductively, where the induction is with respect to the dimension of skeleta 
of $X\times [0,1]$. The images of the vertices under $H$ are mapped to (some of) the nearest points in $M$. To map the edges, we use the contraction map of $S^0$ in $M$ formed by the endpoints of an edge.
Then we proceed by induction. At each induction step, we will need to contract spheres in $M$ corresponding to the boundaries of simplices of $X\times [0,1]$. At each induction step, the sizes of the images of simplices of the triangulation of the considered dimension increase by a factor not exceeding a power of $\Lambda$. The procedure ends at the dimension equal to $dim\ X+1$. This completes the proof of Theorem 1.8 B.

Theorems 1.9 and 1.10 would immediately follow from Theorem 1.8.A if not a slightly stronger conclusion in Theorems 1.9, 1.10 with $K$ which is assumed to be a subset of $M$ (and not just being mapped to $M$ by a continuous map). To prove this stronger assertion, observe that if the Riemannian metric on $M$ is analytic, then the map of $K$ to $M$ produced by our proof is subanalytic, and its image is a polyhedron. So, we can consider a very close analytic approximation of
the Riemannian metric on $M$ and use the corresponding $\phi$ and $K$ for this new metric. If the approximation is sufficiently close, these $K$ and $\phi$ will also satisfy the conclusion of a considered theorem for the original metric.

\forget
Now we would like to ensure that all (somewhat curvilinear) simplices formed by the intersection of $M$ with the chosen ``fat" simplices in the ambient Euclidean space are also
``fat", that is, $c(M)$-bilipshitz homeomorphc to the regular simplices of the same dimension for some $c(M)$ that depends only on $M$ and the chosen triangulation of a domain in the Euclidean space containing $B(x_0,R)$.
To do this, we would like to ensure that the intersections of $M$ with each simplex are not too close to its vertices and do not form small angles with each face of this simplex. If we manage to ensure positive lower bounds for angles with simplices and distances to vertices that depend on $M$ and the chosen triangulation of the intersection of $M$ with a ball somewhat larger than $B(x_0,R)$. An arbitrarily small rigid body motion
will ensure that these distances/angles

Proceeding as in [Wh] we can subdivide $M$ into (very small) ``cubes" (or, more precisely,
subsets $1.1$-bilipschitz equivalent to cubes in a tangent space of $M$ at the ``centers" of the cubes) with slowly varying sizes.
Then we can subdivide these ``cubes" into curvilinear ``fat"
simplices of ``slowly varying sizes". As everything is happening
within the ball of radius $R$ centered in $x_0$, we are interested
only in the part of this triangulation in a ball centered at $x_0$
with slightly larger radius than $R$. Fix this triangulation. This will be a finite triangulation, so regardless of what it is
the ``fatness" of simplices, the bilipshitz constants to Euclidean simplices, and the rate of change of size between adjacent simplices will be controlled by some constants that ultimately depend on on $M$. Now, if $X$ has a very small $HC_m(X)$
bounded by an appropriate constant depending on $M$ (including the choice of the above triangulation), we can do the Federer-Fleming projections as in section 3 until we reach that $(\lceil m\rceil-1)$-skeleton of the chosen triangulation (of a neighborhood of the metric ball of radius $R$ centered at $x_0$ in $M$). This completes the proof of the theorem.

Lemma 3.1 now says that $X$ can be homotoped to a subset of a $(\lceil m\rceil-1)$-dimensional polyhedron $K\subset l^{a(M)}_\infty$ 
by means of a homotopy $H$ with the desired properties. Moreover, the same is true for a larger subset $Y$ containing $X$ defined as the union $Y$ of a collection of metric balls in $M$ such that $\HC_m(Y)$ is only slightly larger than $\HC_m(X)$. Without any loss of generality we can assume that the Riemannian metric on $M$ is analytic (as this can be achieved by an arbitrarily close approximation). Then the distance function on $M$ is subanalytic,
and all metric balls and $Y$ can be smoothly triangulated. Then the image of $H$ will also be a polyhedron in 
of dimension equal to dim $M+1$. When $\HC_m(X)$ and $\HC_m(Y)$ are sufficiently small, the image of $H$ is sufficiently close to $M$ in the ambient Banach space.
We can now modify $H$ into a map to $M$. This new map, $H_{\rm new}$,  will be defined by induction with respect to the dimension of skeleta of $Y\times [0,1]$. For each point $v$ in the $0$-dimensional skeleton of a very fine triangulation of $Y\times [0,1]$ we map $v$ to a nearest point in $M$. Observe that any two vertices of the triangulation of $Y\times [0,1]$ connected by an edge are mapped into $\delta$-close points of $l^{a(M)}_\infty$ for a positive $\delta$ that can be made arbitrarily small by choosing a very fine triangulation of
$Y\times [0,1]$. These vertices will be mapped to a pair of points of $M$
at the distance in $l^{a(M)}_\infty$ that does not exceed $2dist(Im(H), M)+\delta\leq 2c_1(m)\HC_m(Y)^{1\over m}+\delta\leq
2\HC_m(X)^{1\over m}+\epsilon+\delta$ for an arbitrarily small $\delta$. But the distance between two points in the metric of $M$
is at most $2$ times greater than the distance between them
in the ambient space $l^{a(M)}_\infty$. So we obtained a map of $0$-dimensional sphere $S^0$ into $M$ such that its diameter
is only slightly greater than $4\HC_m(X)^{1\over m}$. If this distance is less than $r_0$ we can contract this sphere, and the result of contraction will be a map of the $1$-dimensional disc 
into a ball of radius somewhat larger than $4\Lambda\HC_m(X)^{1\over m}$. These $1$-discs will provide the restriction of the $H_{\rm new}$ to the $1$-skelton of the triangulation of $Y\times [0,1]$. Now we need to extend $H_{\rm new}$ to the $2$-skeleton of the triangulation. In order to do that, we need to contract the newly constructed triangles in $M$ corresponding to the boundaries of $2$-simplices. We observe
that the diameters of each of these triangles are at most $8\Lambda \HC_m(X)^{1\over m}$, so we can contract it by a homotopy in $M$ with an
increase in diameter by the factor of (at most) $\Lambda$, etc. After $n=dim M$ steps we will be done, provided
that $\HC_m(X)\leq {const(n)\over r_0^n}$, and the map image will be in a metric ball of radius at most $const(n)\Lambda^{n+1}\HC_m(X)^{1\over m}$.
This completes the proof of part (1) of Lemma 5.1.

Now we are going to prove part (2). Assume that the image of $H$ was covered by an almost optimal collection of balls $\beta_i\subset l^{a(M)}_\infty$ with radii $r_i$ such that $\Sigma_i r_i^n\leq \HC_m(Y)+\delta$.
We would like to map these balls into metric balls
$\beta'_i\subset M$ with radii $r'_i\leq \Lambda^{q(M)}r_i$
so that the image of the map $\bigcup_i\beta_i\longrightarrow\bigcup_i\beta'_i$ would cover
the image of $H_{\rm new}$.
We can assume that the covering $\beta_i$ had been chosen before we triangulated the domain of $H$, so that the images of all simplices of the triangulation are very small compared to
$\min_i r_i$. 

For simplicity we can assume that $H$ does not map $Y\times [0,1]$
but a finely triangulated mapping cylinder $C_\phi$ of $\phi:Y\longrightarrow K\subset l^{a(M)}_\infty$, $a(M)>2dim M+3$,
$H$ is affine on each simplex of the triangulation and was perturbed to avoid intersections between different open simplices of the triangulation. In addition, we assume that $H_{\rm new}$ is a map from $C_\phi$ to $M$ constructed as above. We can triangulate
each $\beta_i$ to make this triangulation compatible with the triangulation of the image of $H$ (refining some simplices of the triangulation of $C_f$, if necessary). Now we can construct
for each $i$ an extension of the restriction of $H_{\rm new}$
to $\beta_i$ constructed as above. Of course, since the dimension of 
simplices of $\beta_i$ can be as high as $a(M)$, we are going
to obtain balls of
\forgotten


\forget
Each dyadic cube $Q$ is contained in the unique dyadic cube $D(Q)$ with size equal to the twice the size of $Q$. Replace each $Q'_i$ by $D(Q_i)$, and remove all $D(Q_i)$ that are properly contained in other cubes of the sequence, as well as extra copies
of the same cube. Denote the cubes in the resulting sequence as $Q_i$.

\begin{lemma}
\par\noindent
(1) For each $i$ $\HC_m(X\cap Q_i)< \HC_m(Q_i)=rad (Q_i)^m.$ 
\par\noindent
(2) Each connected component of $X\cap Q_i$ contains at most one vertex $v$ of $Q_i$.
If a connected component $C$ of $X\cap Q_i$ contains a vertex $v$ and another point $y$, then $v$ is a vertex if the face of $Q_i$ of minimal possible dimension that contains $y$.
\par\noindent
(3) If a connected component $C$ of $X\cap Q_i$ does not contain any vertices of $Q_i$, then there exists unique vertex $v$ of $Q_i$ such that for each $y\in C$ $v$ is a vertex 
of the face of $Q_i$ of minimal dimension that contains $y$.
\end{lemma}

\begin{proof}
\end{proof}

The lemma implies that for each $i$ there exists a map $\phi_i$ that sends each connected component $C$ of $X\cap Q_i$
to the unique vertex $v$ of $Q_i$ such that for each point $y\in C$ $v$ is a vertex of the face of $Q_i$ of minimal possible dimension that contains $y$. Let us 
reorder $Q_i$ according to their sizes starting from the smallest, and then start applying
$\phi_1$, then $\phi_2$, etc. There is an obvious homotopy $H_i$ within $Q_i$ that connects the identity map of $Q_i$ with $\phi_i$: Each point $y$ moves along the straight
line segment from $y$ to $\phi_i$. This time we do not encounter any problem with continuity due to balls having unequal sizes. When a points moves more than one time,
it will move in cubes with sizes that increase by at least the factor of $2$ each time.
Therefore, the sum of lengths of trajectories of all homotopies applicable to $y$ and its images under composition of $\phi_i$ is bounded by twice the size of the largest
cube, where it moves. The image of $H$ defined as the concatenation of $H_i$ is contained in $\bigcup Q_i$, and will not exceed $\Sigma_i\HC_m(Q_i)\leq 2^m\Sigma_i \HC_m(Q'_i)\leq 2^m(\HC_m(X)+\epsilon).$ This completes the proof of the base of 
induction.
\forgotten
{\bf Acknowledgements.} Sergey Avvakumov was partially supported by NSERC Discovery grant and the grant 765/19 from the Israel Science Foundation (ISF). The research of Alexander Nabutovsky was partially supported
by his NSERC Discovery Grant and Simons Fellowship. This research was partially done during the visit of A.N. to SLMath in Fall 2024,
and Hausdorff Institute for Mathematics in Spring 2025.

\bigskip
\address{S.A.: School of Mathematical Sciences, Tel Aviv University, Tel Aviv 69978, Israel;
s.avvakumov@gmail.com}
\par\noindent
\address{A.N.: Department of Mathematics, Bahen Centre, 40 St. George st., Rm 6290, Toronto, Ontario, M5E2S4, Canada; alex@math.toronto.edu}
\end{document}